\documentclass[a4paper,12pt]{amsart}
\usepackage[frenchb,english]{babel}
\usepackage[latin1]{inputenc}
\usepackage{pdfsync}
\usepackage[T1]{fontenc}   
\usepackage{amsmath}
\usepackage{amssymb}
\usepackage{amsxtra}
\usepackage{amscd}
\usepackage{amsthm}
\usepackage{amsfonts}
\usepackage{eucal}
\usepackage[all]{xy}
\usepackage{graphicx}
\usepackage{comment}
\usepackage{epsfig}
\usepackage{psfrag}
\usepackage{mathrsfs}
\usepackage{amscd}
\usepackage{rotating}
\usepackage{lscape}
\usepackage{amsbsy}
\usepackage{verbatim}
\usepackage{moreverb}
\usepackage{url}
\usepackage{color}
\usepackage{xcolor}

\newcommand{\nc}{\newcommand}
\nc{\renc}{\renewcommand}

\nc\restr[2]{{ 
  \left.\kern-\nulldelimiterspace    #1  
  \vphantom{\big|}  
  \right|_{#2}  
  }}

\newtheorem{cor}[subsubsection]{Corollary}
\newtheorem{lem}[subsubsection]{Lemma}
\newtheorem{prop}[subsubsection]{Proposition}

\newtheorem{construction}[subsubsection]{Construction}
\newtheorem{nota}[subsubsection]{Notation}

\newtheorem{thm}[subsubsection]{Theorem}

\renc{\sec}{\section}
\nc{\ssec}{\subsection}
\nc{\sssec}{\subsubsection}

\theoremstyle{definition}
\newtheorem{defi}[subsubsection]{Definition}
\newtheorem{example}[subsubsection]{Example}

\theoremstyle{remark}
\newtheorem{rem}[subsubsection]{Remark}

\numberwithin{equation}{section}

\nc{\on}{\operatorname}

\nc\wt{\widetilde}
\nc\wh{\widehat}
\nc\ol{\ov}
\nc{\oc}[1]{{\overset{\circ}{#1}}}
\nc{\ov}[1]{{\overline{#1}}}
\nc{\isor}[1]{{\xrightarrow[\raisebox{0.25 em}{\smash{\ensuremath{\sim}}}]{#1}}}
\nc{\modmod}{/ \! \! /}

\nc{\mc}{\mathcal}
\nc{\mf}{\mathfrak}
\nc{\mr}{\mathrm}
\nc{\mb}{\mathbb}
\nc{\mbf}{\mathbf}
\nc{\ms}{\mathscr}

\nc{\R}{{\mathbb R}}
\nc{\Z}{{\mathbb Z}}
\nc{\N}{{\mathbb N}}
\nc{\C}{{\mathbb C}}
\nc{\Q}{{\mathbb Q}}

\nc{\Fq}{{\mathbb F}_q}
\nc{\Fl}{{\mathbb F}_\ell}
\nc{\Fqbar}{\ol{{\mathbb F}_q}}
\nc{\Flbar}{\ol{{\mathbb F}_\ell}}
\nc{\Zl}{{\mathbb Z}_\ell}
\nc{\Zlbar}{\ol{{\mathbb Z}_\ell}}
\nc{\Ql}{E}
\nc{\Qlbar}{\ol{{\mathbb Q}_\ell}}
\nc{\hl}{\overset{\leftarrow}h{}}
\nc{\hr}{\overset{\rightarrow}h{}}
\nc{\Gr}{{\on{Gr}}}
\nc{\Hecke}{\on{Hecke}}
 \nc{\Hom}{\on{Hom}}
 \nc{\Coker}{\on{Coker}}
 \nc{\Ker}{\on{Ker}}
 \nc{\Lie}{\on{Lie}}
\nc{\Loc}{\on{Loc}}
\nc{\Pic}{\on{Pic}}
\nc{\Bun}{\on{Bun}}
\nc{\IC}{\on{IC}}
\nc{\Aut}{\on{Aut}}
\nc{\Perv}{\on{Perv}}
\nc{\pos}{{\on{pos}}}
\nc{\Sym}{\on{Sym}}
\nc{\Oe}{\mc O_E}
\nc{\Hcusp}{\big( \restr{ \mc H_{G, N, I, W}^{j, \, \Oe } }{ \ov{x} } \big)^{\on{cusp}}}

\nc{\ta} {{}^\tau}
\nc {\tu}[1]{{}^{\tau^{#1}}\!}

\nc{\Id}{\on{Id}}
\nc{\Fil}{\on{Fil}}
\nc{\pr}{\on{pr}}
\nc{\Res}{\on{Res}}
\nc{\cusp}{\on{cusp}}
\nc{\Frob}{\on{Frob}}
\nc{\diag}{\Delta}
\nc{\gr}{\on{gr}}
\nc{\Inj}{\on{Inj}}
\nc{\Bl}{\on{Bl}}
\nc{\dem}{\noindent {\bf Proof. }}
\nc{\cqfd}{{\ }\hfill $\square$ \vskip 1mm}
\nc{\s}[1]{\langle #1 \rangle}
\nc{\Cht}{\on{Cht}}
\nc{\isom}{\overset {\thicksim}{\to}}
\nc{\sm}{\smallsetminus}

\begin{document}

\title{Cohomology with integral coefficients of stacks of shtukas}

\author{Cong Xue}
\address{Cong Xue: Université Paris Cité and Sorbonne Université, CNRS, IMJ-PRG, 75005 Paris, France}
\email{cong.xue@imj-prg.fr}
 
\maketitle

\begin{abstract}
We construct the cohomology groups with compact support of stacks of shtukas with $\Zl$-coefficients.
We construct the cuspidal cohomology groups and prove that they are $\Zl$-modules of finite type. We prove that the cohomology groups are modules of finite type over a Hecke algebra with $\Zl$-coefficients. 

As an application, we prove that the cuspidal cohomology groups with $\mb Q_{\ell}$-coefficients are equal to the Hecke-finite cohomology groups with $\mb Q_{\ell}$-coefficients defined by V. Lafforgue.

We also state the Eichler-Shimura relations for cohomology groups with $\Zl$-coefficients and prove the compatibility of the excursion operators and the constant term morphisms.
\end{abstract}

\section*{Introduction}

Let $X$ be a smooth projective geometrically connected curve over a finite field $\Fq$. We denote by $F$ the function field of $X$, $\mb A$ the ring of adèles of $F$ and $\mb O$ the ring of integral adèles. Let $N \subset X$ be a finite subscheme. We denote by $\mc O_N$ the ring of functions on $N$ and write $K_{N}:=\Ker(G(\mb O) \rightarrow G(\mc O_N) )$. 

Let $G$ be a connected split reductive group over $\Fq$.

Let $\ell$ be a prime number not dividing $q$. Let $E$ be a finite extension of $\mb Q_{\ell}$ containing a square root of $q$. We denote by $\Oe$ the ring of integers of $E$.
Let $\wh G$ be the Langlands dual group of $G$ over $\Oe$. We denote by $\wh G_E$ the base change of $\wh G$ over $E$.

Let $I$ be a finite set and $W$ be a finite dimensional $E$-linear representation of $\wh G_E^I$.
Varshavsky (\cite{var}) and V. Lafforgue (\cite{vincent}) defined the stack classifying $G$-shtukas $\Cht_{G, N, I, W}$ over $(X \sm N)^I$ and its degree $j$ cohomology group with compact support in degree $j \in \Z$ with $E$-coefficients: $H_{G, N, I, W}^{j, \, E}$. It is an inductive limit of finite dimensional $E$-vector spaces. 
It is equipped with an action of the Hecke algebra $C_c(K_{N} \backslash G(\mb A) / K_{N}, E)$ via Hecke correspondences.

In \cite{cusp-coho}, for every proper parabolic subgroup $P$ of $G$ with Levi quotient $M$, we defined the degree $j$ cohomology group with compact support with $E$-coefficients of the stack classifying $M$-shtukas (denoted by $H_{M, N, I, W}^{' \, j, \, E}$) and 
we constructed the constant term morphism 
$$C_G^{P, \, E}: H_{G, N, I, W}^{j, \, E} \rightarrow H_{M, N, I, W}^{' \, j, \, E} .$$ 
Note that contrary to $H_{G, N, I, W}^{j, \, E}$ which has finitely many "components", $H_{M, N, I, W}^{' \, j, \, E}$ has infinitely many "components", due to the difference between the center of $M$ and the center of $G$.
We defined the cuspidal cohomology group $H_{G, N, I, W}^{j, \, \on{cusp}, \, E}$ as the intersection of the kernels of $C_G^{P, \, E}$ for all proper parabolic subgroups $P$ of $G$. We proved that $H_{G, N, I, W}^{j, \, \on{cusp}, \, E}$ has finite dimension. We also proved that $H_{G, N, I, W}^{j, \, \on{cusp}, \, E}$ contains the Hecke-finite cohomology group $H_{G, N, I, W}^{j, \, \on{Hf}, \, E}$ defined by V. Lafforgue in \cite{vincent} Section 8.
Besides, in \cite{coho-filt}, we proved that $H_{G, N, I, W}^{j, \, E}$ is of finite type as module over a local Hecke algebra $C_c(G(\mc O_u) \backslash G(F_u) / G(\mc O_u), E)$, where $u$ is any place in $X \sm N$, $\mc O_u$ is the complete local ring at $u$ and $F_u$ is its field of fractions. 

However, to prove that $H_{G, N, I, W}^{j, \, \on{cusp}, \, E}$ is in fact equal to $H_{G, N, I, W}^{j, \, \on{Hf}, \, E}$, we need to construct an $\Oe$-lattice in $H_{G, N, I, W}^{j, \, \on{cusp}, \, E}$ invariant under the action of the Hecke algebra $C_c(K_{N} \backslash G(\mb A) / K_{N}, \Oe)$. This is the initial motivation of this paper.

\quad

In this paper, we show that the results in \cite{cusp-coho} and \cite{coho-filt} still hold for cohomology groups with $\Oe$-coefficients. As an application, we prove that $H_{G, N, I, W}^{j, \, \on{cusp}, \, E}$ is equal to $H_{G, N, I, W}^{j, \, \on{Hf}, \, E}$. However, the cohomology groups with $\Oe$-coefficients have more applications. In the following we give the main results of each section.

Let $W$ be a finite type $\Oe$-linear representation of $\wh G^I$. In Section 1, using the geometric Satake equivalence with $\Oe$-coefficients, 
we define the cohomology group with compact support in degree $j \in \Z$ with $\Oe$-coefficients: $H_{G, N, I, W}^{j, \, \Oe}$. It is an inductive limit of finite type $\Oe$-modules. (Note that this cohomology group modulo torsion was already introduced in \cite{vincent} Section 13.)

In Section 2, we define $H_{M, N, I, W}^{' \, j, \, \Oe} $ and construct the constant term morphism $$C_G^{P, \, \Oe}: H_{G, N, I, W}^{j, \, \Oe} \rightarrow H_{M, N, I, W}^{' \, j, \, \Oe} .$$
Then we define the cuspidal cohomology group $H_{G, N, I, W}^{j, \, \on{cusp}, \, \Oe}$ as the intersection of the kernels of $C_G^{P, \, \Oe}$ for all proper parabolic subgroups $P$ of $G$. In Section 3, we prove 
\begin{thm} (Theorem \ref{thm-cusp-dim-fini-second}) (see \cite{cusp-coho} Theorem 0.0.1 for the $E$-coefficients version) 
$H_{G, N, I, W}^{j, \, \on{cusp}, \, \Oe}$ is an $\Oe$-module of finite type.
\end{thm}

We also prove that for any dominant coweight $\lambda$ of $G$ large enough, the morphism from the Harder-Narasimhan truncation $H_{G, N, I, W}^{j, \, \leq \lambda, \, \Oe} $ to the inductive limit $ H_{G, N, I, W}^{j, \, \Oe}$ is injective (Proposition \ref{prop-finiteness-a-b-c}).

Next, in Section 4 we prove 
\begin{prop} (Proposition \ref{prop-TC-commute-with-Hecke-with-niveau}) (see \cite{cusp-coho} Lemma 6.2.6 for the $E$-coefficients version) 
The constant term morphism $C_G^{P, \, \Oe}$ commutes with the action of the Hecke algebra $C_c(K_{N} \backslash G(\mb A) / K_{N}, \Oe)$.
\end{prop}

The arguments in Sections 1 and 2 of \cite{coho-filt} still hold. In Section 5, we prove
\begin{thm}  \label{thm-H-G-N-I-W-Oe-type-fini-as-Hecke-algebra}
(Theorem \ref{coho-cht-Oe-is-Hecke-mod-type-fini-en-u})  (see \cite{coho-filt} Theorem 0.0.2 for the $E$-coefficients version)
$H_{G, N, I, W}^{j, \, \Oe}$ is of finite type as module over the Hecke algebra $C_c(G(\mc O_u) \backslash G(F_u) / G(\mc O_u), \Oe)$, for any place $u$ in $X \sm N$.
\end{thm}
In particular, it implies that the order of torsion in $H_{G, N, I, W}^{j, \, \Oe}$ is bounded (Corollary \ref{cor-torsion-in-coho-bounded}).

In Section 6, we show that $$H_{G, N, I, W}^{j, \, \Oe} \otimes_{\Oe} E \isom H_{G, N, I, W \otimes_{\Oe} E }^{j, \, E} $$ and that the morphism $$H_{M, N, I, W}^{' \, j, \, \Oe} \otimes_{\Oe} E \rightarrow H_{M, N, I, W \otimes_{\Oe} E }^{' \, j, \, E} $$ is injective (note that this is not evident because $H_{M, N, I, W}^{' \, j, \, \Oe}$ may have infinitely many "components"). We deduce that (Proposition \ref{prop-H-cusp-Zl-H-cusp-Ql})
$$H_{G, N, I, W}^{j, \, \on{cusp}, \, \Oe} \otimes_{\Oe} E \isom H_{G, N, I, W\otimes_{\Oe} E }^{j, \, \on{cusp}, \, E} .$$ 
As a consequence, the image of $H_{G, N, I, W}^{j, \, \on{cusp}, \, \Oe}$ in $H_{G, N, I, W\otimes_{\Oe} E   }^{j, \, \on{cusp}, \, E}$ is an $\Oe$-lattice in $H_{G, N, I, W \otimes_{\Oe} E  }^{j, \, \on{cusp}, \, E}$ invariant under the action of the Hecke algebra $C_c(K_{N} \backslash G(\mb A) / K_{N}, \Oe)$.
Then we prove that the inclusion $H_{G, N, I, W \otimes_{\Oe} E  }^{j, \, \on{cusp}, \, E} \supset H_{G, N, I, W \otimes_{\Oe} E   }^{j, \, \on{Hf}, \, E}$ in \cite{cusp-coho} is an equality:
\begin{prop} (Proposition \ref{prop-H-cusp-equal-H-Hf}, conjectured by V. Lafforgue)
The two sub-$E$-spaces $H_{G, N, I, W \otimes_{\Oe} E  }^{j, \, \on{cusp}, \, E}$ and $H_{G, N, I, W \otimes_{\Oe} E  }^{j, \, \on{Hf}, \, E}$ of $H_{G, N, I, W \otimes_{\Oe} E  }^{j, \, E}$ are equal.
\end{prop}

\quad

Section 7 is logically independent of the previous sections. 
In Section 7, we construct the partial Frobenius morphisms on $H_{G, N, I, W}^{j, \, \Oe}$ and prove the Eichler-Shimura relations. 

Section 8 is an application of Section 5.
In Section 8, we use Drinfeld's lemma to 
obtain an action of $\on{Weil}(\ov F / F)^I$ on $H_{G, N, I, W}^{j, \, \Oe}$, where $\ov F$ is an algebraic closure of $F$. Then we construct the excursion operators.

\quad

%
%

\subsection*{Notations and conventions}

\sssec{} \label{subsection-def-Z-G}
Let $G^{\mr{der}}$ be the derived group of $G$ and $G^{\mr{ab}}:=G/G^{\mr{der}}$ the abelianization of $G$.
Let $Z_G$ be the center of $G$ and $G^{\mr{ad}}$ the adjoint group of $G$ (equal to $G / Z_G$).

\sssec{}  \label{subsection-def-Xi}
We fix a discrete subgroup $\Xi_G$ of $Z_G(\mb A)$ such that $\Xi_G \cap Z_G(\mb O) Z_G(F) = \{ 1 \}$, the quotient $Z_G(F) \backslash Z_G(\mb A) / Z_G(\mb O) \Xi_G $ is finite and the composition of morphisms $\Xi_G \hookrightarrow Z_G(\mb A) \hookrightarrow G(\mb A) \twoheadrightarrow G^{\mr{ab}}(\mb A)$ is injective. Note that the volume of $G(F) \backslash G(\mb A) / G(\mb O) \Xi_G $ is finite. We write $\Xi:=\Xi_G$.

\sssec{} We fix a Borel subgroup $B \subset G$. 
By a parabolic subgroup we will mean a standard parabolic subgroup (i.e. a parabolic subgroup containing $B$), unless explicitly stated otherwise.

\sssec{} \label{subsection-def-Lambda-G}
Let $H$ be a connected split reductive group over $\Fq$ with a fixed Borel subgroup.
Let $\Lambda_H$ (resp. $\wh{\Lambda}_H$) denote the weight (resp. coweight) lattice of $H$. Let $\langle \ , \ \rangle: \wh \Lambda_H \times \Lambda_H \rightarrow \Z$ denote the natural pairing between the two.

Let $\wh \Lambda_H^+ \subset \wh \Lambda_H$ denote the monoid of dominant coweights and $\wh \Lambda_H^{\mr{pos}} \subset \wh \Lambda_H$ the monoid generated by positive simple coroots. Let $\wh \Lambda_H^{\Q}:=\wh \Lambda_H \underset{\Z} \otimes \Q$. Let $\wh \Lambda_H^{\mr{pos}, \Q}$ and $\wh \Lambda_H^{+, \Q}$ denote the rational cones of $\wh \Lambda_H^{\mr{pos}}$ and $\wh \Lambda_H^{+}$. 
We use analogous notation for the weight lattice. 

We use the partial order on $\wh{\Lambda}_{H}^{\Q}$ defined by $\mu_1 \leq^H \mu_2 \Leftrightarrow \mu_2 - \mu_1 \in \wh \Lambda_H^{\mr{pos}, \Q}$
(i.e. $\mu_2 - \mu_1$ is a linear combination of simple coroots of $H$ with coefficients in $\Q_{\geq 0}$).

We will apply these notations to $H = G$, $H=G^{\mr{ad}}$ or $H=$ some Levi quotient $M$ of $G$.

\sssec{} \label{subsection-def-Gamma-G-Gamma-M}
We denote by $\Gamma_G$ the set of vertices of the Dynkin diagram of $G$.
Parabolic subgroups in $G$ are in bijection with subsets of $\Gamma_G$. For a parabolic subgroup $P$ with Levi quotient $M$, we let $\Gamma_M \subset \Gamma_G$ denote the corresponding subset; it identifies with the set of vertices of the Dynkin diagram of $M$.

\sssec{} \label{def-prechat-chat-chat-alg}
We use Definition 3.1 and Definition 4.1 in \cite{chat-alg} for prestacks, stacks and algebraic stacks.

\sssec{}   \label{def-D-c-b-pour-chat-alg}
As in \cite{chat-alg} Section 18, \cite{LO08b} and \cite{LO09}, for $\mathcal{X}$ an algebraic stack locally of finite type over $\Fq$, we denote by $D_c(\mathcal{X}, \Oe)$ (resp. $D_c^{(b)}(\mathcal{X}, \Oe)$, $D_c^{(+)}(\mathcal{X}, \Oe)$, $D_c^{(-)}(\mathcal{X}, \Oe)$) the unbounded (resp. locally bounded, locally bounded below, locally bounded above) derived category of "constructible $\Oe$-modules on $\mathcal{X}$". We have the six operations and the notion of perverse sheaves. 

For a finite type morphism $f: \mc X \rightarrow \mc Y$ of algebraic stacks locally of finite type, we denote by 
$$f_*: D_c^{(+)}(\mathcal{X}, \Oe) \rightarrow D_c^{(+)}(\mathcal{Y}, \Oe), \quad f_!: D_c^{(-)}(\mathcal{X}, \Oe) \rightarrow D_c^{(-)}(\mathcal{Y}, \Oe)$$
$$f^*: D_c(\mathcal{Y}, \Oe) \rightarrow D_c(\mathcal{X}, \Oe), \quad f^!: D_c(\mathcal{Y}, \Oe) \rightarrow D_c(\mathcal{X}, \Oe)$$
the corresponding functors, always understood in the derived sense.

In particular, we have the base change theorem (cf. \cite{LO08b} Section 12 at the level of cohomology sheaves, which would be enough for us, and \cite{LZ14} Sections 1 and 2 at the level of derived categories).

\sssec{}   \label{subsection-reduced-subscheme}
We work with étale cohomology. So for any stack (resp. scheme) (for example $\Cht_{G, N, I, W}$ and $\Gr_{G, I, W}$), we consider only the reduced substack (resp. subscheme) associated to it.

\subsection*{Acknowledgments}
I would like to thank Vincent Lafforgue, Gérard Laumon, Jack Thorne and Zhiyou Wu for stimulating discussions. I would like to thank Weizhe Zheng for answering my questions.

This work has received funding from the European Research Council (ERC) under the European Union's Horizon 2020 research and innovation programme (grant agreement No 714405).

\tableofcontents

\section{Definition of cohomology with integral coefficients of stacks of shtukas}

As in \cite{cusp-coho} Section 1.1, 
we denote by $\Cht_{G, N, I}$ 
the prestack of $G$-shtukas 
over $(X \sm N)^I$, by $\Gr_{G, I}$ 
the Beilinson-Drinfeld affine grassmannian over $X^I$ (which is a global version of the affine grassmannian $G((t)) / G[[t]]$) and by $G_{I, \infty}$ (resp. $\Gr_{I, d}$ for $d \in \N$) the group scheme over $X^I$ which is a global version of the positive loop group $G[[t]]$ (resp. of the jet group $G[t] / t^d$).
In the same way, we define $\Cht_{P, N, I}$, $\Gr_{P, I}$, $P_{I, \infty}$, $P_{I, d}$ and $\Cht_{M, N, I}$, $\Gr_{M, I}$, $M_{I, \infty}$, $M_{I, d}$.

Note that $loc.cit.$ Sections 1.2-1.7 are geometric. 

\quad

In the following, we will follow the constructions in $loc.cit.$ Section 2, with necessary modifications. As in $loc.cit.$, our results are of geometric nature, i.e. we will not consider the action of $\on{Gal}(\Fqbar / \Fq)$. 
From now on, we pass to the base change over $\Fqbar$. 

\subsection{Cohomology with integral coefficients of stacks of $G$-shtukas}   \label{subsection-integral-coho-Cht-G}

\sssec{}
The geometric Satake equivalence for the affine grassmannian is established in \cite{mv} (and reviewed in \cite{BR}) over the ground field $\C$, with coefficients in a Noetherian commutative ring of finite global dimension. By \cite{mv} Section 14, \cite{ga-de-jong} Section 1.6 and \cite{zhu}, the constructions in \cite{mv} can be extended to the case of ground field $\Fqbar$ and with coefficients in $\Oe$ or $\Ql$.

\sssec{}
As in \cite{mv} Section 2 and \cite{ga-de-jong} Section 2,
we denote by $\Perv_{G_{I, \infty}}(\Gr_{G, I}, \Oe)$ the category of $G_{I, \infty}$-equivariant perverse sheaves with $\Oe$-coefficients on $\Gr_{G, I}$, for the standard perverse $t$-structure $p$ (for which $\Perv_{G_{I, \infty}}(\Gr_{G, I}, \Oe)$ is Noetherian) and for the perverse normalization relative to $X^I$. 

\begin{rem}
The Verdier duality interchanges $\Perv_{G_{I, \infty}}(\Gr_{G, I}, \Oe)$ with the category of $G_{I, \infty}$-equivariant perverse sheaves with $\Oe$-coefficients on $\Gr_{G, I}$, for the costandard perverse $t$-structure $p+$ (cf. \cite{bbd} 3.3, 4.0, \cite{juteau}). 

In this paper, we will never use the perverse $t$-structure $p+$. We will never use the Verdier duality (except in Lemma \ref{lem-finite-gp-gerbe-counit-is-isom}, where we work on derived category and derived functors without considering any $t$-structure).
\end{rem}

\sssec{}   \label{subsection-def-G-dual}
Let $\wh G$ be the Langlands dual group of $G$ over $\Oe$. 
We denote by $\on{Rep}_{\Oe}(\wh G^I)$ the category of finite type $\Oe$-linear representations of $\wh G^I$. 
We will need the following theorem, which is a version of the geometric Satake equivalence in family:

\begin{thm} (\cite{mv} Theorem 14.1, \cite{ga-de-jong} Theorem 2.6, \cite{vincent} Théorème 1.17)  \label{thm-Satake-functor-I}
We have a canonical natural fully faithful $\Oe$-linear tensor functor:
\begin{equation}
\on{Sat}_{G, I}^{\Oe}: \on{Rep}_{\Oe}(\wh G^I) \rightarrow \Perv_{G_{I, \infty}}(\Gr_{G, I}, \Oe).
\end{equation}
It satisfies the properties in \cite{vincent} Théorème 1.17.
\cqfd
\end{thm}

\begin{rem}
By \cite{ga-de-jong} 2.5, Theorem 2.6 and the discussion after Theorem 2.6, we denote by $P^{\wh G, I}_{\Oe}$ the category of perverse sheaves with $\Oe$-coefficients on $X^I$ (for the perverse normalization relative to $X^I$) endowed with an extra structure given in {\it loc.cit}. There is a canonical equivalence of categories 
\begin{equation}
F: \Perv_{G_{I, \infty}}(\Gr_{G, I}, \Oe) \isom P^{\wh G, I}_{\Oe}
\end{equation}
compatible with the tensor structures defined in {\it loc.cit}. 

Besides, we have a fully faithful functor 
\begin{equation}
\Psi: \on{Rep}_{\Oe}(\wh G^I) \rightarrow P^{\wh G, I}_{\Oe}:  W \mapsto W \otimes_{\mc O_E} \mc O_{E, X^I} ,
\end{equation}
where $\mc O_{E, X^I}$ is the constant sheaf over $X^I$. The composition $F^{-1} \circ \Psi$ gives the functor $\on{Sat}_{G, I}^{\Oe}$.
\end{rem}

\begin{defi}  \label{def-Gr-G-I-W}
For any $W \in \on{Rep}_{\Oe}(\wh G^I) $, we define $\mc S_{G, I, W}^{\Oe} := \on{Sat}_{G, I}^{\Oe}(W)$.
We define $\Gr_{G, I, W}$ to be the support of $\mc S_{G, I, W}^{\Oe}$. It is a closed subscheme of finite type of $\Gr_{G, I}$.
\end{defi}

\sssec{}  \label{subsection-S-G-I-W-1-plus-W-2}
When $W = W_1 \oplus W_2$, by the functoriality of $\on{Sat}_{G, I}^{\Oe}$, we have $\mc S_{G, I, W}^{\Oe} = \mc S_{G, I, W_1}^{\Oe} \oplus \mc S_{G, I, W_2}^{\Oe}$. Then $\Gr_{G, I, W} = \Gr_{G, I, W_1} \cup \Gr_{G, I, W_2}$.

\sssec{}   \label{subsection-def-F-G-N-I-W}
\cite{cusp-coho} Sections 2.2-2.4 still hold. For any $W \in \on{Rep}_{\Oe}(\wh G^I) $, we define $\Cht_{G, N, I, W}$ to be the inverse image of $[G_{I, \infty} \backslash \Gr_{G, I, W}]$ by 
\begin{equation}   \label{equation-epsilon-G-N-I-infty}
\epsilon_{G, N, I, \infty}: \Cht_{G, N, I}  \rightarrow [G_{I, \infty} \backslash \Gr_{G, I}].
\end{equation}
We define $$\mc{F}_{G, N, I, W}^{\Oe} := (\epsilon_{G, N, I, \infty})^* \mc S_{G, I, W}^{\Oe}$$ 
It is a perverse sheaf on $\Cht_{G, N, I}$, supported on $\Cht_{G, N, I, W} $.

Similarly, we have 
\begin{equation}
\epsilon_{G, N, I, \infty}^{\Xi}: \Cht_{G, N, I} / \Xi  \rightarrow [G_{I, \infty}^{\on{ad}} \backslash \Gr_{G, I}].
\end{equation}
We define $$\mc{F}_{G, N, I, W}^{\Xi, \, \Oe} := (\epsilon_{G, N, I, \infty}^{\Xi})^* \mc S_{G, I, W}^{\Oe}$$
It is a perverse sheaf on $\Cht_{G, N, I} / \Xi$, supported on $\Cht_{G, N, I, W} / \Xi$.

\sssec{}  \label{subsection-F-G-I-W-1-plus-W-2}
When $W = W_1 \oplus W_2$, by \ref{subsection-S-G-I-W-1-plus-W-2} and \ref{subsection-def-F-G-N-I-W}, we have $\mc F_{G, N, I, W}^{\Oe} = \mc F_{G, N, I, W_1}^{\Oe} \oplus \mc F_{G, N, I, W_2}^{\Oe}$ and $\Cht_{G, I, W} = \Cht_{G, I, W_1} \cup \Cht_{G, I, W_2}$.

\quad

\sssec{}
In $loc.cit.$ 1.7, for any $\mu \in \wh{\Lambda}_{G^{\mr{ad}}}^{+, \Q}$, we defined Harder-Narasimhan open substack $\Cht_{G, N, I}^{\leq \mu}$ of $\Cht_{G, N, I}$. Let $\Cht_{G, N, I, W}^{\leq \mu} := \Cht_{G, N, I}^{\leq \mu} \cap \Cht_{G, N, I, W}$. It is a Deligne-Mumford stack of finite type.

In $loc.cit.$ 1.1.7, we defined the morphism of paws $\mf{p}_G: \Cht_{G, I, N} / \Xi \rightarrow (X \sm N)^I.$ 

\begin{defi}  \label{def-mc-H-G-leq-mu}
For any $\mu \in \wh{\Lambda}_{G^{\mr{ad}}}^{+, \Q}$, we define 
$$
\mc H_{G, N, I, W}^{\leq \mu, \, \Oe} := R(\mf{p}_G)_! ( \restr{\mc{F}_{G, N, I, W}^{\Xi, \, \Oe} } {\Cht_{G, N, I, W}^{\leq \mu} / \Xi } ) \in D_c^{(-)}( (X \sm N)^I , \Oe ).
$$
For any $j \in \Z$, we define degree $j$ cohomology sheaf (for the ordinary $t$-structure):
$$
\mc H _{G, N, I, W}^{j, \, \leq\mu, \, \Oe}:=R^j(\mf{p}_G)_! ( \restr{\mc{F}_{G, N, I, W}^{\Xi, \, \Oe} } {\Cht_{G, N, I, W}^{\leq \mu} / \Xi } ).
$$
This is a $\Oe$-constructible sheaf on $ (X \sm N)^I$. 
\end{defi}

\begin{rem}
It is necessary to consider $D_c^{(-)}( (X \sm N)^I , \Oe )$ rather than $D_c^{(b)}( (X \sm N)^I , \Oe ).$ From \cite{LO08a} 4.9.2, we see that $H_c^*(B(\Z / \ell \Z), \Zl)$ is unbounded below, where $B(\Z / \ell \Z)$ is the classifying stack of $\Z / \ell \Z$.
\end{rem}

\sssec{}   \label{subsection-H-G-leq-mu-1-to-leq-mu-2}
Let $\mu_1, \mu_2 \in \wh{\Lambda}_{G^{\mr{ad}}}^{+, \Q}$ and $\mu_1 \leq \mu_2$. We have an open immersion:
\begin{equation}   \label{equation-Cht-G-leq-mu-1-to-leq-mu-2}
\Cht_{G, N, I, W}^{\leq \mu_1} / \Xi  \hookrightarrow \Cht_{G, N, I, W}^{\leq \mu_2} / \Xi.
\end{equation}
For any $j$, morphism (\ref{equation-Cht-G-leq-mu-1-to-leq-mu-2}) induces a morphism of sheaves:
$$\mc H _{G, N, I, W}^{j, \, \leq\mu_1, \, \Oe} \rightarrow \mc H _{G, N, I, W}^{j, \, \leq \mu_2, \, \Oe}.$$

\begin{defi}   \label{def-mc-H-ind-limit}
We define
$$\mc H _{G, N, I, W}^{j, \, \Oe}: = \varinjlim _{\mu}  \mc H _{G, N, I, W}^{j, \, \leq\mu, \, \Oe} $$
as an inductive limit in the category of constructible sheaves on $ (X \sm N)^I$. 
\end{defi}

\subsection{Cohomology with integral coefficients of stacks of $M$-shtukas}      \label{subsection-coho-Cht-M}

Let $P$ be a proper parabolic subgroup of $G$ and let $M$ be its Levi quotient. 

\sssec{}
Let $\wh M$ be the Langlands dual group of $M$ over $\Oe$ defined by the geometric Satake equivalence.
The compatibility between the geometric Satake equivalence and the constant term functor along $P$ (\cite{bg} Theorem 4.3.4, \cite{BR} Proposition 15.3) induces a canonical inclusion $\wh M \hookrightarrow \wh G$ (compatible with pinning).

\sssec{}    \label{subsection-epsilon-M-N-I-d}
We view $W \in \on{Rep}_{\Oe}(\wh G^I)$ as a representation of $\wh M^I$ via $\wh M^I \hookrightarrow \wh G^I$. 
Let $\Gr_{M, I}$ be the Beilinson-Drinfeld affine grassmannian over $X^I$ associated to $M$.
As in Definition \ref{def-Gr-G-I-W} (replacing $G$ by $M$), we define $\Gr_{M, I, W}$ and perverse sheaf $\mc S_{M, I, W}^{\Oe}$ supported on $\Gr_{M, I, W}$.
As in \ref{subsection-def-F-G-N-I-W} (replacing $G$ by $M$), we define $\Cht_{M, N, I, W}$ and perverse sheaf $\mc{F}_{M, N, I, W}^{\Oe}$ (resp. $\mc{F}_{M, N, I, W}^{\Xi, \, \Oe}$) supported on $\Cht_{M, N, I, W}$ (resp. $\Cht_{M, N, I, W} / \Xi$).

\sssec{}       \label{subsection-Cht-M-leq-mu-nu-empty-for-nu-big}
In \cite{cusp-coho} 1.5.11, we defined $\on{pr}_P^{\on{ad}}: \wh \Lambda_{G^{\mr{ad}}}^{\Q} \rightarrow \wh \Lambda_{Z_M / Z_G}^{\Q}$. In $loc.cit.$ 1.5.13, we defined a partial order on $\wh \Lambda_{Z_M / Z_G}^{\Q}$: $\mu_1 \leq^{G^{\mr{ad}}} \mu_2 \Leftrightarrow \mu_2 - \mu_1 \in \on{pr}_P^{\on{ad}} ( \wh \Lambda_{G^{\mr{ad}}}^{\on{pos}, \, \Q})$. In $loc.cit.$ 1.5.20, for any $\mu \in \wh{\Lambda}_{G^{\mr{ad}}}^{+, \Q}$, we defined a translated cone in $\wh \Lambda_{Z_M / Z_G}^{\Q}$:
\begin{equation}
\wh{\Lambda}_{Z_M / Z_G}^{\mu} := \{  \nu \in  \wh \Lambda_{Z_M / Z_G}^{\Q}, \nu \leq^{G^{\mr{ad}}} \on{pr}_P^{\on{ad}}(\mu) \} .
\end{equation}

In \cite{cusp-coho} Section 1.7, for any $\mu \in \wh{\Lambda}_{G^{\mr{ad}}}^{+, \Q}$, we defined an open substack $\Cht_{M, N, I}^{\leq \mu}$ of $\Cht_{M, N, I}$. For any $\nu \in \wh \Lambda_{Z_M / Z_G}^{\Q}$, we defined an open and closed substack $\Cht_{M, N, I}^{\leq \mu, \, \nu}$ of $\Cht_{M, N, I}^{\leq \mu}$.
By $loc.cit.$ 2.6.5, if $\nu \notin \wh{\Lambda}_{Z_M / Z_G}^{\mu}$, then $\Cht_{M, N, I}^{\leq \mu, \, \nu}$ is empty.

Let $\Cht_{M, N, I, W}^{\leq \mu, \, \nu} := \Cht_{M, N, I}^{\leq \mu, \, \nu} \cap \Cht_{M, N, I, W}$.
As in $loc.cit.$ 2.6.3, we have a decomposition
\begin{equation}
\Cht_{M, N, I, W}^{\leq \mu} / \Xi = \underset{ \nu \in \wh{\Lambda}_{Z_M / Z_G}^{\mu}  } \bigsqcup \Cht_{M, N, I, W}^{\leq \mu,  \; \nu}  / \Xi .
\end{equation}
where each $\Cht_{M, N, I, W}^{\leq \mu, \; \nu} / \Xi$ is a Deligne-Mumford stack of finite type.

\quad

We have the morphism of paws $\mf{p}_M: \Cht_{M, I, N} / \Xi \rightarrow (X \sm N)^I.$
\begin{defi}   \label{def-H-M-mu-nu}
For any $\mu \in \wh{\Lambda}_{G^{\mr{ad}}}^{+, \Q}$ and $\nu \in \wh{\Lambda}_{Z_M / Z_G}^{\Q}$, we define 
$$\mc H _{M, N, I, W}^{\leq  \mu, \, \nu, \, \Oe}  := R(\mf{p}_M)_! ( \restr{ \mc{F}_{M, N, I, W}^{\Xi, \, \Oe}    } {\Cht_{M, N, I, W}^{\leq  \mu, \; \nu} / \Xi} )  \in  D_c^{(-)}( (X \sm N)^I , \Oe ) ;$$
For any $j \in \Z$, we define the degree $j$ cohomology sheaf (for the ordinary $t$-structure):
$$\mc H _{M, N, I, W}^{j, \, \leq \mu, \, \nu, \, \Oe}  := R^j(\mf{p}_M)_! ( \restr{ \mc{F}_{M, N, I, W}^{\Xi, \, \Oe}    } {\Cht_{M, N, I, W}^{\leq  \mu, \; \nu} / \Xi} ) .$$ 
This is a $\Oe$-constructible sheaf on $ (X \sm N)^I$. 
\end{defi}

\begin{defi}    \label{def-H-M-j-leq-mu}
Let $\ov{x}$ be a geometric point of $(X \sm N)^I$.
We define 
\begin{equation}
\restr{ \mc H _{M, N, I, W}^{j, \, \leq  \mu, \, \Oe} }{ \ov{x} } := \prod_{\nu \in \wh{\Lambda}_{Z_M / Z_G}^{\mu}} \restr{  \mc H _{M, N, I, W}^{j, \, \leq  \mu, \, \nu, \, \Oe}  }{  \ov{x}  }.
\end{equation}
\end{defi}

\sssec{}   \label{subsection-H-M-leq-mu-nu-zero-for-nu-big}
Note that by \ref{subsection-Cht-M-leq-mu-nu-empty-for-nu-big}, if $\nu \in \wh \Lambda_{Z_M / Z_G}^{\Q}$ but $\nu \notin \wh{\Lambda}_{Z_M / Z_G}^{\mu}$, then $\restr{  \mc H _{M, N, I, W}^{j, \, \leq  \mu, \, \nu, \, \Oe}  }{  \ov{x}  }=0$.

\sssec{}
Let $\mu_1, \mu_2 \in \wh{\Lambda}_{G^{\mr{ad}}}^{+, \Q}$ and $\mu_1 \leq \mu_2$. We have an open immersion:
\begin{equation}    \label{equation-Cht-M-leq-mu-1-leq-mu-2}
\Cht_{M, N, I, W}^{\leq \mu_1} / \Xi \hookrightarrow \Cht_{M, N, I, W}^{\leq \mu_2} / \Xi.
\end{equation}
For any $j$ and $\nu$, morphism (\ref{equation-Cht-M-leq-mu-1-leq-mu-2}) induces a morphism of sheaves:
\begin{equation}
\mc H _{M, N, I, W}^{j, \, \leq\mu_1, \, \nu, \, \Oe} \rightarrow \mc H _{M, N, I, W}^{j, \, \leq \mu_2, \, \nu, \, \Oe}.
\end{equation}

We deduce a morphism of $\Oe$-modules:
\begin{equation}
\restr{ \mc H _{M, N, I, W}^{j, \, \leq  \mu_1, \, \Oe} }{ \ov{x} }  \rightarrow \restr{ \mc H _{M, N, I, W}^{j, \, \leq  \mu_2, \, \Oe} }{ \ov{x} } . 
\end{equation}

\begin{defi}   \label{def-H-M-j}
We define
$$   
\restr{ \mc H _{M, N, I, W}^{j,  \, \Oe} }{ \ov{x} } : = \varinjlim _{\mu}  \restr{ \mc H _{M, N, I, W}^{j, \, \leq  \mu, \, \Oe} }{ \ov{x} } 
$$
as an inductive limit in the category of $\Oe$-modules.
\end{defi}

\begin{defi}    \label{def-H-M-j-nu}
For any $\nu \in \wh{\Lambda}_{Z_M / Z_G}^{\Q}$, we define $\restr{ \mc H _{M, N, I, W}^{j, \, \nu, \, \Oe} }{ \ov{x} } : = \varinjlim _{\mu}  \restr{ \mc H _{M, N, I, W}^{j, \, \leq  \mu, \, \nu, \, \Oe} }{ \ov{x} } $ as an inductive limit in the category of $\Oe$-modules.
\end{defi}

\quad

\section{Constant term morphisms and cuspidal cohomology}

The goal of this section is to extend \cite{cusp-coho} Section 3 to the case of $\Oe$-coefficients. 

Let $P$ be a parabolic subgroup of $G$ and $M$ its Levi quotient. Let $U$ be the unipotent radical of $P$. 
Let $W \in \on{Rep}_{\Oe}(\wh G^I)$. For any geometric point $\ov{x}$ of $(X \sm N)^I$ as in \ref{subsection-ov-x} below,
we will construct a constant term morphism from $\restr{ \mc H _{G, N, I, W}^{j, \, \Oe} }{ \ov{x} } $ to $\restr{ \mc H _{M, N, I, W}^{j, \, \Oe} }{ \ov{x} } $ (in fact, to a variant $\restr{ \mc H _{M, N, I, W}^{' \, j, \, \Oe} }{ \ov{x} } $ of $\restr{ \mc H _{M, N, I, W}^{j, \, \Oe} }{ \ov{x} } $ defined in Definition \ref{def-mc-H-M-prime} below). 

\subsection{Some geometry}   \label{subsection-some-geo}

\sssec{}
The morphisms of groups $G \hookleftarrow P \twoheadrightarrow M$ induce morphisms of ind-schemes over $X^I$: 
\begin{equation}   \label{equation-Gr-G-I-Gr-P-I-Gr-M-I}
\Gr_{G, I} \xleftarrow{i^0} \Gr_{P, I} \xrightarrow{\pi^0} \Gr_{M, I}
\end{equation}
In \cite{cusp-coho} Construction 1.2.4, we constructed morphisms over $(X \sm N)^I$:
\begin{equation}   \label{equation-Cht-G-I-Cht-P-I-Cht-M-I}
\Cht_{G, N, I} \xleftarrow{i} \Cht_{P, N, I} \xrightarrow{\pi} \Cht_{M, N, I}
\end{equation}

\sssec{}
Let $\underline{\omega}=(\omega_i)_{i \in I} \in (\wh \Lambda_G^+)^I$. Let $\Gr_{G, I, \underline{\omega}}$ be the closed substack of $\Gr_{G, I}$ defined in \cite{vincent} Définition 1.12. We have 
\begin{equation*}   
\Gr_{G, I} = \bigcup_{ \underline{\omega} \in (\wh \Lambda_G^+)^I } \Gr_{G, I, \underline{\omega}} .
\end{equation*}

\sssec{}   \label{subsection-def-Gr-M-I-lambda}
Let $\underline{\lambda}=(\lambda)_{i \in I} \in (\wh \Lambda_M^+)^I $. Let $\Gr_{M, I, \underline{\lambda}}$ be the closed substack of $\Gr_{M, I}$ defined in \cite{vincent} Définition 1.12. For any $\underline{\omega}=(\omega_i)_{i \in I} \in (\wh \Lambda_G^+)^I$, let
$$
C_{\underline{\omega}}:=\{  (\lambda_i)_{i \in I} \in (\wh \Lambda_G)^I, \forall i \in I, \lambda_i \text{ is conjugated to a dominant coweight} \leq \omega_i    \}
$$
We define
\begin{equation}
\Gr_{M, I, \underline{\omega}}^* := \bigcup_{ \underline{\lambda} \in  (\wh \Lambda_M^+)^I \cap  C_{\underline{\omega}}   } \Gr_{M, I, \underline{\lambda}}
\end{equation}
We have 
\begin{equation*}    
\Gr_{M, I} = \bigcup_{ \underline{\omega} \in (\wh \Lambda_G^+)^I } \Gr_{M, I, \underline{\omega}}^* .
\end{equation*}

\sssec{}
Let $\underline{\omega} \in (\wh \Lambda_G^+)^I$. 
Let $\Cht_{G, N, I, \underline{\omega}}$ be the closed substack of $\Cht_{G, N, I}$ defined to be the inverse image of $[G_{I, \infty} \backslash \Gr_{G, I, \underline{\omega}}]$ by
\begin{equation*}
\epsilon_{G, N, I, \infty}: \Cht_{G, N, I} \rightarrow [G_{I, \infty} \backslash \Gr_{G, I}]
\end{equation*}
Let $\Cht_{M, N, I, \underline{\omega}}^*$ be the closed substack of $\Cht_{M, N, I}$ defined to be the inverse image of 
$[M_{I, \infty} \backslash \Gr_{M, I, \underline{\omega}}^* ]$ by 
\begin{equation*}
 \epsilon_{M, N, I, \infty}: \Cht_{M, N, I} \rightarrow [M_{I, \infty} \backslash \Gr_{M, I}]
\end{equation*}

\sssec{}
\cite{cusp-coho} Section 3.1 remains valid if we replace everywhere $W$ by $\underline{\omega} \in (\wh \Lambda_G^+)^I$, replace $\Gr_{M, I, W}$ by $\Gr_{M, I, \underline{\omega}}^*$ and replace $\Cht_{M, I, W}$ by $\Cht_{M, I, \underline{\omega}}^*$.
As in $loc.cit.$ Proposition 3.1.1, we have 
\begin{equation}
(i^0)^{-1}( \Gr_{G, I, \underline{\omega}} ) \subset (\pi^0)^{-1} ( \Gr_{M, I, \underline{\omega}}^*  ),
\end{equation}
where the inverse images are in the sense of reduced subschemes in $\Gr_{P, I}$. We define $\Gr_{P, I, \underline{\omega}}: = (i^0)^{-1}( \Gr_{G, I, \underline{\omega}} )$. Morphisms (\ref{equation-Gr-G-I-Gr-P-I-Gr-M-I}) induce morphisms 
\begin{equation}   \label{equation-Gr-G-I-omega-P-I-omega-M-I-omega}
\Gr_{G, I, \underline{\omega}} \xleftarrow{i^0} \Gr_{P, I, \underline{\omega}} \xrightarrow{\pi^0} \Gr_{M, I, \underline{\omega}}^*
\end{equation}

Similarly, we define $\Cht_{P, N, I, \underline{\omega}}: = i^{-1}( \Cht_{G, N, I, \underline{\omega}} )$. Morphisms (\ref{equation-Cht-G-I-Cht-P-I-Cht-M-I}) induce morphisms 
\begin{equation}
\Cht_{G, N, I, \underline{\omega}} \xleftarrow{i} \Cht_{P, N, I, \underline{\omega}} \xrightarrow{\pi} \Cht_{M, N, I, \underline{\omega}}^*
\end{equation}

\begin{rem}
If $W = \boxtimes_{i \in I} W_i$, where $W_i$ is the irreducible $E$-representation of $\wh G_E$ of highest weight $\omega_i$. Then
\begin{equation}   \label{equation-Gr-M-I-omega-*-equal-Gr-M-I-W}
\Gr_{G, I, \underline{\omega}} = \Gr_{G, I, W} \quad \text{and} \quad \Gr_{M, I, \underline{\omega}}^* = \Gr_{M, I, W}
\end{equation}
where $\Gr_{G, I, W}$ and $\Gr_{M, I, W}$ are defined in \cite{cusp-coho} Definition 2.1.8 and 2.6.2.

However, if $W$ is an $\Oe$-representation of $\wh G$, due to possible torsions in the stalks of $\mc S_{G, I, W}^{\Oe}$ and $\mc S_{M, I, W}^{\Oe}$, (\ref{equation-Gr-M-I-omega-*-equal-Gr-M-I-W}) may not be true. In this case, I do not know if $(i^0)^{-1}( \Gr_{G, I, W} ) \subset (\pi^0)^{-1} ( \Gr_{M, I, W}  )$. But we can always find large enough strata in $\Gr_{G, I}$ (resp. $\Gr_{M, I}$) such that the perverse sheaf $\mc S_{G, I, W}^{\Oe}$ (resp. $\mc S_{M, I, W}^{\Oe}$) is supported on them, in the following way.
\end{rem}

\sssec{}      \label{subsection-Gr-G-I-W-included-in-Gr-G-I-omega}
Let $W \in \on{Rep}_{\Oe}(\wh G^I) $. There exists a finite subset $\Sigma \subset (\wh \Lambda_G^+)^I$ (depending on $W$) such that 
\begin{equation}
\Gr_{G, I, W} \subset \bigsqcup_{\underline{\omega} \in \Sigma} \Gr_{G, I, \underline{\omega} } \quad \text{and} \quad  \Gr_{M, I, W} \subset \bigsqcup_{\underline{\omega} \in \Sigma } \Gr_{M, I, \underline{\omega}}^* .
\end{equation}

We define 
\begin{equation}
\Gr_{G, I, \Sigma }:=\bigsqcup_{\underline{\omega} \in \Sigma} \Gr_{G, I, \underline{\omega} }, \quad \Gr_{P, I, \Sigma }:=\bigsqcup_{\underline{\omega} \in \Sigma} \Gr_{P, I, \underline{\omega} }, \quad \Gr_{M, I, \Sigma } := \bigsqcup_{\underline{\omega} \in \Sigma} \Gr_{M, I, \underline{\omega}}^* . 
\end{equation}

Note that $\mc S_{G, I, W}^{\Oe}$ is supported on $\Gr_{G, I, \Sigma }$ and $\mc S_{M, I, W}^{\Oe}$ is supported on $\Gr_{M, I, \Sigma }$.

\sssec{}    \label{subsection-def-Cht-G-Sigma-Cht-M-Sigma}
We define 
$$\Cht_{G, N, I, \Sigma }:=\bigsqcup_{\underline{\omega} \in \Sigma} \Cht_{G,N,  I, \underline{\omega} }, $$
$$ \Cht_{P, N, I, \Sigma }:=\bigsqcup_{\underline{\omega} \in \Sigma} \Cht_{P, N, I, \underline{\omega} }, \quad  \Cht_{M, N, I, \Sigma } := \bigsqcup_{\underline{\omega} \in \Sigma} \Cht_{M, N, I, \underline{\omega}}^* $$

As \cite{cusp-coho} 3.1.4, for $d \in \Z_{\leq 0}$ large enough depending on $\Sigma$ (thus depending on $W$) as in $loc.cit.$ Proposition 2.2.1 applied to $\Gr_{G, I, \Sigma}$ and $\Gr_{M, I, W}$, we have a commutative diagram of algebraic stacks: 
\begin{equation}     \label{diagram-Cht-Gr-G-P-M}
\xymatrix{
 \Cht_{G, N, I, \Sigma }  \ar[d]^{\epsilon_{G, d}} 
&  \Cht_{P, N, I, \Sigma }  \ar[l]_{i}   \ar[d]^{\epsilon_{P, d}}  \ar[r]^{\pi} 
&  \Cht_{M, N, I, \Sigma }    \ar[d]^{ \epsilon_{M, d}} \\
[G_{I, d} \backslash  \Gr_{G, I, \Sigma } ]       
& [ P_{I, d} \backslash   \Gr_{P, I, \Sigma } ]   \ar[l]_{ \ov{i^0_{d}} }  \ar[r]^{ \ov{\pi^0_{d}} } 
&  [  M_{I, d} \backslash  \Gr_{M, I, \Sigma } ]
}
\end{equation}

Note that $\mc F_{G, N, I, W}^{\Oe}$ is supported on $\Cht_{G, N, I, \Sigma }$ and $\mc F_{M, N, I, W}^{\Oe}$ is supported on $\Cht_{M, N, I, \Sigma }$.

\sssec{}    \label{subsection-pi-d-smooth}
The squares in (\ref{diagram-Cht-Gr-G-P-M}) are not Cartesian. Consider the right square. As in \cite{cusp-coho} 3.1.5, we have a commutative diagram, where the square is Cartesian:
\begin{equation}   \label{diagram-wt-Cht-M}
\xymatrix{
\Cht_{P, N, I, \Sigma} \ar[rrd]^{\pi}  \ar[rd]|-{\pi_d}  \ar[rdd]_{\epsilon_{P, d}}  \\
& \wt{\Cht}_{M, N, I, \Sigma }   \ar[r]_{\wt{\pi_d^0}}  \ar[d]^{\wt{\epsilon_{M, d}}}
& \Cht_{M, N, I, \Sigma } \ar[d]^{\epsilon_{M, d}} \\
& [ P_{I, d} \backslash   \Gr_{P, I, \Sigma } ]   \ar[r]^{ \ov{\pi^0_{d}} } 
&  [  M_{I, d} \backslash  \Gr_{M, I, \Sigma } ]
}
\end{equation}
Let $U_{I, d}$ be the group scheme over $X^I$ as in \cite{cusp-coho} Definition 1.1.13 applied to $U$, which is a global version of the jet group $U[t] / t^d$. 
As in $loc.cit.$ Lemma 3.1.8, the morphism $\pi_d$ is smooth of relative dimension $\dim_{X^I} U_{I, d}$.

\subsection{Compatibility of the geometric Satake equivalence and parabolic induction}   \label{subsection-Satake-geo-para-ind}

The goal of this section is to recall Theorem \ref{thm-geo-satake-CT-I-paws} and deduce (\ref{equation-inverse-image-S-M-W-equal-S-M-W-theta-0}), which is the key ingredient for the next section. We rewrite \cite{cusp-coho} Section 3.2 in a slightly more general way.

\sssec{}   \label{subsection-def-Gr-M-I-n}
As in \cite{cusp-coho} 3.2.1, we denote by $\rho_G$ (resp. $\rho_M$) the half sum of positive roots of $G$ (resp. $M$). 
Then $2(\rho_G - \rho_M)$ is a character of $M$ (the determinant of the adjoint action on $\on{Lie} U$).
The morphism $2(\rho_G - \rho_M): M \rightarrow \mb G_m$ induces a morphism $\Gr_{M, I} \rightarrow \Gr_{\mb G_m, I}$ by sending a $M$-bundle $\mc M$ to the $\mb G_m$-bundle $\mc M \overset{M}{\times} \mb G_m$.
We have a morphism $\on{deg}: \Gr_{\mb G_m, I}  \rightarrow \Z$ by taking the degree of a $\mb G_m$-bundle. 
We have the composition of morphisms
\begin{equation}    \label{equation-Gr-M-I-to-Z}
\Gr_{M, I} \rightarrow \Gr_{\mb G_m, I}  \xrightarrow{\on{deg}}  \Z.
\end{equation}
We define $\Gr_{M, I}^n$ to be the inverse image of $n \in \Z$. It is open and closed in $\Gr_{M, I}$.

\sssec{}

We define $\Gr_{P, I}^n := (\pi^0)^{-1} \Gr_{M, I}^n$. Morphism (\ref{equation-Gr-G-I-Gr-P-I-Gr-M-I}) induces a morphism
\begin{equation}   \label{equation-Gr-G-I-Gr-P-I-n-Gr-M-I-n}
\Gr_{G, I} \xleftarrow{i^0_n} \Gr_{P, I}^n \xrightarrow{\pi^0_n} \Gr_{M, I}^n.
\end{equation}

\sssec{}   \label{subsection-def-theta-rho-G-rho-M}
As in $loc.cit.$ 3.2.2,
we define $\wh \Lambda_{M, M}:= \wh \Lambda_M / \wh\Lambda_{[M, M]_{\mr{sc}}},$ where $\wh\Lambda_{[M, M]_{\mr{sc}}} \subset \wh \Lambda_M$ is the sublattice spanned by the simple coroots of $M$. 
By definition, it coincides with $\pi_1(M)$ defined in \cite{var} Lemma 2.2 (which is caconically isomorphic to $\Lambda_{Z_{\wh M}}$, the group of characters of $Z_{\wh M}$). For any $\check{\alpha} \in \wh\Lambda_{[M, M]_{\mr{sc}}}$, we have $\langle \check{\alpha} , 2\rho_G - 2\rho_M \rangle=0$. Thus the pairing $\langle \theta , 2\rho_G - 2\rho_M \rangle$ makes sense for $\theta \in \wh \Lambda_{M, M}$.

\begin{rem}   \label{rem-Gr-M-theta-Gr-M-n}
(We do not need the following fact.) When $I$ is a singleton, the connected components of $\Gr_{M, I}$ are indexed by $\Lambda_{Z_{\wh M}}$. For $\theta \in \Lambda_{Z_{\wh M}}$, we denote by $\Gr_{M, I}^{\theta}$ the connected component corresponding to $\theta$. Then we have 
$$\Gr_{M, I}^n = \bigsqcup_{  \theta \in \Lambda_{Z_{\wh M}} , \, \langle \theta , 2(\rho_G - \rho_M) \rangle = n   } \Gr_{M, I}^{\theta} .$$
\end{rem}

\sssec{}
As in \cite{AHR} Section 4.1, we define a functor from $D^b_c(\Gr_{G, I}, \Oe) $ to $ D^b_c(\Gr_{M, I}, \Oe)$:
$$r^{GM}:= \bigoplus_{n \in \Z} (\pi_{n}^0)_! (i_{n}^0)^*\otimes \big( \Oe [1](\frac{1}{2}) \big)^{\otimes n } .$$

\sssec{}
In Theorem \ref{thm-Satake-functor-I}, we defined a fully faithful functor $\on{Sat}_{G, I}^{\Oe}$. We denote by $\Perv_{G_{I, \infty}}(\Gr_{G, I}, \Oe)^{\on{MV}}$ its essential image. Similarly, we define the functor $\on{Sat}_{M, I}^{\Oe}: \on{Rep}_{\Oe}(\wh M^I)  \rightarrow \Perv_{M_{I, \infty}}(\Gr_{M, I}, \Oe)$ and its essential image $\Perv_{M_{I, \infty}}(\Gr_{M, I}, \Oe)^{\on{MV}}$.

\begin{thm}  \label{thm-geo-satake-CT-I-paws}  (\cite{bd} 5.3.29, \cite{mv} Theorem 3.6 for $M=T$, \cite{AHR} Lemma 4.1, \cite{BR} Proposition 15.2)

(a) The functor $r^{GM}$ sends $\Perv_{G_{I, \infty}}(\Gr_{G, I}, \Oe)^{\on{MV}}$ to $\Perv_{M_{I, \infty}}(\Gr_{M, I}, \Oe)^{\on{MV}}$.

(b) There is a canonical isomorphism of tensor functors 
\begin{equation}   \label{equation-Sat-M-Res-equal-CT-Sat-G}
\on{Sat}_{M, I}^{\Oe} \circ \on{Res}^{\wh G^I}_{\wh M^I} = r^{GM} \circ \on{Sat}_{G, I}^{\Oe}.
\end{equation} 
In other words, the following diagram of categories canonically commutes:
\begin{equation}
\xymatrixrowsep{2pc}
\xymatrixcolsep{6pc}
\xymatrix{
\Perv_{G_{I, \infty}}(\Gr_{G, I}, \Oe)^{\on{MV}}   \ar[r]^{  r^{GM}  } 
& \Perv_{M_{I, \infty}}(\Gr_{M, I}, \Oe)^{\on{MV}}  \\
\on{Rep}_{\Oe}(\wh G^I)   \ar[u]^{ \on{Sat}_{G, I}^{\Oe} }   \ar[r]^{ \on{Res}^{\wh G^I}_{\wh M^I} }
& \on{Rep}_{\Oe}(\wh M^I)   \ar[u]^{ \on{Sat}_{M, I}^{\Oe} }
}
\end{equation}
\cqfd
\end{thm}

\begin{rem}
The references cited above in Theorem \ref{thm-geo-satake-CT-I-paws} are for the case where $I$ is a singleton (i.e. for affine grassmannians). 
The general case where $I$ is arbitrary (i.e. for Beilinson-Drinfeld grassmannians) can be deduced from the case of singleton $I$ by the fact that the constant term functor $r^{GM}$ commutes with fusion (i.e. convolution). The proof for $I = \{1, 2\}$ is already included in the proof of Proposition 15.2 in \cite{BR}. For general $I$ the proof is similar. 
\end{rem}

\sssec{}
Applying (\ref{equation-Sat-M-Res-equal-CT-Sat-G}) to $W \in \on{Rep}_{\Oe}(\wh G^I) $, we obtain a canonical isomorphism
\begin{equation}   \label{equation-TC-satake-geo}
\mc{S}_{M, I, W}^{\Oe} \simeq \bigoplus_{n \in \Z} (\pi_{n}^0)_! (i_{n}^0)^* \mc S_{G, I, W}^{\Oe} [ n] (  n/2 ).
\end{equation}

\quad

\sssec{}     \label{subsection-TC-satake-geo-en-quotient}
For any $n$, we denote by $\Gr_{M, I, \Sigma}^{n}=\Gr_{M, I}^n \cap \Gr_{M, I, \Sigma}$ and $\Gr_{P, I, \Sigma}^n=\Gr_{P, I}^n \cap \Gr_{P, I, \Sigma}$. As \cite{cusp-coho} 3.2.9, we have a commutative diagram:
\begin{equation}  
\xymatrix{
 \Gr_{G, I, \Sigma }  \ar[d]^{\xi_{G, d}} 
&  \Gr_{P, I, \Sigma }^n  \ar[l]_{i^0_{n}}   \ar[d]^{\xi_{P, d}}  \ar[r]^{\pi^0_{n}} 
&   \Gr_{M, I, \Sigma }^{ n}    \ar[d]^{\xi_{M, d}} \\
[G_{I, d} \backslash  \Gr_{G, I, \Sigma }]       
& [ P_{I, d} \backslash   \Gr_{P, I, \Sigma }^n]   \ar[l]_{ \ov{i^0_{d, n}} }  \ar[r]^{ \ov{\pi^0_{d, n}} } 
&  [  M_{I, d} \backslash   \Gr_{M, I, \Sigma }^{n}]
}
\end{equation}

The squares are not Cartesian. Consider the right square.
The morphism $$\Gr_{P, I, \Sigma }^n \rightarrow [ P_{I, d} \backslash   \Gr_{P, I, \Sigma }^n]  \underset{  [  M_{I, d} \backslash   \Gr_{M, I, \Sigma }^{n}] }{\times} \Gr_{M, I, \Sigma }^{ n}  =  [ U_{I, d} \backslash  \Gr_{P, I, \Sigma }^n]$$
is a $U_{I, d}$-torsor. We deduce from the fact that the group scheme $U_{I, d}$ is unipotent over $X^I$ and from the proper base change (see \ref{def-D-c-b-pour-chat-alg}) that
\begin{equation}   \label{equation-pi-xi-shift-m}
(\pi_{n}^0)_! (\xi_{P, d})^* \simeq (\xi_{M, d})^* (\ov{\pi_{d, n}^0})_! [-2m](-m) ,
\end{equation}
where $m= \dim \xi_{P, d} - \dim \xi_{M, d} = \dim_{X^I} U_{I, d}$. 

\sssec{}
As in $loc.cit.$ 2.2.3, let $\mc S_{G, I, W}^{d, \, \Oe}$ (resp. $\mc S_{M, I, W}^{d, \, \Oe} $) be the (shifted) perverse sheaf on $[G_{I, d} \backslash   \Gr_{G, I, W}]$ (resp. $[  M_{I, d} \backslash  \Gr_{M, I, W }]$) such that $\mc S_{G, I, W}^{\Oe} = (\xi_{G, d})^* \mc S_{G, I, W}^{d, \, \Oe}$ (resp. $\mc S_{M, I, W}^{\Oe} = (\xi_{M, d})^* \mc S_{M, I, W}^{d, \, \Oe}$).
Taking into account \ref{subsection-Gr-G-I-W-included-in-Gr-G-I-omega}, we deduce from (\ref{equation-TC-satake-geo}) and (\ref{equation-pi-xi-shift-m}) that
\begin{equation}     \label{equation-TC-satake-geo-en-quotient}
\mc S_{M, I, W}^{d, \, \Oe} \isom \big( \bigoplus_{n \in \Z} (\ov{\pi_{d, n}^0})_! (\ov{i_{d, n}^0})^* \mc S_{G, I, W}^{d, \, \Oe}   [n](n/2) \big) [-2m](-m).
\end{equation} 

\quad

\sssec{}   \label{subsection-Gr-M-I-n-Gr-M-I-omega}
For $\underline{\lambda}=(\lambda_i)_{i \in I} \in (\wh \Lambda_M^+)^I $. 
We denote by $[\sum_{i \in I} \lambda_i]$ the image of $\sum_{i \in I} \lambda_i$ by the projection $\wh \Lambda_M \twoheadrightarrow \wh \Lambda_{M, M}$. By \ref{subsection-def-Gr-M-I-n} and \ref{subsection-def-theta-rho-G-rho-M}, we deduce that
\begin{equation}   
\Gr_{M, I}^n = \bigcup_{ \underline{\lambda}=(\lambda_i)_{i \in I} \in (\wh \Lambda_M^+)^I, \, \langle [\sum_{i \in I} \lambda_i] , 2(\rho_G - \rho_M) \rangle = n    } \Gr_{M, I, \underline{\lambda}} 
\end{equation}
where $\Gr_{M, I, \underline{\lambda}}$ is defined in \ref{subsection-def-Gr-M-I-lambda}.

\sssec{}   \label{subsection-Cht-M-omega-empty}
Let $\Cht_{M, N, I}^n$ (resp. $\Cht_{M, N, I, \underline{\lambda}}$) be the inverse image of $[M_{I, \infty} \backslash \Gr_{M, I}^n]$ (resp. $[M_{I, \infty} \backslash \Gr_{M, I, \underline{\lambda} }]$) by the morphism $\epsilon_{M, N, I, \infty}$. 
By \ref{subsection-Gr-M-I-n-Gr-M-I-omega}, we have 
\begin{equation}   \label{equation-Cht-M-n-Cht-M-omega}
\Cht_{M, N, I}^n = \bigcup_{ \underline{\lambda}=(\lambda_i)_{i \in I} \in (\wh \Lambda_M^+)^I, \, \langle [\sum_{i \in I} \lambda_i] , 2(\rho_G - \rho_M) \rangle = n    } \Cht_{M, N, I, \underline{\lambda}} .
\end{equation}

\begin{lem}  \label{lem-Cht-M-I-omega-non-empty}  (\cite{var} Proposition 2.16 d), \cite{vincent} Proposition 2.6 e))
$\Cht_{M, N, I, \underline{\lambda}}$ is non-empty if and only if $[\sum_{i \in I} \lambda_i]$ is zero. 
\cqfd
\end{lem}

\sssec{}   \label{subsection-Cht-M-I-n-non-empty}
We deduce from Lemma \ref{lem-Cht-M-I-omega-non-empty} and (\ref{equation-Cht-M-n-Cht-M-omega}) that $\Cht_{M, N, I}^n$ is non-empty if and only if $n=0$.
So the image of $\epsilon_{M, N, I, \infty}$ is included in the open and closed substack $[M_{I, \infty} \backslash \Gr_{M, I}^0 ] $. 

\quad

\sssec{}
Using the notation in diagram (\ref{diagram-Cht-Gr-G-P-M}), we deduce
\begin{equation}    \label{equation-inverse-image-S-M-W-equal-S-M-W-theta-0}
\begin{aligned}
(\epsilon_{M, d})^* \mc S_{M, I, W}^{d, \, \Oe}  & \isom (\epsilon_{M, d})^* \big( \bigoplus_{n \in \Z} (\ov{\pi_{d, n}^0})_! (\ov{i_{d, n}^0})^* \mc S_{G, I, W}^{d, \, \Oe}   [n](n/2) \big) [-2m](-m) \\
&  =  (\epsilon_{M, d} )^*   ( \ov{\pi^0_{d, 0}} )_!  ( \ov{i^0_{d, 0}} )^* \mc S_{G, I, W}^{d, \, \Oe} [-2m](-m)  \\
& = (\epsilon_{M, d} )^*  ( \ov{\pi^0_{d}} )_!  ( \ov{i^0_{d}} )^* \mc S_{G, I, W}^{d, \, \Oe} [-2m](-m)  .
\end{aligned}
\end{equation}
The first isomorphism follows form (\ref{equation-TC-satake-geo-en-quotient}). The second equality follows from \ref{subsection-Cht-M-I-n-non-empty}. Note that by definition, $( \ov{\pi^0_{d}} )_!  ( \ov{i^0_{d}} )^* = \bigoplus_{n \in \Z} (\ov{\pi_{d, n}^0})_! (\ov{i_{d, n}^0})^*$, thus the third equality also follows from \ref{subsection-Cht-M-I-n-non-empty}.

\begin{rem}
(We do not need the following fact.) Let $W \in \on{Rep}_{\Oe}(\wh M^I)$. We can view $W$ as a representation of $\mb G_m$ via $\mb G_m \xrightarrow{2\rho_G - 2\rho_M} Z_{\wh M} \hookrightarrow \wh M \xrightarrow{\triangle} \wh M^I$. We have $$W = \bigoplus_{n \in \Z} W^n$$
where $\mb G_m$ acts on $W^n$ by the character $x \mapsto x^n$. Then we have
$$\Gr_{M, I, W} \cap \Gr_{M, I}^n = \Gr_{M, I, W^n}.$$
\end{rem}

\subsection{Construction of the constant term morphisms}   \label{subsection-CT-complexe}

\sssec{}    \label{subsection-def-pi-i-F-G-to-F-M}
By \cite{LO08b}, we have the notion of six operators for $\Oe$-coefficients. Thus the construction in \cite{cusp-coho} Section 3.3 works in the same way for $\Oe$-coefficients. Note that the morphism $i$ is schematic (\cite{cusp-coho} Remark 3.5.4) and the morphism $\pi$ is of finite type. We have a complex $\pi_! i^* \mc{F}_{G, N, I, W}^{\Oe} \in D_c^{(-)}( \Cht_{M, N, I, \Sigma }, \Oe)$.

Once have (\ref{equation-inverse-image-S-M-W-equal-S-M-W-theta-0}), we use the same argument as in $loc.cit.$ 3.3.1 and 3.3.2 to construct a morphism of complexes in $D_c^{(-)}( \Cht_{M, N, I, \Sigma }, \Oe)$:
\begin{equation}   \label{equation-pi-i-F-G-to-F-M}
\pi_! i^* \mc{F}_{G, N, I, W}^{\Oe} \rightarrow \mc{F}_{M, N, I, W}^{\Oe}.
\end{equation}
where $i$ and $\pi$ are defined in (\ref{equation-Cht-G-I-Cht-P-I-Cht-M-I}). 
Concretely, we use the notations in diagrams (\ref{diagram-Cht-Gr-G-P-M}) and (\ref{diagram-wt-Cht-M}). Morphism (\ref{equation-pi-i-F-G-to-F-M}) is constructed as the composition of morphisms
\begin{equation}
\begin{aligned}
\pi_! i^* \mc{F}_{G, N, I, W}^{\Oe} & = \pi_! i^* (\epsilon_{G, d})^* \mc S_{G, I, W}^{d, \, \Oe} \\
& \simeq \pi_! (\epsilon_{P, d})^* (\ov{i_d^0})^* \mc S_{G, I, W}^{d, \, \Oe} \\
& \simeq (\wt{\pi_d^0})_! (\pi_d)_! (\pi_d)^* (\wt{\epsilon_{M, d}})^* (\ov{i_d^0})^* \mc S_{G, I, W}^{d, \, \Oe} \\
& \xrightarrow{\on{Tr}} (\wt{\pi_d^0})_!  (\wt{\epsilon_{M, d}})^* (\ov{i_d^0})^* \mc S_{G, I, W}^{d, \, \Oe} [-2m](-m) \\
(\text{proper base change}) \quad & \simeq (\epsilon_{M, d})^*  (\ov{\pi_d^0})_! (\ov{i_d^0})^* \mc S_{G, I, W}^{d, \, \Oe} [-2m](-m) \\
(\ref{equation-inverse-image-S-M-W-equal-S-M-W-theta-0}) \quad & \simeq (\epsilon_{M, d})^*  \mc S_{M, I, W}^{d, \, \Oe}  = \mc{F}_{M, N, I, W}^{\Oe}
\end{aligned}
\end{equation}
where we use the fact that $\pi_d$ is smooth of dimension $m$ (cf. \ref{subsection-pi-d-smooth}). The trace map $\on{Tr}$ is the composition $(\pi_d)_! (\pi_d)^*[2m](m) \isom (\pi_d)_! (\pi_d)^! \xrightarrow{\on{counit}} \Id$.

\begin{defi}    \label{def-mc-H-M-prime}
In the same way as in $loc.cit.$ Definition 3.4.2, we define 
$$\Cht_{P, N, I, \Sigma}' := \Cht_{P, N, I, \Sigma} \overset{P(\mc O_N)}{\times} G(\mc O_N), \quad \Cht_{M, N, I, \Sigma}' := \Cht_{M, N, I, \Sigma} \overset{P(\mc O_N)}{\times} G(\mc O_N) ,$$
where $\Cht_{P, N, I, \Sigma}$ and $\Cht_{M, N, I, \Sigma}$ are defined in \ref{subsection-def-Cht-G-Sigma-Cht-M-Sigma}.
\end{defi}

\sssec{}
As in $loc.cit.$ 3.4.3, we construct morphisms
\begin{equation}   \label{equation-Cht-G-I-Cht-P-I-Cht-M-I-Xi}
\Cht_{G, N, I, \Sigma}  / \Xi \xleftarrow{i'} \Cht_{P, N, I, \Sigma}' / \Xi \xrightarrow{\pi'} \Cht_{M, N, I, \Sigma}' / \Xi
\end{equation}

\begin{defi}
As in $loc.cit.$ 3.4.7, we define
$\mc F_{M, N, I, W}^{ \, ' \, \Xi, \, \Oe}$ supported on $\Cht_{M, N, I, \Sigma}'$. We define $\mc H _{M, N, I, W}^{' \, \leq \mu, \, \nu, \, \Oe}$,  
$\restr{ \mc H_{M, N, I, W}^{' \, j, \, \leq \mu, \, \nu, \, \Oe} }{\ov{x}},$
$\restr{ \mc H_{M, N, I, W}^{' \, j, \, \leq \mu, \, \Oe}}{\ov{x}}$ and $\restr{ \mc H_{M, N, I, W}^{' \, j, \, \Oe} }{\ov{x}}$.
\end{defi}

All the constructions in \ref{subsection-some-geo}, \ref{subsection-Satake-geo-para-ind} and \ref{subsection-def-pi-i-F-G-to-F-M} are compatible with the quotient by $\Xi$ and the modification with $'$.

\begin{construction}  \label{constr-funtor-TC-for-complex}
We construct a canonical morphism of complexes in $D_c^{(-)}(\Cht_{M, N, I, \Sigma }^{'} / \Xi, \Oe)$:
\begin{equation}    \label{equation-pi-'-i-'-F-G-to-F-M}
(\pi')_! (i')^* \mc{F}_{G, N, I, W}^{\Oe, \, \Xi} \rightarrow \mc{F}_{M, N, I, W}^{' \, \Oe, \, \Xi}.
\end{equation} 
\end{construction}

\sssec{}
As in $loc.cit.$ 3.5.2, for any $\mu \in \wh{\Lambda}_{G^{\mr{ad}}}^{+, \Q}$ and $\nu \in \wh \Lambda_{Z_M / Z_G}^{\Q}$, we have a commutative diagram
\begin{equation}   \label{diagram-Cht-G-P-M-leq-mu-on-open}
\xymatrixrowsep{1pc}
\xymatrixcolsep{1pc}
\xymatrix{
&   \Cht_{P, N, I, \Sigma }^{' \, \leq \mu, \, \nu} / \Xi  \ar[dl]_{i'}  \ar[dr]^{\pi'}  \ar[dd]^{\mf p_P}
&  \\
 \Cht_{G, N, I, \Sigma }^{\leq \mu} / \Xi  \ar[dr]_{\mf{p}_G}
& 
&  \Cht_{M, N, I, \Sigma }^{' \, \leq \mu, \, \nu} / \Xi  \ar[dl]^{\mf{p}_M} \\
& (X \sm N)^I
&
}
\end{equation}
As in \cite{cusp-coho} Remark 3.5.4, $i'$ is schematic.
By \cite{var} Proposition 5.7 (recalled in \cite{cusp-coho} Proposition 3.5.3 and Remark 3.5.4), there exists an open subscheme $\Omega^{\leq \mu, \, \nu}$ of $(X \sm N)^I$ of the form 
$$\Omega(m) = \{ (x_i)_{i \in I} \in (X \sm N)^I, x_i \neq {}^{\tau^r}x_j \text{ for all } i, j \text{ and } r= 1, 2, \cdots, m\} $$
where ${}^{\tau^r}x$ is the image of $x$ by $\Frob^r: X \rightarrow X$ and $m$ is some positive integer depending on $\mu$ and $\nu$,
such that the restriction $$i': \restr{  \Cht_{P, N, I, \Sigma }^{' \, \leq \mu, \, \nu} / \Xi}{\Omega^{\leq \mu, \, \nu}}  \rightarrow  \restr{ \Cht_{G, N, I, \Sigma }^{\leq \mu} / \Xi}{\Omega^{\leq \mu, \, \nu}}$$ is proper. This is a key ingredient of the following construction.

\begin{construction}  \label{constr-funtor-TC-for-cohomology}
As in \cite{cusp-coho} 3.5.6, the cohomological correspondence for diagram (\ref{diagram-Cht-G-P-M-leq-mu-on-open}) restricted to $\Omega^{\leq \mu, \, \nu}$ and for morphism (\ref{equation-pi-'-i-'-F-G-to-F-M}) gives a morphism of complexes in $D_c^{(-)}(\Omega^{\leq \mu, \, \nu}, \Oe)$:
\begin{equation}     \label{equation-mc-H-G-to-mc-H-M-prime-leq-mu-nu}
\mc C_{G, N}^{P, \, \leq \mu, \, \nu, \, \Oe}: \restr{ \mc H_{G, N, I, W}^{\leq \mu, \, \Oe} }{\Omega^{\leq \mu, \, \nu}} \rightarrow \restr{ {\mc H '}_{\! \! \! M, N, I, W}^{ \; \leq \mu, \, \nu, \, \Oe}   }{\Omega^{\leq \mu, \, \nu}}.
\end{equation}

For any $j \in \Z$, we have a morphism of sheaves:
\begin{equation}     \label{equation-mc-H-G-to-mc-H-M-prime-leq-mu-nu-sheaves}
\mc C_{G, N}^{P, \, j, \, \leq \mu, \, \nu, \, \Oe}: \restr{ \mc H_{G, N, I, W}^{j, \, \leq \mu, \, \Oe} }{\Omega^{\leq \mu, \, \nu}} \rightarrow \restr{ {\mc H '}_{\! \! \! M, N, I, W}^{\; j, \, \leq \mu, \, \nu, \, \Oe}   }{\Omega^{\leq \mu, \, \nu}}.
\end{equation}
\end{construction}

\sssec{}   \label{subsection-ov-x}
For $i \in I$, let $\on{pr}_i: X^I \rightarrow X$ be the projection to the $i$-th factor.
Let $\ov{x}$ be a geometric point of $(X \sm N)^I$ such that for every $i, j \in I$, the image of the composition $$\ov{x} \rightarrow (X \sm N)^I \xrightarrow{(\on{pr}_i, \on{pr}_j) } (X \sm N) \times (X \sm N)$$ 
is not included in the graph of any non-zero power of Frobenius morphism $\on{Frob}: X \sm N \rightarrow X \sm N$.

In particular, when $i=j \in I$, the above condition is equivalent to the condition that the composition $\ov{x} \rightarrow (X \sm N)^I \xrightarrow{\on{pr}_i } X \sm N $ 
is over the generic point $\eta$ of $X \sm N$.

One example of geometric point satisfying the above condition is $\ov{x} = \ov{\eta^I}$, a geometric point over the generic point $\eta^I$ of $X^I$. 
Another example of geometric point satisfying the above condition is $\ov{x}=\Delta(\ov{\eta})$, where $\Delta: X \hookrightarrow X^I$ is the diagonal inclusion and $\ov{\eta}$ is a geometric point over $\eta$.

The reason for which we impose such a condition is that for any $\mu$ and $\nu$, we have $\ov{x} \in \Omega^{\leq \mu, \, \nu}$.

\sssec{}  \label{subsection-C-G-P-leq-mu-Oe}
Restrict (\ref{equation-mc-H-G-to-mc-H-M-prime-leq-mu-nu-sheaves}) to $\ov x$, we have a morphism
\begin{equation}    \label{equation-CT-mu-nu-Oe}
C_{G, N}^{P, \, j, \, \leq \mu, \, \nu, \, \Oe}: \restr{ \mc H_{G, N, I, W}^{j, \, \leq \mu, \, \Oe} }{\ov{x}} \rightarrow \restr{ {\mc H '}_{\! \! \! M, N, I, W}^{ \; j, \, \leq \mu, \, \nu, \, \Oe}   }{\ov{x}}.
\end{equation}
   
Taking into account \ref{subsection-H-M-leq-mu-nu-zero-for-nu-big}, we define a morphism
\begin{equation}   \label{equation-C-G-P-leq-mu-prod-en-nu}
\prod_{\nu \in \wh \Lambda_{Z_M / Z_G}^{\mu}} C_{G, N}^{P, \, j, \, \leq \mu, \, \nu, \, \Oe}: \restr{ \mc H_{G, N, I, W}^{j, \, \leq \mu, \, \Oe} }{\ov{x}} \rightarrow \prod_{\nu \in \wh \Lambda_{Z_M / Z_G}^{\mu}} \restr{ {\mc H '}_{\! \! \! M, N, I, W}^{ \; j, \, \leq \mu, \, \nu, \, \Oe}   }{\ov{x}}.
\end{equation}
By definition, $\restr{ {\mc H '}_{\! \! \! M, N, I, W}^{ \; j, \, \leq \mu, \, \Oe}   }{\ov{x}} = \prod_{\nu \in \wh \Lambda_{Z_M / Z_G}^{\mu}} \restr{ {\mc H '}_{\! \! \! M, N, I, W}^{ \; j, \, \leq \mu, \, \nu, \, \Oe}   }{\ov{x}}.$ To shorten the notations, we also write the morphism (\ref{equation-C-G-P-leq-mu-prod-en-nu}) as
\begin{equation}   \label{equation-C-G-P-leq-mu-prod-en-nu-simplify}
C_{G, N}^{P, \, j, \, \leq \mu, \, \Oe} :   \restr{ \mc H_{G, N, I, W}^{j, \, \leq \mu, \, \Oe} }{\ov{x}} \rightarrow \restr{ {\mc H '}_{\! \! \! M, N, I, W}^{ \; j, \, \leq \mu, \, \Oe}   }{\ov{x}}.
\end{equation}
Lemma A.0.8 of $loc.cit.$ still holds for $\Oe$-coefficients. As in $loc.cit.$ 3.5.9, we take inductive limit on $\mu$:

\begin{defi}   \label{def-CT-cohomology}
For every degree $j \in \Z$, we define the constant term morphism of cohomology groups:
\begin{equation}
 C_{G, N}^{P, \, j, \, \Oe}: \restr{ \mc H _{G, N, I, W}^{j, \, \Oe}  }{ \ov{x}  } \rightarrow   \restr{ \mc  {H'}_{\! \! \! M, N, I, W}^{\; j, \, \Oe} }{ \ov{x}  }  . 
\end{equation}
\end{defi}

\begin{rem} \label{rem-CT-cohomology-nu-Oe}
For any $\nu \in \wh \Lambda_{Z_M / Z_G}^{\Q}$, taking inductive limit on $\mu$ in (\ref{equation-CT-mu-nu-Oe}), we obtain the constant term morphism of cohomology groups:
\begin{equation}
 C_{G, N}^{P, \, j, \, \nu, \, \Oe}: \restr{ \mc H _{G, N, I, W}^{j, \, \Oe}  }{ \ov{x}  } \rightarrow   \restr{ \mc  {H'}_{\! \! \! M, N, I, W}^{\; j, \, \nu, \, \Oe} }{ \ov{x}  }  . 
\end{equation}

As in \cite{cusp-coho} Remark 3.5.11, $\restr{ \mc  {H'}_{\! \! \! M, N, I, W}^{\; j, \, \leq \mu, \, \nu, \, \Oe} }{ \ov{x}  }  \rightarrow \restr{ \mc  {H'}_{\! \! \! M, N, I, W}^{\; j, \, \nu, \, \Oe} }{ \ov{x}  } $ for each $\nu$ induces a morphism
\begin{equation}    \label{equation-H-M-leq-mu-to-prod-H-M-nu-Oe} 
\restr{ \mc  {H'}_{\! \! \! M, N, I, W}^{\; j, \, \leq \mu, \, \Oe} }{ \ov{x}  }  \rightarrow \prod_{\nu \in \wh \Lambda_{Z_M / Z_G}^{\Q} } \restr{ \mc  {H'}_{\! \! \! M, N, I, W}^{\; j, \, \nu, \, \Oe} }{ \ov{x}  } 
\end{equation}
Taking inductive limit on $\mu$, we deduce a morphism
\begin{equation}   \label{equation-H-M-to-prod-H-M-nu-Oe}
\restr{ \mc  {H'}_{\! \! \! M, N, I, W}^{\; j, \, \Oe} }{ \ov{x}  }  \rightarrow \prod_{\nu \in \wh \Lambda_{Z_M / Z_G}^{\Q} } \restr{ \mc  {H'}_{\! \! \! M, N, I, W}^{\; j, \, \nu, \, \Oe} }{ \ov{x}  } 
\end{equation}
The following diagram commutes: 
\begin{equation}   \label{diagram-prod-CT-nu-with-CT}
\xymatrixrowsep{2pc}
\xymatrixcolsep{8pc}
\xymatrix{
\restr{ \mc H _{G, N, I, W}^{j, \, \Oe}  }{ \ov{x}  }   \ar[r]^{ \prod_{\nu} C_{G, N}^{P, \, j, \, \nu, \, \Oe}   }  \ar[rd]_{  C_{G, N}^{P, \, j, \, \Oe}   }
& \prod_{\nu \in \wh \Lambda_{Z_M / Z_G}^{\Q} } \restr{ \mc  {H'}_{\! \! \! M, N, I, W}^{\; j, \, \nu, \, \Oe} }{ \ov{x}  }  \\
&  \restr{ \mc  {H'}_{\! \! \! M, N, I, W}^{\; j, \, \Oe} }{ \ov{x}  }  \ar[u]_{(\ref{equation-H-M-to-prod-H-M-nu-Oe})} . 
}
\end{equation}

In Remark \ref{rem-Ker-prod-nu-TC-equal-Ker-TC} in the next section, we will show that $\prod_{\nu} C_{G, N}^{P, \, j, \, \nu, \, \Oe}$ and $C_{G, N}^{P, \, j, \, \Oe} $ have the same kernel.

\end{rem}

\begin{defi}  \label{def-cusp-coho}
For every degree $j \in \Z$, we define the cuspidal cohomology group:
\begin{equation}
 \big( \restr{ \mc H_{G, N, I, W}^{j, \, \Oe } }{ \ov{x} } \big)^{\on{cusp}}:= \underset{P \varsubsetneq G} \bigcap   \Ker C^{P, \,  j, \, \Oe }_{G, N} .
\end{equation}
This is a sub-$\Oe$-module of $\restr{ \mc H _{G, N, I, W}^{j, \, \Oe} }{\ov{x}}$. 
\end{defi}

\begin{example}   \label{exemple-TC-cht-sans-patte}
(Shtukas without paws)
As in $loc.cit.$ Example 3.5.15, when $I = \emptyset$ (empty) and $W=\bf 1$ (trivial representation), the constant term morphism $C_{G, N}^{P, \, j, \, \Oe}$ coincides (up to constants depending on $\nu$) with the classicial constant term morphism:
$$
\begin{aligned}
C_c(G(F) \backslash G(\mb A) / K_{G, N} \Xi, \Oe ) & \rightarrow C( U(\mb A) M(F) \backslash G(\mb A) / K_{G, N}  \Xi, \Oe ) \\
f \quad \quad & \mapsto \quad  f^P: g \mapsto \sum_{u \in U(F) \backslash  U( \mb A)}   f(ug) .
\end{aligned}
$$  

As a consequence, $\big( \restr{ \mc H_{G, N, \emptyset, \bf 1}^{0, \, \Oe } }{ \Fqbar } \big)^{\on{cusp}} = C_c^{\on{cusp}}(G(F) \backslash G(\mb A) / K_{G, N} \Xi, \Oe )$. 
\end{example}

\begin{rem}   \label{rem-image-of-TC-in-cone}
By \ref{subsection-C-G-P-leq-mu-Oe} and Remark \ref{rem-CT-cohomology-nu-Oe}, one important property of the constant term morphism $C_G^{P, \, j, \, \Oe}$ is that every image is supported on the components $\restr{ \mc  {H'}_{\! \! \! M, N, I, W}^{\; j, \, \nu, \, \Oe} }{ \ov{x}  }$ of $\restr{ \mc  {H'}_{\! \! \! M, N, I, W}^{\; j, \, \Oe} }{ \ov{x}  }$ indexed by the translated cone $\wh \Lambda_{Z_M / Z_G}^{\mu} \subset \wh \Lambda_{Z_M / Z_G}^{\Q}$ for some $\mu \in \wh \Lambda_{G^{\mr{ad}}}^{+, \Q}$. 

When $I=\emptyset$, $W=\bf 1$ and coefficients in $E$, this property was already described by Jonathan Wang in [Wan18] Section 5.1.
\end{rem}

\quad

\section{Finiteness of the cuspidal cohomology}

Let $\ov{x}$ be a geometric point of $(X \sm N)^I$ as in \ref{subsection-ov-x}. 

\subsection{Consequences of the contractibility of deep enough horospheres}  

Let $P$ be a parabolic subgroup of $G$ and $M$ its Levi quotient. 
\cite{cusp-coho} Sections 4.1-4.5 are totally geometric. Section 4.6 in $loc.cit.$ still holds for cohomology groups with $\Oe$-coefficients.

\begin{defi} For any $\mu \in \wh{\Lambda}_{G^{\mr{ad}}}^{+, \Q}$, in $loc.cit.$ Definition 4.1.1, we defined a bounded set 
\begin{equation}
S_M(\mu):=\{\lambda \in \wh{\Lambda}_{G^{\mr{ad}}}^{+, \Q} \, | \, \lambda \leq^{\ov M} \mu \}
\end{equation}
where $\ov M = M / Z_G$ and $\leq^{\ov M}$ is defined by applying \ref{subsection-def-Lambda-G} to $H=\ov M$. In $loc.cit.$ Definition 4.1.10, we defined $\Cht_{G, N, I, W}^{S_M(\mu)} $ and $\Cht_{M, N, I, W}^{' S_M(\mu)} $.
We define complexes over $(X \sm N)^I$
$$\mc H_{G, N, I, W}^{S_M(\mu), \, \Oe} = R(\mf{p}_G)_! ( \restr{\mc{F}_{G, N, I, W}^{\Xi, \, \Oe} } {\Cht_{G, N, I, W}^{S_M(\mu)} / \Xi } )$$
$$\mc H_{M, N, I, W}^{' \, S_M(\mu), \, \Oe} =  R(\mf{p}_M ')_! ( \restr{\mc{F}_{M, N, I, W}^{ ' \,\Xi, \, \Oe} } {\Cht_{M, N, I, W}^{' S_M(\mu)} / \Xi } ).$$
\end{defi}

\sssec{}
\cite{cusp-coho} Sections 4.1-4.5 still hold if we replace everywhere index $W$ by $\Sigma$ defined in \ref{subsection-Gr-G-I-W-included-in-Gr-G-I-omega}. In particular, we have morphisms
$$\restr{ \Cht_{G, N, I, \Sigma}^{S_M(\mu)} / \Xi }{\ov x} \xleftarrow{i^{' \, S_M(\mu)}}  \restr{ \Cht_{P, N, I, \Sigma}^{' \, S_M(\mu)} / \Xi }{\ov x} \xrightarrow{\pi^{' \, S_M(\mu)}} \restr{ \Cht_{M, N, I, \Sigma}^{' \, S_M(\mu)} / \Xi }{\ov x} $$

\sssec{}
Let $\wt C(G, X, N, W)$ as in $loc.cit.$ Definition 4.6.1. 
As in $loc.cit.$ 4.6.3, for $\mu \in \wh{\Lambda}_{G^{\mr{ad}}}^{+, \Q}$ such that $\langle  \mu, \alpha
_i \rangle > \wt C(G, X, N, W) \text{ for all }  i \in \Gamma_G - \Gamma_M$, the morphism $i^{' \, S_M(\mu)}$ is schematic and proper.
Similar to Section 2, we construct a constant term morphism in $D_c^{(-)}((X \sm N)^I, \Oe)$:
\begin{equation}    \label{equation-mc-C-G-P-S-M-mu}
\mc C_{G, N}^{P, \, S_M(\mu), \, \Oe}: \restr{ \mc H_{G, N, I, W}^{S_M(\mu), \, \Oe} }{ \ov x } \rightarrow \restr{ \mc H_{M, N, I, W}^{' \, S_M(\mu), \, \Oe} }{  \ov x } .
\end{equation}

\begin{prop}  \label{prop-TC-isom-mu-grand}
For any $\mu \in \wh{\Lambda}_{G^{\mr{ad}}}^{+, \Q}$ such that $\langle  \mu, \alpha_i \rangle > \wt C(G, X, N, W) \text{ for all }  i \in \Gamma_G - \Gamma_M$, morphism (\ref{equation-mc-C-G-P-S-M-mu}) is an isomorphism.
\end{prop}
\dem
The same as $loc.cit.$ Proposition 4.6.4.
\cqfd

\subsection{Finiteness of the cuspidal cohomology}

The goal of this subsection is to prove:
\begin{thm}  \label{thm-cusp-dim-fini-second}
The $\Oe$-module $\Hcusp$ defined in Definition \ref{def-cusp-coho} is of finite type.
\end{thm}

Theorem \ref{thm-cusp-dim-fini-second} will be a direct consequence of the following proposition.
\begin{prop}  \label{prop-H-cusp-inclus-dans-H-leq-mu}
Let $G, X, N, I, W, \ov{x}$ as before.
There exists $\mu_0 \in \wh{\Lambda}_{G^{\mr{ad}}}^{+, \Q}$ (depending on $G, X, N, W, j, \ov{x}$) such that
$$\Hcusp \subset \mathrm{Im}(  \restr{ \mc H_{G, N, I, W}^{j, \, \leq \mu_0, \, \Oe } }{ \ov{x} }  \rightarrow \restr{ \mc H_{G, N, I, W}^{j, \, \Oe} }{ \ov{x}  } )$$
where $\mc H_{G, N, I, W}^{j, \, \Oe} $ is the inductive limit defined in Definition \ref{def-mc-H-ind-limit}.
\end{prop}

\begin{nota}
In the remaining part of this section, to simplify the notations, we will omit the indices $N, I, W, \ov{x}$. For example, $H_G^{j, \, \Oe} := \restr{ \mc H_{G, N, I, W}^{j, \, \Oe} }{ \ov{x} } $ and $H_M^{j, \, \Oe} := \restr{ \mc H_{M, N, I, W}^{' \, j, \, \Oe} }{ \ov{x} } $.
\end{nota}

\sssec{}
Let $\wh R_{G^{\on{ad}}}$ be the coroot lattice of $G^{\mr{ad}}$. Let $\wh R_{G^{\on{ad}}}^+:= \wh \Lambda_{G^{\mr{ad}}}^+ \cap \wh R_{G^{\on{ad}}}$. Let $r \in \N$ as in $loc.cit.$ 5.1.1. For any $i \in \Gamma_G$, we denote by $\alpha_i$ the corresponding simple root and by
 $\check{\alpha}_i$ the corresponding simple coroot. Let $P_{\alpha_i}$ be the maximal parabolic subgroup with Levi quotient $M_{\alpha_i}$ such that $\Gamma_G - \Gamma_{M_{\alpha_i}} = \{ i \}$.

\sssec{}     \label{subsection-H-leq-mu-to-H-j}

For any $\lambda \in \frac{1}{r} \wh R_{G^{\mr{ad}}}^+$, let $$\mc I_{\lambda}: H_G^{j, \, \leq \lambda, \, \Oe } \rightarrow H_G^{j, \, \Oe}$$ be the morphism to the inductive limit. For any simple root $\alpha$ of $G$, let $$ C_G^{P_\alpha, \, j, \, \Oe} \circ \mc I_{\lambda}:  H_G^{j, \, \leq \lambda, \, \Oe } \rightarrow H_{M_\alpha}^{' \, j, \, \Oe}$$ be the composition of morphisms $H_G^{j, \, \leq \lambda, \, \Oe } \xrightarrow{\mc I_{\lambda} } H_G^{j, \, \Oe} \xrightarrow{ C_G^{P_\alpha, \, j, \, \Oe}} H_{M_\alpha}^{' \, j, \, \Oe}$, where the second morphism is the constant term morphism defined in Definition \ref{def-CT-cohomology}.

\sssec{}
Note that for every $c \in H_{G}^{j, \, \Oe } $, there exists $\lambda \in \wh R_{G^{\mr{ad}}}^+$ large enough such that $c \in \mathrm{Im}(  H_{G}^{j, \, \leq \lambda, \, \Oe }  \rightarrow H_{G}^{j, \, \Oe} )$. Thus Proposition \ref{prop-H-cusp-inclus-dans-H-leq-mu} will be a direct consequence of (b) in the following proposition:

\begin{prop}   \label{prop-finiteness-a-b-c}
There exists a constant $C_G^0 \in \Q^{\geq 0}$ (depending on $G, X, N, W, j, \ov{x}$), such that the following properties hold:

(a) Let $\mu \in  \frac{1}{r} \wh R_{G^{\mr{ad}}}^+ $ such that $\langle \mu , \alpha_i \rangle \geq C_G^0$ for all $i \in \Gamma_G$.
Then for any simple root $\alpha$ of $G$ such that $\mu - \frac{1}{r} \check{\alpha} \in \frac{1}{r} \wh R_{G^{\mr{ad}}}^+$ (which is automatic if $C_G^0 > \frac{2}{r}$), the morphism
\begin{equation}
\Ker( H_G^{j, \, \leq \mu - \frac{1}{r} \check{\alpha}, \, \Oe } \rightarrow H_{M_\alpha}^{' \, j, \, \Oe}) \rightarrow \Ker( H_G^{j, \, \leq \mu, \, \Oe } \rightarrow H_{M_\alpha}^{' \, j, \, \Oe})
\end{equation}
is surjective.

(b) There exists $\mu_0 \in \frac{1}{r} \wh R_{G^{\mr{ad}}}^+ $ (depending on $C_G^0$), such that for any $\lambda \in \frac{1}{r} \wh R_{G^{\mr{ad}}}^+$ satisfying $\lambda \geq \mu_0$ and $\langle \lambda , \alpha_i \rangle \geq C_G^0$ for all $i \in \Gamma_G$, the morphism
\begin{equation}
\Ker( H_G^{j, \, \leq \mu_0, \, \Oe  } \rightarrow \prod_{P \subsetneq G} H_{M}^{' \, j, \, \Oe}) \rightarrow \Ker( H_G^{j, \, \leq \lambda, \, \Oe } \rightarrow \prod_{P \subsetneq G} H_{M}^{' \, j, \, \Oe})
\end{equation}
 is surjective.

(c) There exists a constant $C_G \geq C_G^0$, such that for any $\lambda \in \frac{1}{r} \wh R_{G^{\mr{ad}}}^+ $ satisfying $\langle \lambda , \alpha_i \rangle \geq C_G$ for all $i \in \Gamma_G$, the morphism $\mc I_{\lambda}: H_G^{j, \, \leq \lambda, \, \Oe} \rightarrow H_G^{j, \, \Oe}$ is injective.
\end{prop}

\dem
As in \cite{cusp-coho} Section 5, we use an inductive argument on the semisimple rank of the group $G$. Firstly we prove the statements (a), (b) and (c) for every Levi subgroup of $G$ of rank $0$. This is the same as in $loc.cit.$ Section 5.2.

Secondly we prove the key step: for $n \geq 1$, if (c) is true for all Levi subgroups of rank $n-1$, then (a) is true for all Levi subgroups of rank $n$. Then we deduce easily (a) $\Rightarrow$ (b) for all Levi subgroups of rank $n$. These are the same as in $loc.cit.$ Section 5.3.

Finally we prove that (b) $\Rightarrow$ (c) for all Levi subgroups of rank $n$. This is the same as in $loc.cit.$ Section 5.4, except that we 
replace $loc.cit.$ Lemma 5.4.3 by the following Lemma \ref{lem-limite-stationne-Ker}.
\cqfd

\sssec{}
As in \ref{subsection-H-leq-mu-to-H-j}, for $\mu \in \frac{1}{r} \wh R_{G^{\mr{ad}}}^+ $, we have $$\mc I_{\mu}: H_G^{j, \, \leq \mu, \, \Oe } \rightarrow H_G^{j, \, \Oe} .$$
For $\lambda \in \frac{1}{r} \wh R_{G^{\mr{ad}}}^+ $ such that $\lambda \geq \mu$, we denote by $$\mc I_{\mu}^{\lambda}: H_G^{j, \leq \mu, \, \Oe} \rightarrow H_G^{j, \leq \lambda, \, \Oe}$$ the morphism defined in \ref{subsection-H-G-leq-mu-1-to-leq-mu-2}.
We have $\Ker(\mc I_{\mu}^{\lambda}) \subset \Ker(\mc I_{\mu}) \subset H_G^{j, \leq \mu, \, \Oe}.$

For $\lambda_2 \geq \lambda_1 \geq \mu$, we have
$\Ker(\mc I_{\mu}^{\lambda_1}) \subset \Ker(\mc I_{\mu}^{\lambda_2}).$

\begin{lem} \label{lem-limite-stationne-Ker} (see $loc.cit.$ Lemma 5.4.3 for $\Ql$-coefficient)
Let $\mu \in \frac{1}{r} \wh R_{G^{\mr{ad}}}^+ $. There exists $\mu^{\sharp} \in \wh R_{G^{\mr{ad}}}^+$ such that $\mu^{\sharp} \geq \mu$ and $\Ker(\mc I_{\mu}^{\mu^{\sharp}}) = \Ker(\mc I_{\mu})$.
\end{lem}
\dem
We have the filtered system $\{  \Ker(\mc I_{\mu}^{\lambda}) \; | \; \lambda \in \frac{1}{r} \wh R_{G^{\mr{ad}}}^+ , \lambda \geq \mu   \}$ in $\Ker(\mc I_{\mu})$ and $\Ker(\mc I_{\mu}) =\varinjlim _{\lambda} \Ker(\mc I_{\mu}^{\lambda})$. Since $H_G^{j, \, \leq \mu, \, \Oe }$ is a $\Oe$-module of finite type, it is Noetherian. So the sub-$\Oe$-module $\Ker(\mc I_{\mu})$ is Noetherian. Then the result is clear.
\cqfd

\quad

\begin{rem}   \label{rem-Ker-prod-nu-TC-equal-Ker-TC}
The proof of $loc.cit.$ Lemma 5.3.4 still holds for $\Oe$-coefficients. In particular, there exists $\wt{\mu} \in \wh{\Lambda}_{G^{\mr{ad}}}^{+, \Q}$ such that for all $\mu \geq \wt{\mu}$, morphism
(\ref{equation-H-M-leq-mu-to-prod-H-M-nu-Oe}) defined in Remark \ref{rem-CT-cohomology-nu-Oe} is injective.
Thus morphism (\ref{equation-H-M-to-prod-H-M-nu-Oe}) is injective.
We deduce from diagram (\ref{diagram-prod-CT-nu-with-CT}) that 
\begin{equation}
\on{Ker} ( \prod_{\nu} C_{G, N}^{P, \, j, \, \nu, \, \Oe} ) = \on{Ker} ( C_{G, N}^{P, \, j, \, \Oe} ) .
\end{equation}
Hence
\begin{equation}
\big( \restr{ \mc H_{G, N, I, W}^{j, \, \Oe } }{ \ov{x} } \big)^{\on{cusp}} = \bigcap_{P \varsubsetneq G} \, \bigcap_{\quad \nu \in \wh \Lambda_{Z_M / Z_G}^{\Q}} \on{Ker} C_{G, N}^{P, \, j, \, \nu, \, \Oe} .
\end{equation}
\end{rem}

\quad

\section{Commutativity between constant term morphisms and actions of Hecke algebras}

We recall \cite{cusp-coho} Section 6.1-6.2, with necessary modifications.
In this section, all the stacks are restricted to a geometric point $\ov{x}$ of $(X \sm N)^I$ as in \ref{subsection-ov-x}. 

\subsection{Compatibility of constant term morphisms and level change}

\sssec{}   \label{subsection-K-h-with-niveau}
As in $loc.cit.$ 6.1.1, let $K$ be a compact open subgroup of $G(\mb O)$, let $\wt N$ be a level such that $K_{\wt N} \subset K$.
We define $\Cht_{G, K, I, W}:=\Cht_{G, \wt N, I, W} / (K / K_{\wt N})$. It is independent of the choice of $\wt N$.
We define $\mc F_{G, K, I, W}^{\Xi, \, \Oe}$ over $\Cht_{G, K, I, W} / \Xi $
and $$H_{G, K, I, W}^{j, \, \Oe}:=\varinjlim _{\mu} H_c^{j}(  \Cht_{G, K, I, W}^{\leq \mu} / \Xi , \mc F_{G, K, I, W}^{\Xi, \, \Oe}).$$
When $K=K_N$ for some level $N$, we have $H_{G, K_N, I, W}^{j, \, \Oe} = \restr{ \mc H_{G, N, I, W}^{j, \, \Oe}  }{ \ov{x} }$.

As in $loc.cit.$ 6.1.2, let $K'$ be another compact open subgroup of $G(\mb O)$ and $K' \subset K$. 
We have a morphism $\on{pr}^G_{K', K}: \Cht_{G, K', I, W} \rightarrow \Cht_{G, K, I, W}$. It is finite étale of degree $\sharp(K / K')$ (the cardinality). The adjunction morphism $\Id \rightarrow (\on{pr}^G_{K', K})_*(\on{pr}^G_{K', K})^*$ induces a morphism of cohomology groups 
\begin{equation*}
\on{adj}(\on{pr}^G_{K', K}): H_{G, K, I, W}^{j, \, \Oe} \rightarrow H_{G, K', I, W}^{j, \, \Oe}.
\end{equation*}
The counit morphism (in this case equal to the trace map) $(\on{pr}^G_{K', K})_!(\on{pr}^G_{K', K})^! \rightarrow \Id$ induces a (surjective) morphism of cohomology groups
\begin{equation*}
\on{Co}(\on{pr}^G_{K', K}): H_{G, K', I, W}^{j, \, \Oe} \rightarrow H_{G, K, I, W}^{j, \, \Oe}.
\end{equation*}

Although we do not need it, note that the composition of morphisms $\on{Co}(\on{pr}^G_{K', K}) \circ \on{adj}(\on{pr}^G_{K', K})$ is the multiplication by scalar $\sharp(K / K')$.

\sssec{}   \label{subsection-CT-level-K}
Let $K$ be a compact open subgroup of $G(\mb O)$, as in $loc.cit.$ 6.1.4, we define $\Cht_{P, K, I, W}'$. 

As in $loc.cit.$ 6.1.6, let $\mc D$ be the category of discrete sets $S$ equipped with a continuous action of $P(\mb A)$ with finitely many orbits such that the stabilizer of any point is conjugated to some open subgroup of finite index in $P(\mb O)$. For any $S \in \mc D$, we define functorially the cohomology group $H_{M, S, I, W}^{' \, j, \, \Oe}$ in the following way.

When $S$ has only one orbit, choose a point $s \in S$, let $H$ be the stabilizer of $s$. Then $S = P(\mb A) / H$. Let $R$ be a subgroup of finite index in $H \cap U(\mb A)$ and normal in $H$. In $loc.cit.$ 6.1.5, we defined
$\Cht_{M, \infty, I, W} := \varprojlim _{N} \Cht_{M, N, I, W}$
and equipped it with an action of $M(\mb A)$.
As in $loc.cit.$ 6.1.6, we define the Deligne-Mumford stack $\Cht_{M, \infty, I, W} / (H / R) \Xi$ and the perverse sheaf $\mc F_{M, \infty, I, W}^{\Xi, \, \Oe}$ over this stack. We define $$H_{M, H, R, I, W}^{' \, j, \, \Oe}:=\varinjlim _{\mu} \prod_{\nu} H_c^j(\Cht_{M, \infty, I, W}^{\leq \mu, \, \nu} / (H / R) \Xi, \mc F_{M, \infty, I, W}^{\Xi, \, \Oe}),$$

As in $loc.cit.$ 6.1.6, for $R_1 \subset R_2$ two subgroups of finite index in $H \cap U(\mb A)$ and normal in $H$, the projection $H/R_1 \twoheadrightarrow H/R_2$ induces a morphism
$$\mf q_{R_1, R_2}: \Cht_{M, \infty, I, W} / (H / R_1) \rightarrow \Cht_{M, \infty, I, W} / (H / R_2) $$
which is a gerbe for the finite $q$-group $R_2 / R_1$. Note that the cardinality of $R_2 / R_1$ is invertible in $\Oe$. The counit morphism (which is equal to the trace map because $\mf q_{R_1, R_2}$ is smooth of dimension 0) $\on{Co}(\mf q_{R_1, R_2}) : (\mf q_{R_1, R_2})_!(\mf q_{R_1, R_2})^! \rightarrow \Id$ is an isomorphism in $D_c^{(-)}(\Cht_{M, \infty, I, W} / (H / R_2), \Oe)$. Indeed, just as in $loc.cit.$ 6.1.6, by proper base change and the fact that $\mf q_{R_1, R_2}$ is smooth, we reduce to the case of Lemma \ref{lem-finite-gp-gerbe-counit-is-isom} below with $\Gamma = R_2/ R_1$.
The morphism $\on{Co}(\mf q_{R_1, R_2})$ induces an isomorphism of cohomology groups 
\begin{equation}   \label{equation-H-M-S-R-1-to-R-2}
H_{M, H, R_1, I, W}^{' \, j, \, \Oe} \isom H_{M, H, R_2, I, W}^{' \, j, \, \Oe}
\end{equation}
We define $H_{M, S, I, W}^{' \, j, \, \Oe}$ to be any $H_{M, H, R, I, W}^{' \, j, \, \Oe}$, where we identify $H_{M, H, R_1, I, W}^{' \, j, \, \Oe}$ and $H_{M, H, R_2, I, W}^{' \, j, \, \Oe}$ by (\ref{equation-H-M-S-R-1-to-R-2}). As in $loc.cit.$ 6.1.6, $H_{M, S, I, W}^{' \, j, \, \Oe}$ is independent of the choice of $s \in S$. In fact, let $s_1, s_2$ be two point of $S$ and $H_1$ (resp. $H_2$) be the stabilizer of $s_1$ (resp. $s_2$), then $H_2 = p^{-1}H_1p$ for some $p\in P(\mb A)$. The action of $p$ induces an isomorphism $\Cht_{M, \infty, I, W} / (H_1 / R)  \isom \Cht_{M, \infty, I, W} / (p^{-1}H_1p / p^{-1}Rp) $. We deduce an isomorphism of cohomology groups by the adjunction morphism.

In general, for $S = \sqcup_{\alpha \in A} \alpha$ a finite union of orbits, we define $$H_{M, S, I, W}^{' \, j, \, \Oe} := \bigoplus_{\alpha \in A }  H_{M, \alpha, I, W}^{' \, j, \, \Oe} .$$

When $S = G(\mb A) / K$ for some compact open subgroup $K$ in $G(\mb O)$, we define 
\begin{equation}   \label{equation-def-H-M-K-I-W-prime}
H_{M, K, I, W}^{' \, j, \, \Oe}:=H_{M, S, I, W}^{' \, j, \, \Oe} .
\end{equation}

\begin{lem}   \label{lem-finite-gp-gerbe-counit-is-isom}
Let $\Gamma$ be a finite group with cardinality invertible in $\Oe$. We denote by $B\Gamma$ the classifying stack of $\Gamma$ over $\on{Spec} k$, where $k$ is an algebraically closed field over $\Fq$. Let $\mf q: B\Gamma \rightarrow \on{Spec} k$ be the structure morphism. Then the counit morphism (equal to the trace map) $\on{Co}(\mf q): q_! q^! \rightarrow \Id$ of functors in $D_c^{(-)}(\on{Spec} k, \Oe)$ is an isomorphism.
\end{lem}
\dem The same as the proof of \cite{cusp-coho} Lemma 6.1.7.
\cqfd

\begin{rem}
In Lemma \ref{lem-finite-gp-gerbe-counit-is-isom}, the condition on the cardinality of $\Gamma$ is necessary. See \cite{LO08a} 4.9.2 for a counter-example.
\end{rem}

\sssec{}
\cite{cusp-coho} 6.1.8-6.1.14 are still true for cohomology with $\Oe$-coefficients. In particular, let $K$ be a compact open subgroup of $G(\mb O)$. As in $loc.cit.$ 6.1.9, we construct the constant term morphism
\begin{equation}   \label{equation-def-C-G-K-P}
C_{G, K}^{P, \, j, \, \Oe}: H_{G, K, I, W}^{j, \, \Oe} \rightarrow H_{M, K, I, W}^{' \, j, \, \Oe}.
\end{equation}

\sssec{}   \label{subsection-functoriality-S-1-S-2}
Let $S_1, S_2 \in \mc D$ and $f: S_1 \rightarrow S_2$ be a morphism in $\mc D$. 
In $loc.cit.$ 6.1.10, we defined a finite étale morphism
$$\mf q_f^M: \Cht_{M, S_1, R, I, W}' \rightarrow \Cht_{M, S_2, R, I, W}' .$$

The adjunction morphism $\Id \rightarrow (\mf q_f^M)_*(\mf q_f^M)^*$ induces a morphism
$$\on{adj}(\mf q_f^M): H _{M, S_2, I, W}^{' \, j, \, \Oe}  \rightarrow H _{M, S_1, I, W}^{' \, j, \, \Oe}$$
The counit morphism $(\mf q_f^M)_!(\mf q_f^M)^! \rightarrow \Id$ induces a morphism
$$\on{Co}(\mf q_f^M): H _{M, S_1, I, W}^{' \, j, \, \Oe}  \rightarrow H _{M, S_2, I, W}^{' \, j, \, \Oe}$$

\sssec{}   \label{subsection-pr-M-K'-K}
As in $loc.cit.$ 6.1.14, let $K' \subset K$ be two compact open subgroups of $G(\mb O)$. Applying \ref{subsection-functoriality-S-1-S-2} to $S_1 = G(\mb A) / K'$, $S_2 = G(\mb A) / K$ and the projection $G(\mb A) / K' \twoheadrightarrow G(\mb A) / K$, we define a finite étale morphism 
$$\on{pr}^M_{K', K}: \Cht_{M, S_1, R, I, W}'  \rightarrow \Cht_{M, S_2, R, I, W}' .$$
The adjunction morphism $\Id \rightarrow (\on{pr}^M_{K', K})_*(\on{pr}^M_{K', K})^*$ induces
\begin{equation*}
\on{adj}(\on{pr}^M_{K', K}): H _{M, K, I, W}^{' \, j, \, \Oe}  \rightarrow H _{M, K', I, W}^{' \, j, \, \Oe} .
\end{equation*}
The counit morphism $(\on{pr}^M_{K', K})_!(\on{pr}^M_{K', K})^! \rightarrow \Id$ induces
\begin{equation*}
\on{Co}(\on{pr}^M_{K', K}): H _{M, K', I, W}^{' \, j, \, \Oe}  \rightarrow H _{M, K, I, W}^{' \, j, \, \Oe} .
\end{equation*}

Although we do not need it, we can prove that the composition of morphisms $\on{Co}(\on{pr}^M_{K', K}) \circ \on{adj}(\on{pr}^M_{K', K})$ is the multiplication by scalar $\sharp(K / K')$.

\begin{lem}   \label{lem-CT-commutes-with-level-augement}
For $K' \subset K$ as in \ref{subsection-pr-M-K'-K}, the following diagram of cohomology groups commutes:
\begin{equation}     \label{equation-CT-commutes-with-adj}
\xymatrixrowsep{2pc}
\xymatrixcolsep{5pc}
\xymatrix{
  H _{G, K, I, W}^{j, \, \Oe}    \ar[r]^{\on{adj}(\on{pr}^G_{K', K})}  \ar[d]_{  C_{G, K}^{P, \, j, \, \Oe}} 
&   H _{G, K', I, W}^{j, \, \Oe}  \ar[d]^{  C_{G, K'}^{P, \, j, \, \Oe}}  \\
 H _{M, K, I, W}^{' \, j, \, \Oe}  \ar[r]^{\on{adj}(\on{pr}^M_{K', K})}
&  H _{M, K', I, W}^{' \, j, \, \Oe} 
}
\end{equation}
\end{lem}
\dem The same as the proof of $loc.cit.$ Lemma 6.1.15.
\cqfd

\begin{lem}    \label{lem-CT-commutes-with-level-dimunit}
For $K' \subset K$ as in \ref{subsection-pr-M-K'-K}, the following diagram of cohomology groups commutes:
\begin{equation}
\xymatrixrowsep{2pc}
\xymatrixcolsep{5pc}
\xymatrix{
  H _{G, K', I, W}^{j, \, \Oe}    \ar[r]^{\on{Co}(\on{pr}^G_{K', K})}  \ar[d]_{  C_{G, K'}^{P, \, j, \, \Oe}} 
&   H _{G, K, I, W}^{j, \, \Oe}  \ar[d]^{  C_{G, K}^{P, \, j, \, \Oe}}  \\
 H _{M, K', I, W}^{' \, j, \, \Oe}  \ar[r]^{\on{Co}(\on{pr}^M_{K', K}) }
&  H _{M, K, I, W}^{' \, j, \, \Oe} 
}
\end{equation}
\end{lem}
\dem The same as the proof of $loc.cit.$ Lemma 6.1.16.
\cqfd

\subsection{Action of Hecke algebras}

\sssec{}   \label{subsection-g-act-on-H-G-K-infty}
In \cite{cusp-coho} 6.1.3, we defined $\Cht_{G, \infty, I, W} := \varprojlim \Cht_{G, N, I, W} $ and equipped it with an action of $G(\mb A)$.

Let $v$ be a place in $X$. Let $g \in G(F_v)$. Let $\wt K \subset G(\mb O)$ be a compact open subgroup such that $g^{-1} \wt K g \subset G(\mb O)$. In $loc.cit.$ 6.2.1, we proved that the action of $g$ induces an isomorphism
$$\Cht_{G, \infty, I, W} / \wt K \isom \Cht_{G, \infty, I, W} / g^{-1} \wt K g .$$
It induces (by adjunction) an isomorphism of cohomology groups:
\begin{equation}   \label{equation-adj-g-H-G-wt-K}
\on{adj}(g): H_{G, g^{-1} \wt K g, I, W}^{j, \, \Oe} \isom H_{G, \wt K, I, W}^{j, \, \Oe} 
\end{equation}

\sssec{}   \label{subsection-Hecke-G-acts-on-H-G}  
We denote by $\mb O^v$ the ring of integral adèles outside $v$.
Let $K = K^v  K_v  \subset G(\mb O^v) G(\mc O_v) =  G(\mb O) $ be an open compact subgroup.
Let $h =  {\bf 1}_{K_v g K_v} \in C_c(K_v \backslash G(F_v) / K_v, \Oe)$ be the characteristic function of $K_v g K_v$ for some $g \in G(F_v)$. 
As in $loc.cit.$ 6.2.2, the action of $h$ on $H _{G, K, I, W}^{j, \, \Oe}$ is given by the following composition of morphisms
\begin{equation}     \label{equation-T-h-acts-on-H-G}
T(h): H _{G, K, I, W}^{j, \, \Oe} \xrightarrow{\on{adj}} H _{G, K \cap g^{-1} K g, I, W}^{j, \, \Oe} \underset{\sim}{\xrightarrow{  \on{adj}(g)  }} H _{G, g Kg^{-1} \cap K, I, W}^{j, \, \Oe}  \xrightarrow{\on{Co}} H _{G, K, I, W}^{j, \, \Oe} ,
\end{equation}
where $\on{adj}=\on{adj}(\on{pr}^G_{K \cap g^{-1}Kg, K})$ and $\on{Co}=\on{Co}(\pr^G_{gKg^{-1} \cap K, K})$, the isomorphism $\on{adj}(g)$ is induced by (\ref{equation-adj-g-H-G-wt-K}) applied to $\wt K = g Kg^{-1} \cap K$.
Note that (\ref{equation-T-h-acts-on-H-G}) depends only on the class $K_v g K_v$ of $g$ in $G(F_v)$. 
The action of $T(h)$ is equivalent to the one constructed by Hecke correspondence (see \cite{vincent} 2.20 and 4.4).

\sssec{}   \label{subsection-functoriality-K-gKg-M}
Let $\wt K$ and $g$ as in \ref{subsection-g-act-on-H-G-K-infty}.
The right action of $g$ (by right multiplication by $g$) on $G(\mb A)$ induces an isomorphism
\begin{equation}   \label{equation-G-A-K-to-G-A-gKg}
g: G(\mb A) / \wt K \isom G(\mb A) / g^{-1}\wt Kg
\end{equation}
Applying \ref{subsection-functoriality-S-1-S-2} to $S_1=G(\mb A) / \wt K$, $S_2 = G(\mb A) / g^{-1}\wt Kg$ and the isomorphism (\ref{equation-G-A-K-to-G-A-gKg}), we deduce an isomorphism of cohomology groups:
\begin{equation}  \label{equation-adj-g-H-M-wt-K}
\on{adj}(g): H_{M, g^{-1}\wt Kg, I, W}^{' \, j, \, \Oe} \isom H_{M, \wt K, I, W}^{' \, j, \, \Oe}
\end{equation}

\sssec{}   \label{subsection-Hecke-G-acts-on-H-M-'}  
Let $K$ and $h$ as in \ref{subsection-Hecke-G-acts-on-H-G}.  
The action of $h$ on $H _{M, K, I, W}^{' \, j, \, \Oe} $ is given by the following composition of morphisms
\begin{equation}     \label{equation-T-h-acts-on-H-M-'}
T(h): H _{M, K, I, W}^{' \, j, \, \Oe} \xrightarrow{\on{adj}} H _{M, K \cap g^{-1} K g, I, W}^{' \, j, \, \Oe} \underset{\sim}{\xrightarrow{  \on{adj}(g)  }} H _{M, g Kg^{-1} \cap K, I, W}^{' \, j, \, \Oe}  \xrightarrow{\on{Co}} H _{M, K, I, W}^{' \, j, \, \Oe} ,
\end{equation}
where $\on{adj}=\on{adj}(\on{pr}^M_{K \cap g^{-1}Kg, K})$ and $\on{Co}=\on{Co}(\pr^M_{gKg^{-1} \cap K, K})$, the isomorphism $\on{adj}(g)$ is induced by (\ref{equation-adj-g-H-M-wt-K}) applied to $\wt K = g Kg^{-1} \cap K$. Note that (\ref{equation-T-h-acts-on-H-M-'}) depends only on the class $K_v g K_v$ of $g$ in $G(F_v)$. 

We remark that in general $\wt K$ and $g^{-1}\wt Kg$ are not normal in $G(\mb O)$. This is the reason why we define morphism (\ref{equation-adj-g-H-M-wt-K}) using the functoriality in \ref{subsection-functoriality-S-1-S-2}.

\begin{lem}   \label{lem-TC-commute-with-action-g}
Let $\wt K$ and $g$ as in \ref{subsection-g-act-on-H-G-K-infty}. The following diagram of cohomology groups commutes:
\begin{equation}   \label{diagram-TC-commute-avec-g}
\xymatrixrowsep{3pc}
\xymatrixcolsep{4pc}
\xymatrix{
H_{G, g^{-1}\wt Kg, I, W}^{ j, \, \Oe}    \ar[d]_{ C_{G, g^{-1}\wt Kg}^{P, \, j, \, \Oe}} \ar[r]^{\on{adj}(g)}_{\simeq} 
& H_{G, \wt K, I, W}^{j, \, \Oe}   \ar[d]^{ C_{G, \wt K}^{P, \, j, \, \Oe}}   \\
H_{M, g^{-1}\wt Kg, I, W}^{' \, j, \, \Oe} \ar[r]^{\on{adj}(g)}_{\simeq}  
& H_{M, \wt K, I, W}^{' \, j, \, \Oe}    
}
\end{equation}
\end{lem}
\dem
The same as the proof of $loc.cit.$ Lemma 6.2.5.
\cqfd

\begin{prop}   \label{prop-TC-commute-with-Hecke-with-niveau}
(Hecke action at any place commutes with constant term morphism)
For any place $v$ of $X$, any $K$ and $h \in C_c(K_v \backslash G(F_v) / K_v, \Oe)$ as in \ref{subsection-Hecke-G-acts-on-H-G}, the following diagram of cohomology groups commutes:
\begin{equation}   \label{diagram-TC-commute-avec-Hecke}
\xymatrixrowsep{2pc}
\xymatrixcolsep{3pc}
\xymatrix{
  H_{G, K, I, W}^{j, \, \Oe}   \ar[r]^{T(h)}  \ar[d]_{ C_{G, K}^{P, \, j, \, \Oe}} 
&   H_{G, K, I, W}^{j, \, \Oe}  \ar[d]^{ C_{G, K}^{P, \, j, \, \Oe}}  \\
 H_{M, K, I, W}^{' \, j, \, \Oe}  \ar[r]^{T(h)}
&    H_{M, K, I, W}^{' \, j, \, \Oe} 
}
\end{equation}
where the horizontal morphisms are defined in \ref{subsection-Hecke-G-acts-on-H-G} and \ref{subsection-Hecke-G-acts-on-H-M-'},  
the vertical morphisms are the constant term morphism defined in (\ref{equation-def-C-G-K-P}).
\end{prop}
\dem
It is a consequence of Lemma \ref{lem-CT-commutes-with-level-augement}, Lemma \ref{lem-CT-commutes-with-level-dimunit} and Lemma \ref{lem-TC-commute-with-action-g}.
\cqfd

\quad

\sssec{}   \label{subsection-K-h-without-niveau}
From now on let $N \subset X$ be a closed subscheme and $v$ be a place in $X \sm N$. 
As in $loc.cit.$ 6.2.7, we have the (unnormalized) Satake transform:
\begin{equation}   \label{equation-TC-pour-op-de-Hecke}
\begin{aligned}
C_{c}(G(\mc O_v)\backslash G(F_v)/G(\mc O_v), \Oe) & \hookrightarrow C_{c}(M(\mc O_v)\backslash M(F_v) / M(\mc O_v), \Oe) \\
h \quad & \mapsto \quad  h^M: m \mapsto  \sum_{ u \in U( F_v) / U(\mc O_v) }  h(mu).
\end{aligned}
\end{equation}

\sssec{}   \label{subsection-T-h-M-acts-on-H-M-N-I-W}
Applying \ref{subsection-Hecke-G-acts-on-H-G} to $M$ (replacing $G$ by $M$ and $h$ by $h^M$), we obtain an action of $T(h^M)$ on $H _{M, N, I, W}^{j, \, \Oe}$. Since $\Cht_{M, N, I, W} ' = \Cht_{M, N, I, W} \overset{P(\mc O_N)}{\times} G(\mc O_N)$ is a finite union of $ \Cht_{M, N, I, W} $, we deduce an action of $T(h^M)$ on $H _{M, N, I, W}^{' \, j, \, \Oe}$.

\begin{lem}  \label{lem-h-act-on-H-M-equal-h-M-act-on-H-M}
The action of $T(h)$ on $ H _{M, N, I, W}^{' \, j, \, \Oe} = H _{M, K_N, I, W}^{' \, j, \, \Oe} $ (defined in (\ref{equation-T-h-acts-on-H-M-'})) coincides with the action of $T(h^M)$ on $ H _{M, N, I, W}^{' \, j, \, \Oe} $ (defined in \ref{subsection-T-h-M-acts-on-H-M-N-I-W}). 
\end{lem}

\begin{rem}
In \cite{cusp-coho}, the $E$-coefficients version of this lemma is proved by Lemma 6.2.10 in $loc.cit.$. However, the proof of Lemma 6.2.10 in $loc.cit.$ does not work for $\Oe$-coefficients because we no longer have Haar measures. We will give another proof for Lemma \ref{lem-h-act-on-H-M-equal-h-M-act-on-H-M}.
\end{rem}

Since $v \in X \sm N$, $N$ does not play any role in the action of $T(h)$ and $T(h^M)$. To simplify the notations, we give the proof for $N = \emptyset$. The proof for the general case is similar.

\sssec{}   \label{subsection-h-M-scts-on H-M-I-W}
Let $g \in G(F_v)$. Let $h =  {\bf 1}_{G(\mc O_v) g G(\mc O_v)} \in C_c(G(\mc O_v) \backslash G(F_v) / G(\mc O_v), \Oe)$.
We have a correspondence:
$$G(F_v) / G(\mc O_v) \xleftarrow{i^0} P(F_v) / P(\mc O_v) \xrightarrow{\pi^0} M(F_v) / M(\mc O_v)$$
Note that $i^0$ is a bijection.

We view $h$ as a $G(\mc O_v) $-equivariant function on $G(F_v) / G(\mc O_v)$. Thus $h^M$ defined in \ref{subsection-K-h-without-niveau} is given by $(\pi_0)_! (i^0)^{*}h$. Concretely, $h$ is supported on the $G(\mc O_v)$-orbit $G(\mc O_v) g G(\mc O_v) / G(\mc O_v)$. 
The inverse image $(i^0)^{-1} ( G(\mc O_v) g G(\mc O_v) / G(\mc O_v) )$ 
 is a union of $P(\mc O_v)$-orbits 
$$\bigsqcup_{\alpha}  P(\mc O_v) g_P^{\alpha} P(\mc O_v) / P(\mc O_v) $$
where $g_P^{\alpha} \in P(F_v)$ and $ \alpha  \in \{ P(\mc O_v) \text{-orbits in } (i^0)^{-1} ( G(\mc O_v) g G(\mc O_v) / G(\mc O_v) ) \}  $.
For every $\alpha$, the image $\pi^0( P(\mc O_v) g_P^{\alpha} P(\mc O_v) / P(\mc O_v) )$ is a $M(\mc O_v)$-orbit 
$$ M(\mc O_v) g_M^{\alpha} M(\mc O_v) / M(\mc O_v) $$
where $g_M^{\alpha} = \pi^0(g_P^{\alpha}) \in M(F_v)$.
Let $c^{\alpha}$ be the cardinality of the fiber of the projection $P(\mc O_v) g_P^{\alpha} P(\mc O_v) / P(\mc O_v) \xrightarrow{\pi^0} M(\mc O_v) g_M^{\alpha} M(\mc O_v) / M(\mc O_v)$.
We have $$h^M = \sum_{\alpha} c^{\alpha} {\bf 1}_{M(\mc O_v) g_M^{\alpha} M(\mc O_v)} \in C_c(M(\mc O_v) \backslash M(F_v) / M(\mc O_v), \Oe) .$$
Note that for $\alpha_1 \neq \alpha_2$, we have $g_P^{\alpha_1} \neq g_P^{\alpha_2}$, but we may have $g_M^{\alpha_1} = g_M^{\alpha_2}$.

The action of $T(h^M)$ on $ H_{M, I, W}^{j, \, \Oe} $ is given by $\sum_{\alpha} c^{\alpha} T(  {\bf 1}_{M(\mc O_v) g_M^{\alpha} M(\mc O_v)}   )$. For each $\alpha$, apply \ref{subsection-g-act-on-H-G-K-infty} and \ref{subsection-Hecke-G-acts-on-H-G} to $M$ (replacing $G$ by $M$ and $g$ by $g_M^{\alpha}$). The action of $T(  {\bf 1}_{M(\mc O_v) g_M^{\alpha} M(\mc O_v)}   )$ is given by
\begin{equation}   \label{equation-g-M-alpha-act-on-H-M}
H _{M, I, W}^{j, \, \Oe} \xrightarrow{\on{adj}} H _{M, M(\mb O) \cap (g_M^{\alpha})^{-1} M(\mb O) g_M^{\alpha}, I, W}^{j, \, \Oe} \underset{\sim}{\xrightarrow{  \on{adj}(g_M^{\alpha})  }} H _{M, g_M^{\alpha} M(\mb O) (g_M^{\alpha})^{-1} \cap M(\mb O), I, W}^{j, \, \Oe}  \xrightarrow{\on{Co}} H _{M, I, W}^{j, \, \Oe} 
\end{equation}
The middle morphism is induced by 
$$\Cht_{M, \infty, I, W} / g_M^{\alpha} M(\mb O) (g_M^{\alpha})^{-1} \cap M(\mb O) \underset{\sim}{\xrightarrow{  g_M^{\alpha} } } \Cht_{M, \infty, I, W} / M(\mb O) \cap (g_M^{\alpha})^{-1} M(\mb O) g_M^{\alpha} $$

\sssec{}
We view $g \in G(F_v)$ in \ref{subsection-h-M-scts-on H-M-I-W} as an element in $G(\mb A)$. 
By \ref{subsection-Hecke-G-acts-on-H-M-'}, the action of $T(h)$ on $ H_{M, I, W}^{j, \, \Oe} $ is given by
\begin{equation}   \label{equation-T-h-act-on-H-M-I-W}
H _{M, I, W}^{j, \, \Oe} \xrightarrow{\on{adj}} H _{M, G(\mb O) \cap g^{-1} G(\mb O) g, I, W}^{' \, j, \, \Oe} \underset{\sim}{\xrightarrow{  \on{adj}(g)  }} H _{M, g G(\mb O) g^{-1} \cap G(\mb O), I, W}^{' \, j, \, \Oe}  \xrightarrow{\on{Co}} H _{M, I, W}^{j, \, \Oe} 
\end{equation}

The first morphism is induced (by the functoriality in \ref{subsection-functoriality-S-1-S-2}) by the projection
\begin{equation}   \label{equation-G-A-quotient-K-to-G-A-quotient-G-O}
G(\mb A) / G(\mb O) \cap g^{-1} G(\mb O) g \twoheadrightarrow G(\mb A) / G(\mb O)
\end{equation}
The middle morphism is induced by the right multiplication by $g$:
\begin{equation}   \label{equation-G-A-quotient-K-isom-G-A-quotient-gKg}
G(\mb A) /g G(\mb O) g^{-1} \cap G(\mb O) \underset{\sim}{\xrightarrow{  g  }} G(\mb A) / G(\mb O) \cap g^{-1} G(\mb O) g
\end{equation}
The third morphism is induced by the projection 
\begin{equation}    \label{equation-G-A-quotient-gKg-1-to-G-A-quotient-G-O}
G(\mb A) / g G(\mb O) g^{-1} \cap G(\mb O) \twoheadrightarrow G(\mb A) / G(\mb O)
\end{equation}
The $P(\mb A)$ action on all the sets is by left multiplication. 
We have
$$
\begin{aligned}
\Omega:= & \{ P(\mb O) \text{-orbit in } G(\mb O) g G(\mb O) / G(\mb O) \}  \\
\overset{(a)}{=} &  \{ P(\mb O) \text{-orbit in } G(\mb O) / g G(\mb O) g^{-1} \cap G(\mb O) \} \\
\overset{(b)}{=} &  \{ P(\mb A) \text{-orbit in } G(\mb A) / g G(\mb O) g^{-1} \cap G(\mb O) \} \\
\overset{(c)}{=} &  \{ P(\mb A) \text{-orbit in } G(\mb A) / G(\mb O) \cap g^{-1} G(\mb O) g \}  \\
\overset{(d)}{=} &  \{ P(\mb O) \text{-orbit in } G(\mb O) / G(\mb O) \cap g^{-1} G(\mb O) g \} \\
\overset{(e)}{=} &  \{ P(\mb O) \text{-orbit in } G(\mb O) g^{-1} G(\mb O) / G(\mb O) \} 
\end{aligned}
$$
where (a) follows from the fact that $G(\mb O) g G(\mb O) / G(\mb O) = G(\mb O) / \on{Stab} (g)$ and the stabilizer $\on{Stab} (g) = g G(\mb O) g^{-1} \cap G(\mb O)$, (b) and (d) follows from $G(\mb A) / G(\mb O) = P(\mb A) / P(\mb O)$, (c) follows from (\ref{equation-G-A-quotient-K-isom-G-A-quotient-gKg}).

Note that since $g \in G(F_v)$, we have $G(\mb O) g^{-1} G(\mb O) / G(\mb O) = G(\mc O_v) g^{-1} G(\mc O_v) / G(\mc O_v)$.

\sssec{}
Let $S_1 = G(\mb A) / G(\mb O) \cap g^{-1} G(\mb O) g$ and $S_2 = G(\mb A) / G(\mb O)$. Let $f$ be the morphism (\ref{equation-G-A-quotient-K-to-G-A-quotient-G-O}).
In the setting of \cite{cusp-coho} 6.1.10 and Remark 6.1.11 (recalled in  \ref{subsection-functoriality-S-1-S-2} above), all $P(\mb A)$-orbits in $S_1$ are sent to the unique $P(\mb A)$-orbit in $S_2$. Choose a representative $s_2$ in $G(\mb O) / G(\mb O)$ in $S_2$. Its stabilizer $H_2 = P(\mb O)$. 
For each $\alpha$, choose the representative $s_1^{\alpha} \in G(\mb O) / G(\mb O) \cap g^{-1} G(\mb O) g = G(\mb O) g^{-1} G(\mb O) / G(\mb O)$ to be $(g_P^{\alpha})^{-1}$. Then its stabilizer $H_1^{\alpha}=P(\mb O) \cap (g_P^{\alpha})^{-1} P(\mb O) g_P^{\alpha}$. 
The projection $P(\mb A) \twoheadrightarrow M(\mb A)$ induces a morphism
\begin{equation}   \label{equation-g-1Pg-to-g-1Mg}
P(\mb O) \cap (g_P^{\alpha})^{-1} P(\mb O) g_P^{\alpha} \rightarrow M(\mb O) \cap (g_M^{\alpha})^{-1} M(\mb O) g_M^{\alpha} 
\end{equation}
Let $R_1^{\alpha}$ be the kernel of (\ref{equation-g-1Pg-to-g-1Mg}). We have $R_1^{\alpha} = U(\mb O) \cap H_1^{\alpha}$. 
The morphism $H_1^{\alpha} / R_1^{\alpha} \hookrightarrow H_2 / R_1^{\alpha}$ induces a morphism
\begin{equation}   
\Cht_{M, \infty, I, W} / (H_1^{\alpha} / R_1^{\alpha}) \rightarrow \Cht_{M, \infty, I, W} / (H_2 / R_1^{\alpha})
\end{equation}
By adjunction, it induces
\begin{equation}
H_{M, H_2, R_1^{\alpha}, I, W}^{' \, j} \xrightarrow{ \on{adj} } H_{M, H_1^{\alpha}, R_1^{\alpha}, I, W}^{' \, j}
\end{equation}

Moreover, 
$H_1^{\alpha} / R_1^{\alpha} \subset M(\mb O) \cap (g_M^{\alpha})^{-1} M(\mb O) g_M^{\alpha}$ induces morphism
$$ \Cht_{M, \infty, I, W} / (H_1^{\alpha} / R_1^{\alpha})  \rightarrow \Cht_{M, \infty, I, W} / \big( M(\mb O) \cap (g_M^{\alpha})^{-1} M(\mb O) g_M^{\alpha} \big) .$$
We have morphisms
\begin{equation}   \label{equation-g-1Mg-adj-co} 
H_{M, M(\mb O) \cap (g_M^{\alpha})^{-1} M(\mb O) g_M^{\alpha}, I, W}^{ j}  \xrightarrow{\on{adj}} H_{M, H_1^{\alpha}, R_1^{\alpha}, I, W}^{' \, j} \xrightarrow{\on{Co}} H_{M, M(\mb O) \cap (g_M^{\alpha})^{-1} M(\mb O) g_M^{\alpha}, I, W}^{ j} .
\end{equation}
The composition is the multiplication by $\sharp \big( (  M(\mb O) \cap (g_M^{\alpha})^{-1} M(\mb O) g_M^{\alpha} ) / (  H_1^{\alpha} / R_1^{\alpha}   )  \big)$.

Note that $\Cht_{M, I, W} = \Cht_{M, \infty, I, W} / ( P(\mb O) / U(\mb O) ).$
Since $R_1^{\alpha} \subset U(\mb O)$, as in \ref{subsection-CT-level-K}, the projection $P(\mb O) / R_1^{\alpha} \twoheadrightarrow P(\mb O) / U(\mb O)$ induces a morphism 
$$\Cht_{M, \infty, I, W} / (P(\mb O) / R_1^{\alpha}) \rightarrow \Cht_{M, \infty, I, W} / (P(\mb O) / U(\mb O)) $$
which is a gerbe for the finite $q$-group $U(\mb O) / R_1^{\alpha}$. 
We deduce
\begin{equation}   \label{equation-R-adj-co}
H_{M, I, W}^{j} = H_{M, H_2, U(\mb O), I, W}^{' \, j} \underset{\sim}{\xrightarrow{  \on{adj}  } } H_{M, H_2, R_1^{\alpha}, I, W}^{' \, j} \underset{\sim}{\xrightarrow{  \on{Co}  } }  H_{M, H_2, U(\mb O), I, W}^{' \, j} = H_{M, I, W}^{j} 
\end{equation}
where the composition is the multiplication by $(\sharp U(\mb O) / R_1^{\alpha} )^{-1}$.

We have a commutative diagram
\begin{equation}  \label{equation-commute-adj}
\xymatrix{
H_{M, H_2, R_1^{\alpha}, I, W}^{' \, j}  \ar[r]^{\on{adj}}
& H_{M, H_1^{\alpha}, R_1^{\alpha}, I, W}^{' \, j}   \\
H_{M, I, W}^{j} \ar[r]^{\on{adj} \quad \quad \quad }   \ar[u]^{\on{adj} }_{\simeq}
& H_{M, M(\mb O) \cap (g_M^{\alpha})^{-1} M(\mb O) g_M^{\alpha}, I, W}^{ j}    \ar[u]^{\on{adj} }
}
\end{equation}

\sssec{}

Similarly, let $S_3 = G(\mb A) /g G(\mb O) g^{-1} \cap G(\mb O)$. Apply \cite{cusp-coho} 6.1.10 to (\ref{equation-G-A-quotient-gKg-1-to-G-A-quotient-G-O}): $S_3 \rightarrow S_2$. All $P(\mb A)$-orbits in $S_3$ are sent to the unique $P(\mb A)$-orbit in $S_2$. For each $\alpha$, 
choose the representative $s_3^{\alpha} \in G(\mb O) /g G(\mb O) g^{-1} \cap G(\mb O) = G(\mb O) g G(\mb O) / G(\mb O)$ to be $g_P^{\alpha}$. Then its stabilizer $H_3^{\alpha}=g_P^{\alpha} P(\mb O) (g_P^{\alpha})^{-1} \cap P(\mb O)$. 
The projection $P(\mb A) \twoheadrightarrow M(\mb A)$ induces a morphism
\begin{equation}   \label{equation-gPg-1-to-gMg-1}
g_P^{\alpha} P(\mb O) (g_P^{\alpha})^{-1} \cap P(\mb O)  \rightarrow   g_M^{\alpha} M(\mb O) (g_M^{\alpha})^{-1} \cap M(\mb O)
\end{equation}
Let $R_3^{\alpha}$ be the kernel of (\ref{equation-gPg-1-to-gMg-1}). We have $R_3^{\alpha} = U(\mb O) \cap H_3^{\alpha}$. 
The morphism $H_3^{\alpha} / R_3^{\alpha} \hookrightarrow H_2 / R_3^{\alpha}$ induces a morphism
\begin{equation}   
\Cht_{M, \infty, I, W} / (H_3^{\alpha} / R_3^{\alpha} ) \rightarrow \Cht_{M, \infty, I, W} / (H_2/ R_3^{\alpha}  )
\end{equation}
By counit morphism, it induces
\begin{equation}
H_{M, H_3^{\alpha}, R_3^{\alpha}, I, W}^{' \, j} \xrightarrow{ \on{Co} } H_{M, H_2, R_3^{\alpha}, I, W}^{' \, j}
\end{equation}

Moreover, $H_3^{\alpha} / R_3^{\alpha}  \subset g_M^{\alpha} M(\mb O) (g_M^{\alpha})^{-1} \cap M(\mb O) $ induces 
$$\Cht_{M, \infty, I, W} / ( H_3^{\alpha} / R_3^{\alpha} ) \rightarrow \Cht_{M, \infty, I, W} / \big( g_M^{\alpha} M(\mb O) (g_M^{\alpha})^{-1} \cap M(\mb O) \big) $$
We have morphisms 
\begin{equation}     \label{equation-gMg-1-adj-co} 
H_{M, g_M^{\alpha} M(\mb O) (g_M^{\alpha})^{-1} \cap M(\mb O) , I, W}^j   \xrightarrow{ \on{adj} }
H_{M, H_3^{\alpha}, R_3^{\alpha}, I, W}^{' \, j} \xrightarrow{ \on{Co} } H_{M, g_M^{\alpha} M(\mb O) (g_M^{\alpha})^{-1} \cap M(\mb O) , I, W}^j.
\end{equation}
The composition is the multiplication by $\sharp\big( (  g_M^{\alpha} M(\mb O) (g_M^{\alpha})^{-1} \cap M(\mb O) )  / ( H_3^{\alpha} / R_3^{\alpha}  )  \big)$.

Since $R_3^{\alpha} \subset U(\mb O)$, as in \ref{subsection-CT-level-K}, the projection $P(\mb O) / R_3^{\alpha} \twoheadrightarrow P(\mb O) / U(\mb O)$ induces a morphism 
$$\Cht_{M, \infty, I, W} / (P(\mb O) / R_3^{\alpha}) \rightarrow \Cht_{M, \infty, I, W} / (P(\mb O) / U(\mb O)) $$
which is a gerbe for the finite $q$-group $U(\mb O) / R_3^{\alpha}$. We deduce
\begin{equation}   \label{equation-R-3-adj-co}
H_{M, I, W}^{j} = H_{M, H_2, U(\mb O), I, W}^{' \, j} \underset{\sim}{\xrightarrow{  \on{adj}  } } H_{M, H_2, R_3^{\alpha}, I, W}^{' \, j}   \underset{\sim}{\xrightarrow{  \on{Co}  } }  H_{M, H_2, U(\mb O), I, W}^{' \, j} = H_{M, I, W}^{j}    
\end{equation}
where the composition is the multiplication by $(\sharp U(\mb O) / R_3^{\alpha} )^{-1}$.

We have a commutative diagram
\begin{equation}
\xymatrix{ 
H_{M, H_3^{\alpha} , R_3^{\alpha}, I, W}^{' \, j}  \ar[r]^{\on{Co}}    \ar[d]^{\on{Co}}
& H_{M, H_2, R_3^{\alpha}, I, W}^{' \, j}  \ar[d]^{\on{Co}}_{\simeq}  \\
H_{M, g_M^{\alpha} M(\mb O) (g_M^{\alpha})^{-1} \cap M(\mb O) , I, W}^j   \ar[r]^{\quad \quad \quad \on{Co}}  
&  H_{M, I, W}^{j} 
}
\end{equation}

\sssec{}
Note that $H_3^{\alpha} = g_P^{\alpha} H_1^{\alpha} (g_P^{\alpha})^{-1}$ and $R_3^{\alpha} = g_P^{\alpha} R_1^{\alpha} (g_P^{\alpha})^{-1}$.
The action of $g_P^{\alpha}$ (where $g_P^{\alpha}$ acts on $\Cht_{M, \infty, I, W}$ by $g_M^{\alpha}$) induces an isomorphism
$$ \Cht_{M, \infty, I, W} / (H_3^{\alpha} / R_3^{\alpha}) \underset{\sim}{\xrightarrow{  g_P^{\alpha}  }}  \Cht_{M, \infty, I, W} / (H_1^{\alpha} / R_1^{\alpha} )$$
We deduce an isomorphism of cohomology groups by the adjunction morphism.
\begin{equation}  
H_{M, H_1^{\alpha} , R_1^{\alpha}, I, W}^{' \, j}   \underset{\sim}{\xrightarrow{  \on{adj}(g_P^{\alpha})  }}  
H_{M, H_3^{\alpha} , R_3^{\alpha}, I, W}^{' \, j} 
\end{equation}

\quad

\noindent {\bf Proof of Lemma \ref{lem-h-act-on-H-M-equal-h-M-act-on-H-M}:} 

For each $\alpha$, we have the following diagram, whose commutativity will be proven next:
\begin{equation}   \label{equation-diam-alpha-action-h-vs-h-M}
{  \resizebox{15cm}{!}{ 
\!\!\!\!\!\!\!\!\!\!\!\!\!
\xymatrixrowsep{2pc}
\xymatrixcolsep{3pc}
\xymatrix{
H_{M, H_2, R_1^{\alpha}, I, W}^{' \, j}  \ar[r]^{\on{adj}}
& H_{M, H_1^{\alpha}, R_1^{\alpha}, I, W}^{' \, j} \ar[r]^{\on{adj}(g_P^{\alpha})}   
& H_{M, H_3^{\alpha} , R_3^{\alpha}, I, W}^{' \, j}  \ar[r]^{\on{Co}}     
& H_{M, H_2, R_3^{\alpha}, I, W}^{' \, j}  \ar[d]^{\on{Co}}_{\simeq}  \\
H_{M, I, W}^{j} \ar[r]^{\on{adj} \quad \quad \quad}   \ar[u]^{\on{adj} }_{\simeq}
& H_{M, M(\mb O) \cap (g_M^{\alpha})^{-1} M(\mb O) g_M^{\alpha}, I, W}^{ j}  \ar[r]^{ c^{\alpha} \on{adj}(g_M^{\alpha})}    
& H_{M, g_M^{\alpha} M(\mb O) (g_M^{\alpha})^{-1} \cap M(\mb O) , I, W}^j   \ar[r]^{\quad \quad \quad \on{Co}}  
&  H_{M, I, W}^{j} 
}
 } }
\end{equation}

where $c^{\alpha}$ is defined in \ref{subsection-h-M-scts-on H-M-I-W}.

The morphism
$$P(\mc O_v) g_P^{\alpha} P(\mc O_v) / P(\mc O_v) \twoheadrightarrow M(\mc O_v) g_M^{\alpha} M(\mc O_v)  / M(\mc O_v)$$
coincides with the morphism 
$$ P(\mc O_v) / g_P^{\alpha} P(\mc O_v) (g_P^{\alpha})^{-1} \cap P(\mc O_v) \twoheadrightarrow M(\mc O_v) / g_M^{\alpha} M(\mc O_v) (g_M^{\alpha})^{-1} \cap M(\mc O_v) $$
i.e.
$$ P(\mb O) / g_P^{\alpha} P(\mb O) (g_P^{\alpha})^{-1} \cap P(\mb O) \twoheadrightarrow M(\mb O) / g_M^{\alpha} M(\mb O) (g_M^{\alpha})^{-1} \cap M(\mb O) $$
We deduce that 
$$c^{\alpha} = \sharp ( U(\mb O) / R_1^{\alpha} ) \cdot \sharp \big( ( M(\mb O) \cap (g_M^{\alpha})^{-1} M(\mb O) g_M^{\alpha} ) / (H_1^{\alpha} / R_1^{\alpha})  \big)$$
Similarly, $c^{\alpha} = \sharp ( U(\mb O) / R_3^{\alpha} ) \cdot \sharp \big(  ( g_M^{\alpha} M(\mb O) (g_M^{\alpha})^{-1} \cap  M(\mb O) )  / (H_3^{\alpha} / R_3^{\alpha})  \big)$.
Moreover, since $H_3^{\alpha} = g_P^{\alpha} H_1^{\alpha} (g_P^{\alpha})^{-1}$ and $R_3^{\alpha} = g_P^{\alpha} R_1^{\alpha} (g_P^{\alpha})^{-1}$, we have 
$$\sharp \big( ( M(\mb O) \cap (g_M^{\alpha})^{-1} M(\mb O) g_M^{\alpha} ) / (H_1^{\alpha} / R_1^{\alpha})  \big) = \sharp \big(  ( g_M^{\alpha} M(\mb O) (g_M^{\alpha})^{-1} \cap  M(\mb O) )  / (H_3^{\alpha} / R_3^{\alpha})  \big)$$
$$\sharp ( U(\mb O) / R_1^{\alpha} ) = \sharp ( U(\mb O) / R_3^{\alpha} ) .$$ 
Taking into account (\ref{equation-g-1Mg-adj-co}), (\ref{equation-R-adj-co}), (\ref{equation-gMg-1-adj-co}) and (\ref{equation-R-3-adj-co}),
we deduce that diagram (\ref{equation-diam-alpha-action-h-vs-h-M}) is commutative.

Taking direct sum on $\alpha \in \Omega$, we deduce that the following diagram is commutative 
\begin{equation}
{  \resizebox{14cm}{!}{ 
\!\!\!\!\!\!\!\!\!\!\!
\xymatrixrowsep{2pc}
\xymatrixcolsep{3pc}
\xymatrix{
& \oplus_{\alpha} H_{M, H_1^{\alpha}, R_1^{\alpha}, I, W}^{' \, j} \ar[r]^{\on{adj}(g_P^{\alpha})}   
& \oplus_{\alpha} H_{M, H_3^{\alpha} , R_3^{\alpha}, I, W}^{' \, j}  \ar[rd]^{\on{Co}}  
&  \\
H_{M, I, W}^{j} \ar[r]^{\on{adj} \quad \quad \quad}   \ar[ur]^{\on{adj} }
& \oplus_{\alpha} H_{M, M(\mb O) \cap (g_M^{\alpha})^{-1} M(\mb O) g_M^{\alpha}, I, W}^{ j}  \ar[r]^{ c^{\alpha} \on{adj}(g_M^{\alpha})}  
& \oplus_{\alpha} H_{M, g_M^{\alpha} M(\mb O) (g_M^{\alpha})^{-1} \cap M(\mb O) , I, W}^j   \ar[r]^{\quad \quad \quad \on{Co}}  
&  H_{M, I, W}^{j} 
}
} }
\end{equation}

Taking into account (\ref{equation-def-H-M-K-I-W-prime}), the composition of the upper line is the action of $T(h)$ given in (\ref{equation-T-h-act-on-H-M-I-W}), the composition of the lower line is the action of $T(h^M)$ given in \ref{subsection-h-M-scts-on H-M-I-W}. We deduce the lemma.
\cqfd

\quad

\begin{lem}    \label{lem-TC-commute-with-Hecke-without-niveau}
(Hecke action at any unramified place commutes with constant term morphism)
For any place $v$ of $X \sm N$ and any $h \in C_{c}( G(\mc O_v) \backslash G(F_v)/ G(\mc O_v), \Oe)$, the following diagram of cohomology groups is commutative:
\begin{equation}   \label{diagram-TC-commute-avec-Hecke}
\xymatrixrowsep{2pc}
\xymatrixcolsep{3pc}
\xymatrix{
  H_{G, N, I, W}^{j, \, \Oe}  \ar[r]^{T(h)}  \ar[d]_{ C_{G, N}^{P, \, j, \, \Oe}} 
&    H_{G, N, I, W}^{j, \, \Oe}   \ar[d]^{ C_{G, N}^{P, \, j, \, \Oe}}  \\
 H_{M, N, I, W}^{' \, j, \, \Oe}   \ar[r]^{T(h^M)}
&   H_{M, N, I, W}^{' \, j, \, \Oe} 
}
\end{equation}
where the vertical morphisms are the constant term morphism defined in Definition \ref{def-CT-cohomology}.
\end{lem}
\dem
Follows from Proposition \ref{prop-TC-commute-with-Hecke-with-niveau} and Lemma \ref{lem-h-act-on-H-M-equal-h-M-act-on-H-M}.
\cqfd

\begin{rem}
The direct argument in the proof of \cite{these} Lemme 8.1.1 also works for $\Oe$-coefficients. This gives another proof of Lemma \ref{lem-TC-commute-with-Hecke-without-niveau}.
\end{rem}

\quad

\section{Finiteness of cohomology with integral coefficients as module over Hecke algebras}

Let $\ov{x}$ be a geometric point of $(X \sm N)^I$ as in \ref{subsection-ov-x}. Sections 1 and 2 of \cite{coho-filt} still hold for $\Oe$-coefficients. 

\sssec{}
For any place $u$ of $X \sm N$, we denote by $\ms H_{G, u}^{\Oe}:=C_{c}(G(\mc O_u)\backslash G(F_u)/G(\mc O_u), \Oe)$ the Hecke algebra at $u$.

\begin{prop}   \label{prop-H-G-equal-Hecke-H-G-leq-mu-Zl}
There exists $\mu_1 \in \wh \Lambda_{G^{\on{ad}}}^{+, \Q}$ large enough, such that
\begin{equation}
\restr{ \mc H_{G, N, I, W}^{j, \, \Oe} }{\ov{x}} = \ms H_{G, u}^{\Oe} \cdot \on{Im} \big( \restr{ \mc H_{G, N, I, W}^{j, \, \leq \mu_1, \, \Oe} }{\ov{x}}  \rightarrow \restr{ \mc H_{G, N, I, W}^{j, \, \Oe} }{\ov{x}}   \big)
\end{equation}
\end{prop}
\dem
Similar to \cite{coho-filt} Definition 2.2.2, we define a filtration $F_G^{\bullet}$ of $\restr{ \mc H_{G, N, I, W}^{j,  \, \Oe} }{\ov{x}} $:
$$\mu \in \frac{1}{r}\wh R_{G^{\mr{ad}}}^+, \quad F_G^{\mu}:=\on{Im}(  \restr{ \mc H_{G, N, I, W}^{j, \, \leq\mu, \, \Oe} }{\ov{x}}  \rightarrow \restr{ \mc H_{G, N, I, W}^{j, \, \Oe} }{\ov{x}}      ) .$$
For every maximal parabolic subgroup $P$ of $G$ with Levi quotient $M$,
similar to $loc.cit.$ Definition 2.4.2, 
we define a filtration $F_M^{\bullet}$ of $\restr{ \mc H_{M, N, I, W}^{' \, j,  \, \Oe} }{\ov{x}} $:
$$\mu \in \frac{1}{r}\wh R_{G^{\mr{ad}}}^+, \quad F_M^{\mu}:=\on{Im}(  \restr{ \mc H_{M, N, I, W}^{' \, j, \, \leq\mu, \, \Oe} }{\ov{x}}  \rightarrow \restr{ \mc H_{M, N, I, W}^{' \, j, \, \Oe} }{\ov{x}}      ) .$$
For any $\mu_1, \mu_2 \in \frac{1}{r}\wh R_{G^{\mr{ad}}}^+$ and $\mu_1 \leq \mu_2$, $F_{G}^{\mu_1}$ (resp. $F_{M}^{\mu_1}$ ) is a sub-$\Oe$-module of $F_G^{\mu_2}$ (resp. $F_{M}^{\mu_2}$).

As in $loc.cit.$, let $\omega$ be the quasi-fundamental coweight of $G$ in the center of $M$ and $h_{\omega}^G \in \ms H_{G, u}^{\Oe}$ be the characteristic function of $G(\mc O_u) \varpi^{\omega} G(\mc O_u)$, where $\varpi$ is a uniformizer of $\mc O_u$.
Let $\check{\alpha}$ be the simple coroot of $G$ which is not a simple coroot of $M$. Let $l = \deg(u)$. 
Taking into account Lemma \ref{lem-TC-commute-with-Hecke-without-niveau}, as in $loc.cit.$ 2.5, we have a commutative diagram
\begin{equation}   \label{equation-filt-quotient-commut-diag}
\xymatrix{
F_{G}^{\mu- \omega l} / F_{G}^{\mu- \omega l - \frac{1}{r}\check{\alpha} }  \ar[r]^{\quad h^G_{\omega}}   \ar[d]_{\ov{C_G^P}}
& F_{G}^{\mu} / F_{G}^{\mu - \frac{1}{r} \check{\alpha} }   \ar[d]_{\ov{C_G^P}} \\
F_{M}^{\mu- \omega l} / F_{M}^{\mu- \omega - \frac{1}{r}\check{\alpha} l}   \ar[r]^{\quad h^G_{\omega}}_{ \quad \simeq }
& F_{M}^{\mu } / F_{M}^{\mu - \frac{1}{r} \check{\alpha}}  
}
\end{equation}
where the horizontal morphisms are induced by the action of the Hecke operator $h^G_{\omega}$, and the vertical morphisms are induced by the constant term morphisms. The lower line is an isomorphism for the same reason as in $loc.cit.$ Lemma 2.4.3.

By Proposition \ref{prop-TC-isom-mu-grand} and Proposition \ref{prop-finiteness-a-b-c} in this paper, $loc.cit.$ Proposition 2.3.14 still hold for cohomology groups with $\Oe$-coefficients. Thus for $\mu$ large enough as in $loc.cit.$ Lemma 2.2.11, the vertical morphisms in (\ref{equation-filt-quotient-commut-diag}) are isomorphisms. So the upper line of (\ref{equation-filt-quotient-commut-diag}) is an isomorphism.

Then we use the same argument as in the proof of $loc.cit.$ Proposition 2.2.4 for the filtrations $F_G^{\bullet}$ defined above. This proves Proposition \ref{prop-H-G-equal-Hecke-H-G-leq-mu-Zl}.
\cqfd

\sssec{}
Proposition \ref{prop-H-G-equal-Hecke-H-G-leq-mu-Zl} and the fact that $\restr{ \mc H_{G, N, I, W}^{j, \, \leq \mu_1, \, \Oe} }{\ov{x}}$ is a $\Oe$-module of finite type imply the following

\begin{thm}     \label{coho-cht-Oe-is-Hecke-mod-type-fini-en-u}
The cohomology group $\restr{ \mc H_{G, N, I, W}^{j, \, \Oe} }{\ov{x}}$ is of finite type as a $\ms H_{G, u}^{\Oe}$-module.
\end{thm}

\begin{cor}   \label{cor-torsion-in-coho-bounded}
The order of torsion in $\restr{ \mc H_{G, N, I, W}^{j, \, \Oe} }{\ov{x}} $ as a $\Oe$-module is bounded.
\end{cor}
\dem
Since $\ms H_{G, u}^{\Oe}$ is Noetherian, $\restr{ \mc H_{G, N, I, W}^{j, \, \Oe} }{\ov{x}}$ is also Noetherian. Let $K_n$ be the sub-$\Oe$-module of $\restr{ \mc H_{G, N, I, W}^{j, \, \Oe} }{\ov{x}} $ of elements $a$ such that $\ell^n a =0$. Then $K_1, K_2, \cdots , K_n, \cdots $ form an ascending chain. Thus it is stationary. 
\cqfd

\quad

\section{Cuspidal cohomology and Hecke-finite cohomology}

In this section, let $\ov{x}$ be a geometric point of $(X \sm N)^I$ as in \ref{subsection-ov-x}. 

In \cite{cusp-coho}, we proved that
$$\big( \restr{ \mc H_{G, N, I, W}^{j, \, \Ql}}{ \ov{x} } \big)^{\on{cusp}} = \big( \restr{ \mc H_{G, N, I, W}^{j,  \, \Ql} }{ \ov{x} } \big)^{\on{Hf-rat}} \supset \big( \restr{ \mc H_{G, N, I, W}^{j,  \, \Ql} }{ \ov{x} } \big)^{\on{Hf}} .$$
In this section, firstly we will prove in Proposition \ref{prop-cusp-egal-Hf-Oe} that
$$\big( \restr{ \mc H_{G, N, I, W}^{j, \, \Oe}}{ \ov{x} } \big)^{\on{cusp}} = \big( \restr{ \mc H_{G, N, I, W}^{j,  \, \Oe} }{ \ov{x} } \big)^{\on{Hf}} .$$
Then we will prove that the image of $\big( \restr{ \mc H_{G, N, I, W}^{j, \, \Oe}}{ \ov{x} } \big)^{\on{cusp}}$ in $\big( \restr{ \mc H_{G, N, I, W}^{j, \, \Ql}}{ \ov{x} } \big)^{\on{cusp}} $ is an $\Oe$-lattice stable under the action of Hecke algebras. Finally we will prove in Proposition \ref{prop-H-cusp-equal-H-Hf} that
$$\big( \restr{ \mc H_{G, N, I, W}^{j, \, \Ql}}{ \ov{x} } \big)^{\on{cusp}} = \big( \restr{ \mc H_{G, N, I, W}^{j,  \, \Ql} }{ \ov{x} } \big)^{\on{Hf}} .$$

\subsection{Cuspidal cohomology with $\Oe$-coefficients and Hecke-finite cohomology with $\Oe$-coefficients}

\begin{defi}  \label{def-H-Hf-rat}
Similarly to Définition 8.19 of \cite{vincent}, we define $\big( \restr{ \mc H_{G, N, I, W}^{j, \, \Oe} }{ \ov{x}  } \big)^{\on{Hf}}:=$
$$\{  c \in \restr{ \mc H_{G, N, I, W}^{j, \, \Oe} }{ \ov{x} } , \; \text{the } \Oe \text{-module } C_{c}(K_{N}\backslash G(\mb A)/K_{N}, \Oe) \cdot c \text{ is of finite type} \}.$$
\end{defi}

\begin{prop}   \label{prop-cusp-egal-Hf-Oe}
The two sub-$\Oe$-modules $\Hcusp$ and $\big( \restr{ \mc H_{G, N, I, W}^{j, \, \Oe} }{ \ov{x}  } \big)^{\on{Hf}}$ of $\restr{ \mc H_{G, N, I, W}^{j, \, \Oe} }{ \ov{x} }$ are equal.
\end{prop} 

\dem
It follows from Lemma \ref{lem-cusp-inclus-dans-Hf-rat} and Lemma \ref{lem-Hf-rat-inclus-dans-cusp} below.
\cqfd

\begin{lem}   \label{lem-cusp-inclus-dans-Hf-rat}
We have an inclusion
\begin{equation}
\Hcusp \subset \big( \restr{ \mc H_{G, N, I, W}^{j, \, \Oe} }{ \ov{x}  } \big)^{\on{Hf}}
\end{equation}
\end{lem}

\dem
By Theorem \ref{thm-cusp-dim-fini-second}, the $\Oe$-module $\Hcusp$ is of finite type. 
By Proposition \ref{prop-TC-commute-with-Hecke-with-niveau}, it is stable under the action of the Hecke algebra $C_{c}(K_{N}\backslash G(\mb A)/K_{N}, \Oe)$. 
We conclude by Definition \ref{def-H-Hf-rat}.
\cqfd

\begin{lem}   \label{lem-Hf-rat-inclus-dans-cusp}
We have an inclusion 
\begin{equation}
\Hcusp \supset \big( \restr{ \mc H_{G, N, I, W}^{j, \, \Oe} }{ \ov{x}  } \big)^{\on{Hf}}
\end{equation}
\end{lem}

\dem
The proof is similar to the proof of \cite{cusp-coho} Lemma 6.3.3. 
Let $a \in \big( \restr{ \mc H_{G, N, I, W}^{j, \, \Oe} }{ \ov{x}  } \big)^{\on{Hf}}$. 
We argue by contradiction. Suppose that $a \notin \Hcusp$, then there exists a maximal parabolic subgroup $P$ with Levi quotient $M$ such that $b: =C_G^{P, \, j, \, \Oe}(a)$ is non zero. Let $v$ be a place in $X \sm N$.

(1) On the one hand, $C_c(G(\mc O_v) \backslash G(F_v) / G(\mc O_v), \Oe) \cdot a$ is a $\Oe$-module of finite type. By Lemma \ref{lem-TC-commute-with-Hecke-without-niveau} and Lemma \ref{lem-Hecke-M-fini-sur-Hecke-G} below, $C_c(M(\mc O_v) \backslash M(F_v) / M(\mc O_v), \Oe) \cdot b$ is also a $\Oe$-module of finite type.

(2) On the other hand, let $g \in Z_M(F)$ such that $g \notin Z_M(\mc O_v)Z_G(F_v)$, the action of the Hecke operator $T(g)$ associated to $g$ induces an isomorphism of stacks $\Cht_{M, N, I, W}^{' \, \nu} / \Xi \isom \Cht_{M, N, I, W}^{' \, \nu+\xi(g)} / \Xi $, where $\xi(g) \in \wh \Lambda_{Z_M / Z_G}^{\Q} \simeq \Q$ and $\xi(g)>0$ (if not, take $g^{-1}$ instead of $g$). It induces an isomorphism of cohomology groups $\restr{ \mc H_{M, N, I, W}^{' \, j, \, \nu, \, \Oe} }{ \ov{x}  } \isom \restr{ \mc H_{M, N, I, W}^{' \, j, \, \nu+\xi(g), \, \Oe} }{ \ov{x}  }$.

There exists $\mu \in \wh \Lambda_{G^{\mr{ad}}}^{+, \Q}$ such that $a \in \on{Im}( \restr{ \mc H_{G, N, I, W}^{j, \, \leq \mu \, \Oe} }{ \ov{x}  }  \rightarrow \restr{ \mc H_{G, N, I, W}^{j, \, \Oe} }{ \ov{x}  }   )$. By \ref{subsection-C-G-P-leq-mu-Oe}, $b \in \on{Im}( \restr{ \mc H_{M, N, I, W}^{' \, j, \, \leq \mu, \, \Oe} }{ \ov{x}  } \rightarrow \restr{ \mc H_{M, N, I, W}^{' \, j, \, \Oe} }{ \ov{x}  } ).$
By Remark \ref{rem-CT-cohomology-nu-Oe} and Remark \ref{rem-Ker-prod-nu-TC-equal-Ker-TC}, we have an injective morphism
$$\restr{ \mc H_{M, N, I, W}^{' \, j, \, \Oe} }{ \ov{x}  } \hookrightarrow \prod_{\nu} \restr{ \mc H_{M, N, I, W}^{' \, j, \, \nu, \, \Oe} }{ \ov{x}  }$$
We consider the image of $b$ in $\prod_{\nu} \restr{ \mc H_{M, N, I, W}^{' \, j, \, \nu, \, \Oe} }{ \ov{x}  }$. It is supported on the components $\restr{ \mc H_{M, N, I, W}^{' \, j, \, \nu, \, \Oe} }{ \ov{x}  } $ indexed by $\nu$ in the translated cone $\wh \Lambda_{Z_M / Z_G}^{\mu} \subset \wh \Lambda_{Z_M / Z_G}^{\Q}$. Since $P$ is maximal, $\wh \Lambda_{Z_M / Z_G}^{\Q} \simeq \Q$. The translated cone $\wh \Lambda_{Z_M / Z_G}^{\mu}$ is of the form $\{\nu \in \Q, \nu \leq \mu\}$.

Since $b\neq 0$, there exists $m \in \Z_{>0}$ such that $T(g)^m (b)$ is supported on the cone $\wh \Lambda_{Z_M / Z_G}^{\mu + m\xi(g)} \supset \wh \Lambda_{Z_M / Z_G}^{\mu}$, but not supported on $\wh \Lambda_{Z_M / Z_G}^{\mu}$. Therefore $T(g)^{2m} (b)$ is supported on the cone $\wh \Lambda_{Z_M / Z_G}^{\mu + 2m\xi(g)} $, but not supported on $\wh \Lambda_{Z_M / Z_G}^{\mu + m\xi(g)}$, etc. We deduce that for any $k>0$, $T(g)^{(k+1)m}(b)$ is not included in the sub-$\Oe$-module of $\restr{ \mc H_{M, N, I, W}^{' \, j, \, \Oe} }{ \ov{x}  } $ generated by $b, T(g)^m(b), T(g)^{2m}(b), \cdots, T(g)^{km}(b)$. 

Thus the sub-$\Oe$-module of $C_c(M(\mc O_v) \backslash M(F_v) / M(\mc O_v), \Oe) \cdot b$  generated by $T(g)^{\Z}(b)$ is not of finite type. Since $\Oe$ is Noetherian, we deduce that $C_c(M(\mc O_v) \backslash M(F_v) / M(\mc O_v), \Oe) \cdot b$ is not of finite type.

(3) We deduce from (1) and (2) a contradiction. 
\cqfd

\begin{lem} \label{lem-Hecke-M-fini-sur-Hecke-G} (\cite{gross})
Under the Satake transformation (\ref{equation-TC-pour-op-de-Hecke}), the algebra $C_{c}(M(O_v)\backslash M(F_v) / M(O_v), \Oe)$ is finite over $C_{c}(G(O_v)\backslash G(F_v)/G(O_v), \Oe)$.
\cqfd
\end{lem}

\subsection{Compatibility of cohomology with $\Oe$-coefficients and with $E$-coefficients}

\sssec{}
Let $$ \on{Sat}_{G, I}^{\Oe}: \on{Rep}_{\Oe}(\wh G^I) \rightarrow \Perv_{G_{I, \infty}}(\Gr_{G, I}, \Oe) $$ be the functor defined in Theorem \ref{thm-Satake-functor-I}.
Let $$ \on{Sat}_{G, I}^{\Ql}: \on{Rep}_{\Ql}(\wh G^I_E) \rightarrow \Perv_{G_{I, \infty}}(\Gr_{G, I}, \Ql) $$ be the functor defined in \cite{cusp-coho} Corollary 2.1.7.

\sssec{}   \label{subsection-Zl-lattice-in-W}
The canonical functor $\on{Rep}_{\Oe}(\wh G^I) \otimes_{\Oe} \Ql \rightarrow \on{Rep}_{\Ql}(\wh G^I_E)$ is an equivalence of categories. The inverse functor sends any $W \in \on{Rep}_{\Ql}(\wh G^I_E)$ to $W^{\Oe} \otimes_{\Oe} \Ql $, where $W^{\Oe}$ is any $\Oe$-lattice in $W$ stable by $\wh G^I$. Such a lattice exists for the following reason.
Let $W^*$ be the dual of $W$.
Let $\xi_1, \cdots, \xi_n$ be a base of $W^*$. For any $w \in W$ and any $i \in \{1, \cdots, n\}$, let $f_{w, \xi_i}: \wh G^I_E \rightarrow \Ql$ be the function of matrix coefficient associated to $w$ and $\xi_i$. Let
$$W^{\Oe}:= \{ w \in W \, | \, \forall 1 \leq  i \leq n, f_{w, \xi_i} \text{ comes from a function } \wh G^I \rightarrow \Oe  \} .$$
Then $W^{\Oe}$ is a $\Oe$-lattice in $W$ stable under the action of $\wh G^I$.

\sssec{}
By definition, we have a canonical equivalence of categories $$\Perv_{G_{I, \infty}}(\Gr_{G, I}, \Oe) \otimes_{\Oe} \Ql  \isom \Perv_{G_{I, \infty}}(\Gr_{G, I}, \Ql) .$$

\begin{prop}    \label{prop-Satake-compatible-Zl-Ql}
The following diagram of categories is commutative:
\begin{equation}
\xymatrixrowsep{2pc}
\xymatrixcolsep{7pc}
\xymatrix{
\on{Rep}_{\Oe}(\wh G^I) \otimes_{\Oe} \Ql  \ar[r]^{\on{Sat}_{G, I}^{\Oe} \otimes_{\Oe} \Ql \quad }   \ar[d]^{\simeq}
& \Perv_{G_{I, \infty}}(\Gr_{G, I}, \Oe) \otimes_{\Oe} \Ql  \ar[d]^{\simeq}  \\
\on{Rep}_{\Ql}(\wh G^I_E)  \ar[r]^{ \on{Sat}_{G, I}^{\Ql} }
& \Perv_{G_{I, \infty}}(\Gr_{G, I}, \Ql)  
}
\end{equation}
\end{prop}
\dem
\cite{mv} Theorem 14.1.
\cqfd

\sssec{}  \label{subsection-S-W-Zl-S-W-Ql}
Let $W \in \on{Rep}_{\Oe}(\wh G^I)$. By Theorem \ref{prop-Satake-compatible-Zl-Ql}, we have a canonical isomorphism
\begin{equation}   \label{equation-S-G-Zl-S-G-Ql}
 \on{Sat}_{G, I}^{\Oe}(W) \otimes_{\Oe} \Ql \isom \on{Sat}_{G, I}^{\Ql}(W \otimes_{\Oe} \Ql) .
\end{equation}

\sssec{}   \label{subsection-H-G-Zl-H-G-Ql}
Let $W \in \on{Rep}_{\Oe}(\wh G^I)$. Let $\mc F_{G, N, I, W \otimes_{\Oe} \Ql}^{\Ql}$ be the perverse sheaf on $\Cht_{G, N, I}$ associated to $W \otimes_{\Oe} \Ql$ defined in \cite{cusp-coho} Definition 2.4.5.
We deduce from (\ref{equation-S-G-Zl-S-G-Ql}) that
\begin{equation}   \label{equation-F-G-Zl-F-G-Ql}
\mc F_{G, N, I, W}^{\Oe} \otimes_{\Oe} \Ql \isom \mc F_{G, N, I, W \otimes_{\Oe} \Ql}^{\Ql} .
\end{equation}

In $loc.cit.$ Definition 2.5.3, we defined the sheaf $\mc H_{G, N, I, W \otimes_{\Oe}\Ql }^{j, \, \Ql}$ over $(X \sm N)^I$. 
(\ref{equation-F-G-Zl-F-G-Ql}) induces a canonical isomorphism
\begin{equation}   \label{equation-H-G-Zl-H-G-Ql-leq-mu}
 \mc H_{G, N, I, W}^{j, \, \leq \mu, \, \Oe}  \otimes_{\Oe} \Ql \isom  \mc H_{G, N, I, W \otimes_{\Oe}\Ql }^{j, \, \leq \mu, \, \Ql}  .
\end{equation}
Take limit on $\mu$. Since $\otimes_{\Oe} E$ commutes with inductive limit, we deduce
\begin{equation}   \label{equation-H-G-Zl-H-G-Ql}
 \mc H_{G, N, I, W}^{j, \, \Oe}  \otimes_{\Oe} \Ql \isom  \mc H_{G, N, I, W \otimes_{\Oe}\Ql }^{j, \, \Ql}  .
\end{equation}

\quad

\sssec{}    \label{subsection-H-M-Zl-H-M-Ql}
Applying Proposition \ref{prop-Satake-compatible-Zl-Ql} to $M$ (replacing $G$ by $M$), we have a canonical isomorphism
\begin{equation}  \label{equation-S-M-Zl-S-M-Ql}
\on{Sat}_{M, I}^{\Oe}(W) \otimes_{\Oe} \Ql \isom \on{Sat}_{M, I}^{\Ql}(W \otimes_{\Oe} \Ql) .
\end{equation}
Let $\mc F_{M, N, I, W \otimes_{\Oe} \Ql}^{' \, \Ql}$ be the perverse sheaf on $\Cht_{M, N, I}'$ associated to $W \otimes_{\Oe} \Ql$ defined in \cite{cusp-coho} Definition 3.4.7. We deduce from (\ref{equation-S-M-Zl-S-M-Ql}) that
\begin{equation}   \label{equation-F-M-Zl-F-M-Ql}
\mc F_{M, N, I, W}^{' \, \Oe} \otimes_{\Oe} \Ql \isom \mc F_{M, N, I, W \otimes_{\Oe} \Ql}^{' \, \Ql} .
\end{equation}
It induces a canonical isomorphism
\begin{equation}   
 \mc H_{M, N, I, W}^{' \, j, \, \leq \mu, \, \nu, \, \Oe}  \otimes_{\Oe} \Ql \isom  \mc H_{M, N, I, W \otimes_{\Oe}\Ql }^{' \, j, \, \leq \mu, \, \nu, \, \Ql}  .
\end{equation}
Take limit on $\mu$. Since $\otimes_{\Oe} E$ commutes with inductive limit, we deduce
\begin{equation}   
 \mc H_{M, N, I, W}^{' \, j, \, \nu, \, \Oe}  \otimes_{\Oe} \Ql \isom  \mc H_{M, N, I, W \otimes_{\Oe}\Ql }^{' \, j, \, \nu, \, \Ql}  .
\end{equation}

\begin{lem}   \label{lem-order-in-H-M-leq-mu-nu-bounded}
The order of torsion in $\restr{ \mc H_{M, N, I, W}^{' \, j, \, \nu, \, \Oe} }{\ov x}$ is bounded when $\nu$ varies.  
\end{lem}
\dem
By Corollary \ref{cor-torsion-in-coho-bounded} applied to $M$, for every $\nu \in \wh \Lambda_{Z_M / Z_G}^{\Q}$, the order of torsion in $\restr{ \mc H_{M, N, I, W}^{' \, j, \, \nu, \, \Oe} }{\ov x}$ is bounded. 

For any $\nu_1, \nu_2 \in \wh \Lambda_{Z_M / Z_G}$, there exists $g \in Z_M(F)$, such that the action of the Hecke operator associated to $g$ induces an isomorphism of stacks $\Cht_{M, N, I, W}^{' \, \nu_1} / \Xi \isom \Cht_{M, N, I, W}^{' \, \nu_2} / \Xi $. It induces an isomorphism of cohomology groups
\begin{equation}   \label{equation-H-M-nu-1-H-M-nu-2}
\restr{ \mc H_{M, N, I, W}^{' \, j, \, \nu_1, \, \Oe} }{\ov x} \isom \restr{ \mc H_{M, N, I, W}^{' \, j, \, \nu_2, \, \Oe} }{\ov x} .
\end{equation}

By \cite{cusp-coho} 1.5.7, for any $\nu \in \wh \Lambda_{Z_M / Z_G}^{\Q}$, the stack $\Cht_{M, N, I, W}^{' \, \nu}$ is non-empty if and only if $\nu$ is included in the image of $\pi_0(\Bun_M)$ in $\wh \Lambda_{Z_M / Z_G}^{\Q}$, which is a $\Z$-lattice in $\wh \Lambda_{Z_M / Z_G}^{\Q}$. Thus among all $\restr{ \mc H_{M, N, I, W}^{' \, j, \, \nu, \, \Oe} }{\ov x}$, $\nu \in \wh \Lambda_{Z_M / Z_G}^{\Q}$, there are only a finite number of isomorphism classes in the sense of (\ref{equation-H-M-nu-1-H-M-nu-2}). We deduce the lemma.
\cqfd

\sssec{}
By Lemma \ref{lem-order-in-H-M-leq-mu-nu-bounded}, the morphism
\begin{equation}   \label{equation-prod-H-M-nu-Zl-to-H-M-nu-Ql}
\big( \prod_{\nu \in \wh \Lambda_{Z_M / Z_G}^{\Q}} \restr{ \mc H_{M, N, I, W}^{' \, j, \, \nu, \, \Oe} }{\ov x} \big)  \otimes_{\Oe} \Ql  \rightarrow  \prod_{\nu \in \wh \Lambda_{Z_M / Z_G}^{\Q} } \restr{ \mc H_{M, N, I, W \otimes_{\Oe}\Ql }^{j, \, \nu, \, \Ql} }{ \ov x} 
\end{equation}
is injective.
Consider the morphism
\begin{equation}   \label{equation-H-M-Zl-H-M-Ql}
 \mc H_{M, N, I, W}^{' \, j, \, \Oe}  \otimes_{\Oe} \Ql \rightarrow \mc H_{M, N, I, W \otimes_{\Oe}\Ql }^{' \, j, \, \Ql}  .
\end{equation}
We have a commutative diagram
$$
\xymatrixrowsep{2pc}
\xymatrixcolsep{6pc}
\xymatrix{
\big( \prod_{\nu} \mc H_{M, N, I, W}^{j, \, \nu, \, \Oe} \big)  \otimes_{\Oe} \Ql  \ar@{^{(}->}[r]^{ (\ref{equation-prod-H-M-nu-Zl-to-H-M-nu-Ql})  }
&  \prod_{\nu} \mc H_{M, N, I, W \otimes_{\Oe}\Ql }^{j, \, \nu, \, \Ql} \\
 \mc H_{M, N, I, W}^{' \, j, \, \Oe}  \otimes_{\Oe} \Ql    \ar@{^{(}->}[u]   \ar[r]^{  (\ref{equation-H-M-Zl-H-M-Ql})  }
 & \mc H_{M, N, I, W \otimes_{\Oe}\Ql }^{' \, j, \, \Ql}    \ar@{^{(}->}[u]
}
$$
where the left vertical morphism is the tensor product with $E$ of morphism (\ref{equation-H-M-to-prod-H-M-nu-Oe}) defined in Remark \ref{rem-CT-cohomology-nu-Oe}, whose injectivity is proved in Remark \ref{rem-Ker-prod-nu-TC-equal-Ker-TC}. The right vertical morphism is defined in \cite{cusp-coho} Remark 3.5.11.

We deduce that morphism (\ref{equation-H-M-Zl-H-M-Ql}) is injective.

\sssec{}
We have a commutative diagram:
\begin{equation}   \label{diagram-CT-Ql-CT-Zl}
\xymatrixrowsep{2pc}
\xymatrixcolsep{7pc}
\xymatrix{
\restr{ \mc H_{G, N, I, W \otimes_{\Oe}\Ql }^{j, \, \Ql} }{\ov{x}}    \ar[r]^{  C_{G, N}^{P, \, j, \, \Ql}  }
& \restr{ \mc H_{M, N, I, W \otimes_{\Oe}\Ql }^{' \, j, \, \Ql} }{\ov{x}} \\
\restr{ \mc H _{G, N, I, W}^{j, \, \Oe}  }{ \ov{x}  } \otimes_{\Oe} \Ql    \ar[r]^{  C_{G, N}^{P, \, j, \, \Oe} \otimes_{\Oe} \Ql  }   \ar[u]_{\simeq}^{(\ref{equation-H-G-Zl-H-G-Ql})}
&  \restr{ \mc  {H'}_{\! \! \! M, N, I, W}^{\; j, \, \Oe} }{ \ov{x}  } \otimes_{\Oe} \Ql   \ar@{^{(}->}[u]_{(\ref{equation-H-M-Zl-H-M-Ql})} ,
}
\end{equation}
where $C_{G, N}^{P, \, j, \, \Oe} \otimes_{\Oe} \Ql$ is induced by $C_{G, N}^{P, \, j, \, \Oe}$ (defined in Definition \ref{def-CT-cohomology}), $C_{G, N}^{P, \, j, \, \Ql}$ is defined in \cite{cusp-coho} Definition 3.5.10.

\sssec{}
In \cite{cusp-coho} Definition 3.5.14, we defined 
\begin{equation}
\big( \restr{ \mc H_{G, N, I, W \otimes_{\Oe}\Ql }^{j,  \, \Ql} }{\ov{x}} \big)^{\on{cusp}} = \underset{P \varsubsetneq G} \bigcap \Ker C_{G, N}^{P, \, j, \, \Ql}.
\end{equation}

\begin{prop}  \label{prop-H-cusp-Zl-H-cusp-Ql}
Let $W \in \on{Rep}_{\Oe}(\wh G^I)$. Then we have a canonical isomorphism
$$\Hcusp \otimes_{\Oe} \Ql \isom \big( \restr{ \mc H_{G, N, I, W \otimes_{\Oe}\Ql }^{j,  \, \Ql} }{\ov{x}} \big)^{\on{cusp}}.$$
\end{prop}
\dem
For every proper parabolic subgroup $P$ of $G$ with Levi quotient $M$, we have
\begin{equation}
\begin{aligned}
& \Ker \big( \restr{ \mc H_{G, N, I, W}^{j, \, \Oe} }{\ov{x}} \rightarrow \restr{ \mc H_{M, N, I, W}^{' \, j, \, \Oe} }{\ov{x}} \big) \otimes_{\Oe} \Ql  \\
\isom & \Ker \big( \restr{ \mc H_{G, N, I, W}^{j, \, \Oe} }{\ov{x}} \otimes_{\Oe} \Ql  \rightarrow \restr{ \mc H_{M, N, I, W}^{' \, j, \, \Oe} }{\ov{x}}  \otimes_{\Oe} \Ql  \big) \\
\isom & \Ker \big( \restr{ \mc H_{G, N, I, W \otimes_{\Oe}\Ql}^{j, \, \Ql} }{\ov{x}} \rightarrow \restr{ \mc H_{M, N, I, W \otimes_{\Oe}\Ql}^{' \, j, \, \Ql} }{\ov{x}} \big)
\end{aligned}
\end{equation}
The first isomorphism follows from the fact that $\Ql$ is flat over $\Oe$ (see \cite{matsumura} Theorem 4.4 and Theorem 4.5), the second isomorphism follows from diagram (\ref{diagram-CT-Ql-CT-Zl})

We deduce the proposition.
\cqfd

\begin{rem} 
We have a commutative diagram
\begin{equation}
\xymatrix{
\big( \restr{ \mc H_{G, N, I, W \otimes_{\Oe} \Ql}^{j,  \, \Ql } }{\ov{x}} \big)^{\on{cusp}}  \ar@{^{(}->}[r]
&  \restr{ \mc H_{G, N, I, W \otimes_{\Oe} \Ql }^{j,  \, \Ql} }{\ov{x}}  \\
\Hcusp   \ar@{^{(}->}[r]   \ar[u]
& \restr{ \mc H_{G, N, I, W}^{j,  \, \Oe} }{\ov{x}}   \ar[u]
}
\end{equation}
We do not know if it is Cartesian. We have an injective morphism 
\begin{equation}   \label{equation-H-cusp-Oe-to-H-cusp-E-intersect-H-Oe}
\Hcusp \rightarrow \big( \restr{ \mc H_{G, N, I, W \otimes_{\Oe} \Ql}^{j,  \, \Ql } }{\ov{x}} \big)^{\on{cusp}}  \underset{ \restr{ \mc H_{G, N, I, W \otimes_{\Oe} \Ql }^{j,  \, \Ql} }{\ov{x}} }{\times} \restr{ \mc H_{G, N, I, W}^{j,  \, \Oe} }{\ov{x}} 
\end{equation}
We do not know if (\ref{equation-H-cusp-Oe-to-H-cusp-E-intersect-H-Oe}) is surjective, due to the possible torsion in $\restr{ \mc H_{M, N, I, W}^{' \, j,  \, \Oe} }{\ov{x}}$ (thus in $\prod_{\nu} \restr{ \mc H_{M, N, I, W}^{j, \, \nu, \, \Oe} }{ \ov x }$). 
\end{rem}

\subsection{Cuspidal cohomology with $E$-coefficients and Hecke-finite cohomology with $E$-coefficients}

\sssec{}
Let $W \in \on{Rep}_{\Ql}(\wh G^I_E)$.
In \cite{cusp-coho} Definition 6.3.1, we defined
$$\big( \restr{ \mc H_{G, N, I, W}^{j, \, \Ql}}{ \ov{x} } \big)^{\on{Hf-rat}} :=\{  c \in \restr{ \mc H_{G, N, I, W}^{j, \, \Ql} }{ \ov{x} } , \; \dim_{\Ql} C_{c}(K_{N}\backslash G(\mb A)/K_{N}, \Ql) \cdot c < + \infty  \}.$$
In $loc.cit.$ Proposition 6.0.1, we proved that
\begin{equation}   \label{equation-H-cusp-equal-H-Hf-rat-Ql}
\big( \restr{ \mc H_{G, N, I, W}^{j, \, \Ql}}{ \ov{x} } \big)^{\on{cusp}} = \big( \restr{ \mc H_{G, N, I, W}^{j,  \, \Ql} }{ \ov{x} } \big)^{\on{Hf-rat}}
\end{equation}
as sub-$\Ql$-vector spaces of $\restr{ \mc H_{G, N, I, W}^{j, \, \Ql} }{ \ov{x} }$.

\quad

\sssec{}
In \cite{vincent} Définition 8.19, V. Lafforgue defined $\big( \restr{ \mc H_{G, N, I, W}^{j,  \, \Ql} }{ \ov{x} } \big)^{\on{Hf}} :=$
$$\{  c \in \restr{ \mc H_{G, N, I, W}^{j, \, \Ql} }{\ov{x}} , \; \text{the } \Oe \text{-module } C_{c}(K_{N}\backslash G(\mb A)/K_{N}, \Oe) \cdot c \text{ is of finite type} \} .$$
By Definition, 
\begin{equation}    \label{equation-H-Hf-includ-in-H-Hf-rat-Ql}
\big( \restr{ \mc H_{G, N, I, W}^{j,  \, \Ql} }{ \ov{x} } \big)^{\on{Hf-rat}} \supset \big( \restr{ \mc H_{G, N, I, W}^{j,  \, \Ql} }{ \ov{x} } \big)^{\on{Hf}}  .
\end{equation}

Thus (\ref{equation-H-cusp-equal-H-Hf-rat-Ql}) and (\ref{equation-H-Hf-includ-in-H-Hf-rat-Ql}) imply
\begin{equation}   \label{equation-H-cusp-contain-H-Hf-Ql}
\big( \restr{ \mc H_{G, N, I, W}^{j, \, \Ql}}{ \ov{x} } \big)^{\on{cusp}}  \supset \big( \restr{ \mc H_{G, N, I, W}^{j,  \, \Ql} }{ \ov{x} } \big)^{\on{Hf}}  .
\end{equation}

\begin{prop} \label{prop-H-cusp-equal-H-Hf} (see $loc.cit.$ Proposition 8.23 for the case $I = \emptyset$ and $W = \bf 1$)
(\ref{equation-H-cusp-contain-H-Hf-Ql}) is an equality
as sub-$\Ql$-vector spaces of $\restr{ \mc H_{G, N, I, W}^{j, \, \Ql} }{\ov{x}}$.
\end{prop}

\dem
It is enough to prove that 
\begin{equation}   \label{equation-H-Hf-contain-H-cusp-Ql}
\big( \restr{ \mc H_{G, N, I, W}^{j,  \, \Ql} }{ \ov{x} } \big)^{\on{cusp}} \subset \big( \restr{ \mc H_{G, N, I, W}^{j,  \, \Ql} }{ \ov{x} } \big)^{\on{Hf}}.
\end{equation}

As in \ref{subsection-Zl-lattice-in-W}, let $W^{\Oe} \in \on{Rep}_{\Oe}(\wh G)$ such that $W \simeq W^{\Oe} \otimes_{\Oe} \Ql$.
By Proposition \ref{prop-H-cusp-Zl-H-cusp-Ql}, we have $\big( \restr{ \mc H_{G, N, I, W^{\Oe}}^{j, \, \Oe} }{\ov{x}} \big)^{\on{cusp}} \otimes_{\Oe} \Ql \isom \big( \restr{ \mc H_{G, N, I, W  }^{j, \, \Ql} }{\ov{x}} \big)^{\on{cusp}}.$ By Theorem \ref{thm-cusp-dim-fini-second}, the $\Oe$-module $\big( \restr{ \mc H_{G, N, I, W^{\Oe}}^{j, \, \Oe} }{\ov{x}} \big)^{\on{cusp}}$ is of finite type. By Proposition \ref{prop-TC-commute-with-Hecke-with-niveau}, it is stable under the action of Hecke algebra $C_{c}(K_{N}\backslash G(\mb A)/K_{N}, \Oe)$. We deduce (\ref{equation-H-Hf-contain-H-cusp-Ql}).
\cqfd

\quad

\section{Partial Frobenius morphisms and the Eichler-Shimura relations}
\label{section-partial-Frob-ES-integral-coef}

In \ref{subsection-partial-Frob-coef-integral}, we recall the construction of the partial Frobenius morphisms for cohomology with coefficients in $\mc O_E$. In \ref{subsection-ES-relation-integral}-\ref{subsection-proof-ES-general}, we prove the Eichler-Shimura relations (Proposition \ref{prop-Eichler-Shimura-general-M}).

\subsection{Partial Frobenius morphisms}   \label{subsection-partial-Frob-coef-integral}

The goal of this subsection is for $W \in \on{Rep}_{\mc O_E}(\wh G^I)$, to equip $\mc H_{G, N, I, W}^{j, \, \mc O_E}$ with the partial Frobenius morphisms. For this, we need to construct a canonical morphism (\ref{equation-Frob-partial-F-G-N-I-W}).
First in \ref{subsection-Frob-partial-F-G-N-I-times-W-j} we treat the case $W = \boxtimes_{i \in I} W_i$ where each $W_i$ is a representation of $\wh G$ on a free $\mc O_E$-module of finite type. This case is the same as \cite{vincent} Section 3 but we prefer to write it down. Then in \ref{subsection-S-times-R-times-W}-\ref{subsection-Frob-partial-S-W-arbitaire} we treat the general case.

\begin{rem}  \label{rem-rep-O-E-not-semisimple}
By the functoriality coming from geometric Satake equivalence, for $W = W_1 \oplus W_2$ in $\on{Rep}_{E}(\wh G^I)$, we have $\mc F_{G, N, I, W} = \mc F_{G, N, I, W_1} \oplus \mc F_{G, N, I, W_2}$. For $E$-coefficients, the category $\on{Rep}_{E}(\wh G^I)$ is semisimple. Thus to construct the partial Frobenius morphisms, it is enough to consider $W$ irreducible. All $W$ irreducible are of the form $\boxtimes_{i \in I} W_i$ (up to enlarging $E$), thus it is enough to consider only the case $W = \boxtimes_{i \in I} W_i$. 

For $\mc O_E$-coefficients, the category $\on{Rep}_{\mc O_E}(\wh G^I)$ is no longer semisimple. Thus to construct the partial Frobenius morphisms, we need to consider general $W \in \on{Rep}_{\mc O_E}(\wh G^I)$.
\end{rem}

\sssec{}   \label{subsection-partial-Frob-semismall}
To construct the partial Frobenius morphisms, we need stacks of shtukas with intermediate modifications.
Let $k \in \N$ and $(I_1, \cdots, I_k)$ be a partition of $I$. Let $\Gr_{G, I}^{(I_1, \cdots, I_k)}$ be the ind-scheme defined in \cite{vincent} Définition 1.6 and $\Cht_{G, N, I}^{(I_1, \cdots, I_k)}$ be the stack defined in $loc.cit.$ Définition 2.1. When $k=1$, we have $\Gr_{G, I}^{(I)} = \Gr_{G, I}$ and $\Cht_{G, N, I}^{(I)} = \Cht_{G, N, I}$.

Section \ref{subsection-integral-coho-Cht-G} still holds if we replace $\Gr_{G, I}$ by $\Gr_{G, I}^{(I_1, \cdots, I_k)}$ and $\Cht_{G, N, I}$ by $\Cht_{G, N, I}^{(I_1, \cdots, I_k)}$. In particular, for any $W \in \on{Rep}_{\mc O_E}(\wh G^I)$, we have a $G_{I, \infty}$-equivariant perverse sheaf $\mc S_{G, I, W}^{(I_1, \cdots, I_k) , \, \mc O_E}$ (for the perverse normalization relative to $X^I$) on $\Gr_{G, I}^{(I_1, \cdots, I_k)}$ whose support is $\Gr_{G, I, W}^{(I_1, \cdots, I_k)}$. We have a perverse sheaf $\mc F_{G, N, I, W}^{(I_1, \cdots, I_k) , \, \mc O_E}$ (for the perverse normalization relative to $X^I$) on $\Cht_{G, N, I}^{(I_1, \cdots, I_k)}$ whose support is $\Cht_{G, N, I, W}^{(I_1, \cdots, I_k)}$.

We have the morphism forgetting the intermediate modifications:
$$
\xymatrixrowsep{1pc}
\xymatrixcolsep{3pc}
\xymatrix{
\Gr_{G, I, W}^{(I_1, \cdots, I_k)}   \ar[rr]^{ \pi^{(I_1, \cdots, I_k)}_{(I)}  }   \ar[rd]
& & \Gr_{G, I, W}^{(I)}   \ar[ld]  \\
& (X \sm N)^I
}
$$
Note that the morphism $\pi^{(I_1, \cdots, I_k)}_{(I)}$ is projective and small relatively to $(X \sm N)^I$.
We have a canonical isomorphism (cf. \cite{vincent} Théorème 1.17 b))
\begin{equation}
( \pi^{(I_1, \cdots, I_k)}_{(I)} )_! \mc S_{G, I, W}^{(I_1, \cdots, I_k) , \, \mc O_E} \simeq \mc S_{G, I, W}^{\mc O_E}
\end{equation}

We have a Cartesian diagram
$$
\xymatrixrowsep{2pc}
\xymatrixcolsep{5pc}
\xymatrix{
\Cht_{G, N, I}^{(I_1, \cdots, I_k)}   \ar[r]^{ \pi^{(I_1, \cdots, I_k)}_{(I)}  }   \ar[d]^{\epsilon^{ (I_1, \cdots, I_k)  }}
& \Cht_{G, N, I}^{(I)}   \ar[d]^{\epsilon^{ (I)  }}  \\
[G_{I, d} \backslash \Gr_{G, I, W}^{(I_1, \cdots, I_k)} ]   \ar[r]^{ \pi^{(I_1, \cdots, I_k)}_{(I)}  }  
& [G_{I, d} \backslash \Gr_{G, I, W}^{(I)}  ]
}
$$
By definition $\mc F_{G, N, I, W}^{(I), \, \mc O_E} = (\epsilon^{ (I) })^* \mc S_{G, I, W}^{(I), \, \mc O_E}$ and $\mc F_{G, N, I, W}^{(I_1, \cdots, I_k), \, \mc O_E} = (\epsilon^{ (I_1, \cdots, I_k)  })^* \mc S_{G, I, W}^{(I_1, \cdots, I_k), \, \mc O_E}$. By proper base change, we deduce 
\begin{equation}
( \pi^{(I_1, \cdots, I_k)}_{(I)} )_! \mc F_{G, I, W}^{(I_1, \cdots, I_k) , \, \mc O_E} \simeq \mc F_{G, I, W}^{\mc O_E}
\end{equation}
As a consequence, $\mc F_{G, I, W}^{(I_1, \cdots, I_k) , \, \mc O_E}$ and $\mc F_{G, I, W}^{\mc O_E}$ have the same direct image over $(X \sm N)^I$, i.e. define the same cohomology sheaf with compact support $\mc H_{G, N, I, W}^{j, \, \mc O_E}$.

\sssec{}   \label{subsection-def-Frob-I}
Recall the construction in \cite{vincent} Section 3.
We have a partial Frobenius morphisms 
\begin{equation}
\on{Fr}_{I_1}^{(I_1, \cdots, I_k)}: \Cht_{G, N, I}^{(I_1, \cdots, I_k)}  \rightarrow \Cht_{G, N, I}^{(I_2, \cdots, I_k, I_1)}
\end{equation}

We define $\Frob_{I_1}: (X \sm N)^I \rightarrow (X \sm N)^I$ by sending $(x_i)_{i \in I}$ to $(x_i')_{i \in I}$ with $x_i' = \Frob(x_i)$ if $i \in I_1$ and $x_i' = x_i$ if $i \notin I_1$.
The following diagram is commutative:
$$
\xymatrixrowsep{2pc}
\xymatrixcolsep{5pc}
\xymatrix{
\Cht_{G, N, I}^{(I_1, \cdots, I_k)}    \ar[r]^{ \on{Fr}_{I_1}^{(I_1, \cdots, I_k)}  }  \ar[d]
& \Cht_{G, N, I}^{(I_2, \cdots, I_k, I_1)} \ar[d] \\
(X \sm N)^I   \ar[r]^{ \Frob_{I_1}  }
& (X \sm N)^I
}
$$
where the vertical morphisms are the morphisms of paws.

Let $W \in \on{Rep}_{\mc O_E}(\wh G^I)$. We want to construct a canonical isomorphism
\begin{equation}    \label{equation-Frob-partial-F-G-N-I-W}
( \on{Fr}_{I_1}^{(I_1, \cdots, I_k)} )^* (  \mc F_{G, N, I, W}^{(I_2, \cdots, I_k, I_1), \, \mc O_E} )  \simeq \mc F_{G, N, I, W}^{(I_1, \cdots, I_k), \, \mc O_E}
\end{equation}

\sssec{}   \label{subsection-Frob-partial-F-G-N-I-times-W-j}
Let $W = \boxtimes_{j \in \{1, \cdots, k\}} W_j$, where each $W_j$ is a representation of $\wh G^{I_j}$ on a free $\mc O_E$-module of finite type. Consider the composition 
$$\wt{\epsilon}^{(I_1, \cdots, I_k)} : \Cht_{G, N, I, W}^{(I_1, \cdots, I_k)} \xrightarrow{ \epsilon^{(I_1, \cdots, I_k)}  } [ G_{I, d} \backslash \Gr_{G, I, W}^{(I_1, \cdots, I_k)}  ] \xrightarrow{ \wt{\kappa}^{(I_1, \cdots, I_k)}_I }
\prod_{j=1}^k [ G_{I_j, d} \backslash \Gr_{G, I_j, W_j}^{(I_j)} ]  $$
where the first morphism follows from \cite{vincent} (2.3) and the second morphism follows from $loc.cit.$ (1.13). The composition is $loc.cit.$ (2.5).
By definition we have
$$\mc F_{G, N, I, W}^{(I_1, \cdots, I_k), \, \mc O_E} = (\epsilon^{(I_1, \cdots, I_k)})^* \mc S_{G, I, W}^{(I_1, \cdots, I_k), \, \mc O_E} $$
By $loc.cit.$ Théorème 1.17 c) we have a canonical isomorphism
$$\mc S_{G, I, W}^{(I_1, \cdots, I_k), \, \mc O_E} \simeq (\wt{\kappa}^{(I_1, \cdots, I_k)}_I)^* ( \boxtimes_{j \in \{1, \cdots, k\} } \mc S_{G, I_j, W_j}^{(I_j), \, \mc O_E} )$$
We deduce a canonical isomorphism (see also $loc.cit.$ Corollaire 2.16)
\begin{equation}   \label{equation-F-1-k-product-S}
\mc F_{G, N, I, W}^{(I_1, \cdots, I_k), \, \mc O_E} \simeq (\wt{\epsilon}^{(I_1, \cdots, I_k)})^* ( \boxtimes_{j \in \{1, \cdots, k\} } \mc S_{G, I_j, W_j}^{(I_j), \, \mc O_E} )
\end{equation} 
Similarly we define 
$$\wt{\epsilon}^{(I_2, \cdots, I_k, I_1)} : \Cht_{G, N, I, W}^{(I_2, \cdots, I_k, I_1)} \rightarrow
\prod_{j=1}^k [ G_{I_j, d} \backslash \Gr_{G, I_j, W_j}^{(I_j)} ]  $$
and we have a canonical isomorphism
\begin{equation}    \label{equation-F-2-k-product-S}
\mc F_{G, N, I, W}^{(I_2, \cdots, I_k, I_1), \, \mc O_E} \simeq (\wt{\epsilon}^{(I_2, \cdots, I_k, I_1)})^* ( \boxtimes_{j \in \{1, \cdots, k\} } \mc S_{G, I_j, W_j}^{(I_j), \, \mc O_E} )
\end{equation}

As in $loc.cit.$ Proposition 3.4, we have a Cartesian diagram:
\begin{equation}   \label{equation-Cht-prod-Gr-Frob-partiel}
\xymatrixrowsep{2pc}
\xymatrixcolsep{6pc}
\xymatrix{
\Cht_{G, N, I, W}^{(I_1, \cdots, I_k)}    \ar[r]^{\on{Fr}_{I_1}^{(I_1, \cdots, I_k)}}  \ar[d]^{\wt{\epsilon}^{(I_1, \cdots, I_k)}}
& \Cht_{G, N, I, W}^{(I_2, \cdots, I_k, I_1)}    \ar[d]^{\wt{\epsilon}^{(I_2, \cdots, I_k, I_1)}}   \\
\prod_{j=1}^k [ G_{I_j, d} \backslash \Gr_{G, I_j, W_j}^{(I_j)} ]    \ar[r]^{\on{Frob_{1}} \times \Id}
& \prod_{j=1}^k [ G_{I_j, d}  \backslash  \Gr_{G, I_j, W_j}^{(I_j)} ] 
}
\end{equation}
where 
the lower line is the product of the total Frobenius morphism $\on{Frob}_1$ on $[ G_{I_1, d}  \backslash  \Gr_{G, I_1, W_1}^{(I_1)} ]$ and identity on the other $[ G_{I_j, d}  \backslash  \Gr_{G, I_j, W_j}^{(I_j)} ]$ for $j=2, \cdots, k$.

Since $[ G_{I_1, d} \backslash \Gr_{G, I_1, W_1}^{(I_1)} ]$ is defined over $\on{Spec} \Fq$, we have a canonical isomorphism 
\begin{equation}   \label{equation-Frob-S-to-S}
F_1: \on{Frob}_1^* ( \mc S_{G, I_1, W_1}^{(I_1), \, \mc O_E}  )  \simeq   \mc S_{G, I_1, W_1}^{(I_1), \, \mc O_E} 
\end{equation}

We deduce a canonical isomorphism
\begin{equation}   \label{equation-Frob-F-2-k-equal-F-1-k}
\begin{aligned}
( \on{Fr}_{I_1}^{(I_1, \cdots, I_k)} )^* (  \mc F_{G, N, I, W}^{(I_2, \cdots, I_k, I_1), \, \mc O_E} )  &  \overset{(\ref{equation-F-2-k-product-S})}{\simeq} ( \on{Fr}_{I_1}^{(I_1, \cdots, I_k)} )^* (\wt{\epsilon}^{(I_2, \cdots, I_k, I_1)})^* (\boxtimes_{j=1}^k \mc S_{G, I_j, W_j}^{(I_j), \, \mc O_E} ) \\
& \simeq (\wt{\epsilon}^{(I_1, \cdots, I_k)})^* (\on{Frob}_1 \times \Id) ^* ( \boxtimes_{j=1}^k \mc S_{G, I_j, W_j}^{(I_j), \, \mc O_E}  )  \\
& \overset{(\ref{equation-Frob-S-to-S})}{\simeq} (\wt{\epsilon}^{(I_1, \cdots, I_k)})^* ( \boxtimes_{j=1}^k \mc S_{G, I_j, W_j}^{(I_j), \, \mc O_E}  ) \\
& \overset{(\ref{equation-F-1-k-product-S})}{\simeq} \mc F_{G, N, I, W}^{(I_1, \cdots, I_k), \, \mc O_E}
\end{aligned}
\end{equation}

Moreover, since each isomorphism is canonical and functorial in each $W_j$, the composition (\ref{equation-Frob-F-2-k-equal-F-1-k}) is functorial in each $W_j$.

\sssec{}    \label{subsection-S-times-R-times-W}
Now let $W$ be an arbitrary representation of $\wh G^I$ on a $\mc O_E$-module of finite type (not necessarily free). 

As in \cite{vincent} proof of Théorème 1.17, let $\mc R$ be the algebra of regular functions on $\wh G$ with coefficients in $\mc O_E$, considered as an ind-object of $\on{Rep}_{\mc O_E}(\wh G)$ by the action of $\wh G$ on itself by the left multiplication. 
We view $\boxtimes_{i \in I} \mc R$ as the regular representation of $\wh G^I$ by the action of $\wh G^I$ on itself by the left multiplication. 
We define $\mc S_{G, I, \boxtimes_{i \in I} \mc R}^{(I_1, \cdots, I_k), \, \mc O_E}$ as an ind-object of $\Perv_{G_{I, \infty}}(\Gr_{G, I}, \mc O_E)$.

By \cite{mv} 4.1 (or \cite{BR} Lemma 1.10.9), the complex $\mc S_{G, I, \boxtimes_{i \in I} \mc R}^{(I_1, \cdots, I_k), \, \mc O_E} \overset{L}{\otimes}_{\mc O_E} W$ is perverse. Thus it is equal to $\mc S_{G, I, \boxtimes_{i \in I} \mc R}^{(I_1, \cdots, I_k), \, \mc O_E} \otimes_{\mc O_E} W$ in the abelian category of perverse sheaves.
The action of $\wh G^I$ on $\boxtimes_{i \in I} \mc R$ by the right multiplication induces an action of $\wh G^I$ on $\mc S_{G, I, \boxtimes_{i \in I} \mc R}^{(I_1, \cdots, I_k), \, \mc O_E}$ by the functoriality of the geometric Satake functor. We deduce a diagonal action of $\wh G^I$ on $\mc S_{G, I, \boxtimes_{i \in I} \mc R}^{(I_1, \cdots, I_k), \, \mc O_E} \otimes_{\mc O_E} W$. 
We have (cf. \cite{vincent} proof of Théorème 1.17)
\begin{equation}    \label{equation-S-W-equal-S-times-R-times-W-inv-by-G}
\mc S_{G, I, W}^{(I_1, \cdots, I_k), \, \mc O_E} = \big(  \mc S_{G, I, \boxtimes_{i \in I} \mc R}^{(I_1, \cdots, I_k), \, \mc O_E} \otimes_{\mc O_E} W \big)^{\wh G^I}  
\end{equation}


\sssec{}    \label{subsection-Frob-partial-S-W-arbitaire}
We denote by $\mc R_j:= \boxtimes_{i \in I_j} \mc R$. 
We have a Cartesian diagram
\begin{equation}   \label{equation-Cht-prod-Gr-Frob-partiel-R}
\xymatrixrowsep{2pc}
\xymatrixcolsep{6pc}
\xymatrix{
\Cht_{G, N, I, \boxtimes_{i \in I} \mc R}^{(I_1, \cdots, I_k)} \ar[r]^{\on{Fr}_{I_1}^{(I_1, \cdots, I_k)}}  \ar[d]^{\wt{\epsilon}^{(I_1, \cdots, I_k)}}
& \Cht_{G, N, I, \boxtimes_{i \in I} \mc R}^{(I_2, \cdots, I_k, I_1)} \ar[d]^{\wt{\epsilon}^{(I_2, \cdots, I_k, I_1)}} \\
\prod_{j=1}^k [ G_{I_j, \infty} \backslash \Gr_{G, I_j, \mc R_j}^{(I_j)} ]    \ar[r]^{\on{Frob_{1}} \times \Id}
& \prod_{j=1}^k [ G_{I_j, \infty} \backslash \Gr_{G, I_j, \mc R_j}^{(I_j)} ] 
}
\end{equation}
where the lower line is the product of the total Frobenius morphism $\on{Frob}_1$ on $[ G_{I_1, d}  \backslash  \Gr_{G, I_1, \mc R_1}^{(I_1)} ]$ and identity on the other $[ G_{I_j, d}  \backslash  \Gr_{G, I_j, \mc R_j}^{(I_j)} ]$ for $j=2, \cdots, k$.

We have a canonical isomorphism
\begin{equation}   \label{equation-Frob-total-S-R}
(\on{Frob}_1 \times \Id)^*(\boxtimes_{j \in \{1, \cdots, k\} } \mc S_{G, I_j, \mc R_j}^{(I_j), \, \mc O_E} ) \simeq   \boxtimes_{j \in \{1, \cdots, k\} } \mc S_{G, I_j, \mc R_j}^{(I_j), \, \mc O_E}  
\end{equation}
By \cite{vincent} Corollaire 2.16, 
\begin{equation}   \label{equation-F-R-equal-prod-S-R-i}   
\mc F_{G, N, I, \boxtimes_{i \in I} \mc R}^{(I_1, \cdots, I_k), \, \mc O_E} = (\wt{\epsilon}^{(I_1, \cdots, I_k)})^* ( \boxtimes_{j \in \{1, \cdots, k\} } \mc S_{G, I_j, \mc R_j}^{(I_j), \, \mc O_E} )
\end{equation}
$$
\mc F_{G, N, I, \boxtimes_{i \in I} \mc R}^{(I_2, \cdots, I_k, I_1), \, \mc O_E} = (\wt{\epsilon}^{(I_2, \cdots, I_k, I_1)})^* ( \boxtimes_{j \in \{1, \cdots, k\} } \mc S_{G, I_j, \mc R_j}^{(I_j), \, \mc O_E} )
$$
As in \ref{subsection-Frob-partial-F-G-N-I-times-W-j} we deduce a canonical isomorphism
\begin{equation}
(\on{Fr}_{I_1}^{(I_1, \cdots, I_k)})^* (  \mc F_{G, N, I, \boxtimes_{i \in I} \mc R}^{(I_2, \cdots, I_k, I_1), \, \mc O_E}   )  \simeq \mc F_{G, N, I, \boxtimes_{i \in I} \mc R}^{(I_1, \cdots, I_k), \, \mc O_E} 
\end{equation}

It induces a canonical isomorphism of ind-perverse sheaves:
\begin{equation}     \label{equation-Frob-partial-S-times-R-times-W}
\begin{aligned}
(\on{Fr}_{I_1}^{(I_1, \cdots, I_k)})^* (  \mc F_{G, N, I, \boxtimes_{i \in I} \mc R}^{(I_2, \cdots, I_k, I_1), \, \mc O_E} \otimes_{\mc O_E} W   )
& \simeq  \big(  (\on{Fr}_{I_1}^{(I_1, \cdots, I_k)})^* (  \mc F_{G, N, I, \boxtimes_{i \in I} \mc R}^{(I_2, \cdots, I_k, I_1), \, \mc O_E}  ) \big)  \otimes_{\mc O_E} W  \\
&  \simeq \mc F_{G, N, I, \boxtimes_{i \in I} \mc R}^{(I_1, \cdots, I_k), \, \mc O_E} \otimes_{\mc O_E} W
\end{aligned}
\end{equation}

By (\ref{equation-F-R-equal-prod-S-R-i}), (\ref{equation-Frob-total-S-R}) and the fact that the action of $\wh G^I$ commutes with the total Frobenius morphism, the isomorphism (\ref{equation-Frob-partial-S-times-R-times-W}) commutes with the action of $\wh G^I$ on $\mc F_{G, N, I, \boxtimes_{i \in I} \mc R}^{(I_1, \cdots, I_k), \, \mc O_E} \otimes_{\mc O_E} W$. Thus we have morphisms
\begin{equation}    \label{equation-Frob-partial-RGamma-S-times-R-times-W}
\begin{aligned}
& (\on{Fr}_{I_1}^{(I_1, \cdots, I_k)})^* \big( ( \mc F_{G, N, I, \boxtimes_{i \in I} \mc R}^{(I_2, \cdots, I_k, I_1), \, \mc O_E} \otimes_{\mc O_E} W )^{\wh G^I}  \big)  \\
 \simeq &    \big(  (\on{Fr}_{I_1}^{(I_1, \cdots, I_k)})^* ( \mc F_{G, N, I, \boxtimes_{i \in I} \mc R}^{(I_2, \cdots, I_k, I_1), \, \mc O_E} \otimes_{\mc O_E} W )   \big)^{\wh G^I}  \\
  \simeq &   (  \mc F_{G, N, I, \boxtimes_{i \in I} \mc R}^{(I_1, \cdots, I_k), \, \mc O_E} \otimes_{\mc O_E} W  )^{\wh G^I}
\end{aligned}
\end{equation}

We deduce from (\ref{equation-S-W-equal-S-times-R-times-W-inv-by-G}) and (\ref{equation-Frob-partial-RGamma-S-times-R-times-W}) a canonical isomorphism
\begin{equation}  \label{equation-Fr-I-1-F}
 (\on{Fr}_{I_1}^{(I_1, \cdots, I_k)})^* ( \mc F_{G, N, I, W}^{(I_2, \cdots, I_k, I_1), \, \mc O_E} )  \simeq   \mc F_{G, N, I, W}^{(I_1, \cdots, I_k), \, \mc O_E}
\end{equation}

\quad

\sssec{}   \label{subsection-Frob-partial-H-G-N-I-W}
Once have (\ref{equation-Frob-partial-F-G-N-I-W}) (which is also (\ref{equation-Fr-I-1-F})), taking into account \ref{subsection-partial-Frob-semismall}, as in \cite{vincent} Section 4.3, we deduce that there exists $\kappa$ depending on $W$, such that for any $\mu$ the cohomological correspondence induces a partial Frobenius morphism:
\begin{equation}
F_{I_1}: \on{Frob}_{I_1}^*(\mc H_{G, N, I, W}^{\leq \mu, \, \mc O_E}) \rightarrow \mc H_{G, N, I, W}^{\leq \mu + \kappa, \, \mc O_E} 
\end{equation}
where $\on{Frob}_{I_1}: X^I \rightarrow X^I$ is defined in \ref{subsection-def-Frob-I}.

Similarly, we have $\on{Frob}_{I_2}, \cdots, \on{Frob}_{I_k}$. The composition $\Frob_{I_k} \circ \cdots \circ \Frob_{I_1}$ is the total Frobenius morphism followed by an increase of the Harder-Narasimhan truncation $\on{Frob}^*(\mc H_{G, N, I, W}^{\leq \mu, \, \mc O_E}) \rightarrow \mc H_{G, N, I, W}^{\leq \mu + \kappa', \, \mc O_E} $.  


\sssec{}   \label{subsection-partial-Frob-functorial}
For any $W, W' \in \on{Rep}_{\mc O_E}(\wh G^I)$ and $W \xrightarrow{u} W'$, by functoriality it induces a morphism $\mc H_{G, N, I, W}^{\leq \mu, \, \mc O_E} \xrightarrow{u} \mc H_{G, N, I, W'}^{\leq \mu, \, \mc O_E}$. The following diagram is commutative:
$$
\xymatrixrowsep{2pc}
\xymatrixcolsep{5pc}
\xymatrix{
\on{Frob}_{I_1}^*(\mc H_{G, N, I, W}^{\leq \mu, \, \mc O_E})  \ar[r]^{ F_{I_1}  }  \ar[d]^{\on{Frob}_{I_1}^*(u)}
& \mc H_{G, N, I, W}^{\leq \mu + \kappa, \, \mc O_E}  \ar[d]^{u} \\
\on{Frob}_{I_1}^*(\mc H_{G, N, I, W'}^{\leq \mu, \, \mc O_E})   \ar[r]^{ F_{I_1}  }  
& \mc H_{G, N, I, W'}^{\leq \mu + \kappa, \, \mc O_E}  
}
$$
The reason is that the following diagram is commutative
$$
\xymatrixrowsep{2pc}
\xymatrixcolsep{5pc}
\xymatrix{
(\on{Fr}_{I_1}^{(I_1, \cdots, I_k)})^* ( \mc F_{G, N, I, W}^{(I_2, \cdots, I_k, I_1), \, \mc O_E} )   \ar[r]^{ \quad \quad \simeq }  \ar[d]^{\on{Frob}_{I_1}^*(u)}
&  \mc F_{G, N, I, W}^{(I_1, \cdots, I_k), \, \mc O_E} \ar[d]^{u} \\
(\on{Fr}_{I_1}^{(I_1, \cdots, I_k)})^* ( \mc F_{G, N, I, W'}^{(I_2, \cdots, I_k, I_1), \, \mc O_E} )  \ar[r]^{ \quad \quad \simeq }  
&  \mc F_{G, N, I, W'}^{(I_1, \cdots, I_k), \, \mc O_E}
}
$$
In fact, when $W = \boxtimes_{j \in \{1, \cdots, k\}} W_j$, 
we have (\ref{equation-F-1-k-product-S}) and (\ref{equation-F-2-k-product-S}).
The total Frobenius morphism (\ref{equation-Frob-S-to-S}) is canonical and functorial on $W$. More generally in \ref{subsection-Frob-partial-S-W-arbitaire}, (\ref{equation-Frob-partial-S-times-R-times-W}) is functorial in $W$.

\begin{rem}
The commutativity in \ref{subsection-partial-Frob-functorial} does not come from a general property. It comes from special property of Satake sheaves. 

In general, if $\mc F$ (resp. $\mc G$) is a sheaf over $X^I$ equipped with partial Frobenius morphisms $f_i: \on{Frob}_{\{i\}}^* \mc F \isom \mc F$ (resp. $g_i: \on{Frob}_{\{i\}}^* \mc G \isom \mc G$) commuting with each other and the composition is the total Frobenius morphism, and if we have a morphism $\mc F \xrightarrow{u} \mc G$, the diagram 
$$
\xymatrixrowsep{2pc}
\xymatrixcolsep{5pc}
\xymatrix{
\on{Frob}_{\{i\}}^*\mc F  \ar[r]^{ f_i  }_{\sim}  \ar[d]^{\on{Frob}_{\{i\}}^*(u)}
& \mc F  \ar[d]^{u} \\
\on{Frob}_{\{i\}}^* \mc G \ar[r]^{ g_i  }_{\sim}  
& \mc G 
}
$$
may not be commutative. 
However the diagram with total Frobenius is always commutative:
$$
\xymatrixrowsep{2pc}
\xymatrixcolsep{5pc}
\xymatrix{
\on{Frob}^*\mc F  \ar[r]^{ \simeq }  \ar[d]^{\on{Frob}^*(u)}
& \mc F  \ar[d]^{u} \\
\on{Frob}^* \mc G \ar[r]^{ \simeq }  
& \mc G 
}
$$
\end{rem}

\quad

\subsection{Eichler-Shimura relation}    \label{subsection-ES-relation-integral}

\sssec{}
For cohomology with coefficients in $E$, the Eichler-Shimura relation in \cite{vincent} Proposition 7.1 is enough for the applications: the category $\on{Rep}_{E}(\wh G^I)$ is semisimple so it is enough to consider only $W \in \on{Rep}_{E}(\wh G^I)$ of the form $W = \boxtimes_{i \in I} W_i$.

For cohomology with coefficients in $\mc O_E$ and $W = \boxtimes_{i \in I} W_i$ with $W_i \in \on{Rep}_{\mc O_E}(\wh G)^{\on{free}}$, this proposition still holds. However, the category $\on{Rep}_{\mc O_E}(\wh G^I)$ is no longer semisimple (see Remark \ref{rem-rep-O-E-not-semisimple}) so it is not enough to consider only $W = \boxtimes_{i \in I} W_i$. We need the Eichler-Shimura relation for more general $W \in \on{Rep}_{\mc O_E}(\wh G^I)$, which is Proposition \ref{prop-Eichler-Shimura-general-M} below. Let us begin by recall some definitions.

\begin{nota}
In the following, we write $\mc H_{I, W}$ for $\mc H_{G, N, I, W}^{\mc O_E}$ to shorten the notation.
\end{nota}

\sssec{}   \label{subsection-partial-Frob-i}
For any $i \in I$, apply \ref{subsection-Frob-partial-H-G-N-I-W} to $I_1 = \{i\}$. We have the partial Frobenius morphism over $(X \sm N)^I$
$$F_{\{i\}}: \Frob_{\{i\}}^* \mc H_{I, W}^{\leq \mu} \rightarrow \mc H_{I, W}^{\leq \mu + \kappa}$$
Let $v$ be a place of $X \sm N$. Note that $\Frob^{\on{deg}(v)}(v) = v$. Restrict $F_{\{i\}}^{\on{deg}(v)}$ to $(X \sm N)^{I-\{i\}} \times v$, we obtain a morphism
\begin{equation}
F_{\{i\}}^{\on{deg}(v)}: \restr{ \mc H_{I, W}^{\leq \mu} }{ (X \sm N)^{I-\{i\}} \times v  } \rightarrow \restr{ \mc H_{I, W}^{\leq \mu + \on{deg}(v) \kappa} }{ (X \sm N)^{I-\{i\}} \times v }
\end{equation}

\sssec{}   \label{subsection-def-S-V-v}
The construction in \cite{vincent} Section 6.1 still works for cohomology with coefficients in $\mc O_E$. Let's recall the main construction. Let $v$ be a place of $X \sm N$. We denote by $\on{Rep}_{\mc O_E}(\wh G)^{\on{free}}$ the category of representations of $\wh G$ on a free $\mc O_E$-module of finite type. Let $V \in \on{Rep}_{\mc O_E}(\wh G)^{\on{free}}$. We will construct an operator $S_{V, v}$.

%
%
%

Let $V^*$ be the dual representation of $V$. We denote by $V \boxtimes V^*$ the representation of $\wh G \times \wh G$ and by $V \otimes V^*$ the diagonal representation of $\wh G$. 
We have morphisms of $\wh G$-representations:
$$\delta_V: {\bf 1} \rightarrow V \otimes V^*, \quad 1 \mapsto \sum_j e_j \otimes e_j^*$$
$$\on{ev}_V: V \otimes V^* \rightarrow {\bf 1}, \quad
 x \otimes \xi \mapsto \xi(x)$$
where $\{e_j\}_j$ is a basis of $V$ and $\{e_j^*\}_j$ its dual basis in $V^*$. Let $\Delta: X \rightarrow X \times X$ be the diagonal morphism.

As in \cite{vincent} Section 6.1, for any finite set $I$ and $W \in \on{Rep}_{\mc O_E}(\wh G^{I})$, let $S_{V, v}$ be the composition of morphisms
$${  \resizebox{15cm}{!}{ 
\!\!\!\!\!\!\!\!\!\!\!
\xymatrixrowsep{2pc}
\xymatrixcolsep{5pc}
\xymatrix{
(\mc O_E)_v \boxtimes \mc H_{I, W}^{\leq \lambda}  \simeq\restr{\mc H_{\{*\} \cup I, {\bf 1} \boxtimes W}^{\leq \lambda} }{v \times (X \sm N)^I}  \ar[r]^{\quad \quad  \delta_{V} \boxtimes \Id_W} 
& \restr{ \mc H_{\{*\} \cup I,  (V \otimes V^*)\boxtimes W}^{\leq \lambda} }{v \times (X \sm N)^I}   \ar[r]_{\simeq  }^{  \chi_{\zeta_{\{1, 2\}}}^{-1} }
 & \restr{ \mc H_{\{1, 2\} \cup I, V \boxtimes V^* \boxtimes W}^{\leq \lambda}  }{\Delta(v) \times (X \sm N)^I} \ar[d]^{F_{\{1\}}^{\on{deg}(v)}} \\
(\mc O_E)_v   \boxtimes  \mc H_{I, W}^{\leq \lambda + \kappa}   
\simeq \restr{\mc H_{\{*\} \cup I, {\bf 1} \boxtimes W}^{\leq \lambda + \kappa} }{v \times (X \sm N)^I}    
& \restr{ \mc H_{\{*\} \cup I,  (V \otimes V^*)\boxtimes W}^{\leq \lambda + \kappa} }{v \times (X \sm N)^I}  \ar[l]_{ \quad \quad \on{ev}_{V} \boxtimes \Id_W}
&  \restr{ \mc H_{\{1, 2\} \cup I, V \boxtimes V^* \boxtimes W}^{\leq \lambda + \kappa}  }{\Delta(v) \times (X \sm N)^I}  \ar[l]^{\simeq}_{ \chi_{\zeta_{\{1, 2\}}}^{-1}  }
}
} }$$
where $\chi_{\zeta_{\{1, 2\}}}$ is the fusion isomorphism (cf. \cite{vincent} Proposition 4.12) associated to 
$\zeta_{\{1, 2\}}: \{1, 2\} \rightarrow \{*\}$.
As in $loc.cit.$ Section 6.1, the composition commutes with the partial Frobenius morphism on $(\mc O_E)_v$, thus $S_{V, v}$ descends to a morphism:
\begin{equation}    \label{equation-S-V-v-H-G-N-I-W}
S_{V, v}: \mc H_{I, W}^{\leq \lambda} \rightarrow \mc H_{I, W}^{\leq \lambda + \kappa}  \quad \text{ in }  D_c^{(-)}((X \sm N)^{I}, \mc O_E)
\end{equation}

\sssec{}
Let $M \in \on{Rep}_{\mc O_E}(\wh G)^{\on{free}}$. We define $\wedge^j M$ to be the quotient of $\otimes^j M$ by all $\otimes_{i \in \{1, 2, \cdots, j\}} x_i$ where two of $x_i$'s are equal. If $\{ e_i \}_i$ is a basis of $M$, then $\wedge^j M$ is free as $\mc O_E$-module with basis $\{ e_{i_1} \wedge \cdots \wedge e_{i_j}, \text{ where } i_1 < \cdots < i_j \}$.

\begin{prop}  (see \cite{vincent} Proposition 7.1 for the version of coefficients in $E$)  \label{prop-Eichler-Shimura-general-M}
For any finite set $I = \wt I \cup \{0\}$ and $W \in \on{Rep}_{\mc O_E}(\wh G^I)$,
there exists $M \in \on{Rep}_{\mc O_E}(\wh G)^{\on{free}}$, such that
\begin{equation}  \label{equation-ES-coef-O-E}
\sum_{\alpha=0}^{\on{rk}M} (-1)^{\alpha} S_{\wedge^{\on{rk} M - \alpha} M, v} \circ (F_{\{0\}}^{\on{deg}(v)})^{\alpha} =0 \quad \text{ in }
\end{equation}
$$\on{Hom}_{D^{(-)}_c( (X \sm N)^{\wt I} \times v , \mc O_E)}(\restr{\mc H_{I, W}^{\leq \mu}}{(X \sm N)^{\wt I} \times v}, \restr{\mc H_{I, W}^{\leq \mu + \kappa}}{(X \sm N)^{\wt I} \times v}).$$
\end{prop}

The proof will be given in \ref{subsection-proof-ES-singleton} and \ref{subsection-proof-ES-general}.

\sssec{}
In particular, when $W$ is of the form $W= \boxtimes_{i \in I} W_i$ with $W_i$ free, for each $i \in I$, we can write $I = \wt I \cup \{i\}$ and take $M = W_i$. We recover \cite{vincent} Proposition 7.1.

\quad

\begin{rem}
We call Proposition \ref{prop-Eichler-Shimura-general-M} the Eichler-Shimura relation because the operators $S_{V, v}$ extend the Hecke operators (see Proposition \ref{prop-S-equal-T-O-E} below). But the proof of Proposition \ref{prop-Eichler-Shimura-general-M} does not need this fact.
\end{rem}

\sssec{}
Let $v$ be a place of $X \sm N$ and $\ms H_{G, v}^{\mc O_E}:=C_c(G(\mc O_v) \backslash G(F_v) / G(\mc O_v), \mc O_E)$ the Hecke algebra with coefficients in $\mc O_E$. 
Let $V \in \on{Rep}_{\mc O_E}(\wh G)^{\on{free}}$. Let $h_{V, v} \in \ms H_{G, v}^{\mc O_E}$ be the function corresponding to $V$ via the Satake isomorphism. \cite{vincent} Section 4.4 still works for coefficients in $\mc O_E$. We have the Hecke operator defined by Hecke correspondence
$$T(h_{V, v}) : \restr{\mc H_{I, W}^{\leq \mu, \, \mc O_E} }{ (X \sm (N \cup v))^I} \rightarrow \restr{ \mc H_{I, W}^{\leq \mu+\kappa, \, \mc O_E}}{(X \sm (N \cup v))^I }$$ 
$loc.cit.$ Section 6 works for coefficients in $\mc O_E$. In particular, we have
\begin{prop} (\cite{vincent} Proposition 6.2) \label{prop-S-equal-T-O-E}
The operator $S_{V, v}$ defined in Section \ref{subsection-def-S-V-v}, which is a morphism of sheaves over $(X \sm N)^I$, extends the Hecke operator $T(h_{V, v})$, which is a morphism of sheaves over $(X \sm (N \cup v))^I$.
\cqfd
\end{prop}

\quad

\begin{rem}
For $E$-coefficients, there are two proofs of Proposition \ref{prop-Eichler-Shimura-general-M}. The one in \cite{vincent} Section 7 does not work for $\mc O_E$-coefficients (because that the denominator may not be invertible in $\mc O_E$). The one in \cite{xiao-zhu} Section 6 works for $\mc O_E$-coefficients. 

In \ref{subsection-proof-ES-singleton} and \ref{subsection-proof-ES-general} below, we will recall the construction in \cite{xiao-zhu} Section 6.2. There is no new idea here, everything is already contained in $loc.cit.$
\end{rem}

\sssec{}
The idea of the proof of Proposition \ref{prop-Eichler-Shimura-general-M} is that we will construct a morphism
$$\Hom_{ \wh G^I }(W, \mc R \otimes W)  \xrightarrow{\Theta} \on{End}_{D^{(-)}_c( (X \sm N)^{\wt I} \times v , \mc O_E)}(\restr{\mc H_{I, W}}{(X \sm N)^{\wt I} \times v})$$
On the left hand side, there is the Cayley-Hamilton theorem. This will give (\ref{equation-ES-coef-O-E}) on the right hand side.

We will first trait the case when $I$ is a singleton. Then the general case.

\begin{rem}
To simplify the notation, we will not write down the HN truncations. Just keep in mind that the partial Frobenius morphisms augment the HN truncations, and the morphisms induced by functoriality preserve the HN truncations.
\end{rem}

\subsection{Proof of Proposition \ref{prop-Eichler-Shimura-general-M}: case $I$ singleton}   \label{subsection-proof-ES-singleton}

\sssec{}    \label{subsection-def-R-otimes-W-singleton}
Recall that $\mc R$ (defined in \ref{subsection-S-times-R-times-W}) is the algebra of regular functions on $\wh G$ with coefficients in $\mc O_E$. It is equipped with an action of $\wh G^{\{r\}} \times \wh G^{\{l\}}$, where $\wh G^{\{r\}}$ acts by the right action and $\wh G^{\{l\}}$ acts by the left action:
$$
\wh G^{\{r\}} \times \wh G^{\{l\}} \times \mc R \rightarrow \mc R, \quad
\big( (h_1, h_2), f \big) \mapsto [g \mapsto f(h_2^{-1} g h_1)]
$$


Let $I = \{0\}$ be a singleton. Let $W \in \on{Rep}_{\mc O_E}(\wh G)$.
We equip $\mc R \boxtimes W$ with an action of $\wh G^{\{r\}} \times \wh G^{\{l\}} \times \wh G^{\{0\}}$, where $\wh G^{\{r\}} \times \wh G^{\{l\}}$ acts on $\mc R$ and $\wh G^{\{0\}}$ acts on $W$. In other words, if we view $\mc R \boxtimes W$ as the algebra of regular functions on $\wh G$ with coefficients in $W$, then this action is 
$$
\begin{aligned}
\wh G^{\{r\}} \times \wh G^{\{l\}} \times \wh G^{\{0\}} \times (\mc R \boxtimes W) & \rightarrow  \mc R \boxtimes W, \\
\big( (h_1, h_2, h), f: \wh G \rightarrow W \big) & \mapsto  [g \mapsto h f(h_2^{-1} g h_1)]
\end{aligned}
$$

We equip $\mc R \otimes W$ with a diagonal action of $\wh G$, where the action of $\wh G$ on $\mc R$ is by conjugation. In other words, if we view $\mc R \otimes W$ as the algebra of regular functions on $\wh G$ with coefficients in $W$, then this action of $\wh G$ is 
$$
\begin{aligned}
\wh G \times (\mc R \otimes W) & \rightarrow  \mc R \otimes W, \\
\big( h, f: \wh G \rightarrow W \big) & \mapsto  [g \mapsto h f(h^{-1} g h)]
\end{aligned}
$$

\begin{construction}    \label{construction-Theta-I-singleton}
For any $V_1, V_2 \in \on{Rep}_{\mc O_E}(\wh G)$,
we construct a morphism
\begin{equation}
\Hom_{ \wh G }(V_1, \mc R \otimes V_2)  \xrightarrow{\Theta} \on{Hom}_{D^{(-)}_c(  v , \mc O_E)}(\restr{\mc H_{\{0\}, V_1}}{ v}, \restr{\mc H_{\{0\}, V_2}}{ v})
\end{equation}
in the following way, where the action of $\wh G$ on $\mc R \otimes V_2$ is defined in \ref{subsection-def-R-otimes-W-singleton}.
For any $u \in \Hom_{ \wh G }(V_1, \mc R \otimes V_2) $, let $\Theta(u)$ be the composition of morphisms
\begin{equation}    \label{equation-Theta-u-I-singleton}
\xymatrixrowsep{2pc}
\xymatrixcolsep{4pc}
\xymatrix{
\restr{\mc H_{\{0\}, V_1}}{ v}   \ar[r]^{u}
& \restr{\mc H_{\{0\}, \mc R \otimes V_2}}{ v}  \ar[r]_{\simeq \quad \quad}^{\chi_{\zeta}^{-1} \quad \quad }
& \restr{\mc H_{\{r, l, 0\} , \mc R \boxtimes V_2}}{\Delta(v) }  \ar[d]^{F_{\{r\}}^{\on{deg}(v)}} \\
\restr{\mc H_{ \{0\}, V_2}}{ v}   
& \restr{\mc H_{ \{0\}, \mc R \otimes V_2}}{ v}   \ar[l]_{\on{ev}_{\mc R} \otimes \Id_{V_2}}
& \restr{\mc H_{\{r, l, 0\} , \mc R \boxtimes V_2}}{\Delta(v) }   \ar[l]^{\simeq \quad \quad }_{\chi_{\zeta} \quad \quad}
}
\end{equation}
where 
\begin{itemize}
\item the $\wh G$-representation $\mc R \otimes V_2$ and the $\wh G \times \wh G \times \wh G$-representation $\mc R \boxtimes V_2$ are defined in \ref{subsection-def-R-otimes-W-singleton}.

\item $\Delta: X^{\{0\}} \hookrightarrow X^{\{r, l, 0\}}$ is the diagonal morphism, $\chi_{\zeta}$ is the fusion isomorphism (cf. \cite{vincent} Proposition 4.12) associated to $\zeta: \{r, l, 0\} \rightarrow \{0\}, $

\item $\on{ev}_{\mc R}: \mc R \rightarrow {\bf 1}, f \mapsto f(1)$

\item $F_{\{r\}}$ is the partial Frobenius morphism defined in \ref{subsection-partial-Frob-i}
\end{itemize}
\end{construction}

\sssec{}   \label{subsection-Theta-linear-compatible-with-composition-singleton}
The morphism $\Theta$ is $\mc O_E$-linear (by the functoriality of $\mc H_{I, W}$ on $W$). 
The morphism $\Theta$ is compatible with composition: let $V_3 \in \on{Rep}_{\mc O_E}(\wh G)$, $u \in \Hom_{ \wh G }(V_1, \mc R \otimes V_2)$ and $u' \in \Hom_{ \wh G }(V_2, \mc R \otimes V_3)$. We define 
$$u'': V_1 \xrightarrow{u} \mc R \otimes V_2 \xrightarrow{\Id_{\mc R} \otimes u'} \mc R \otimes \mc R \otimes V_3 \xrightarrow{\on{res} \otimes \Id_{V_3}} \mc R \otimes V_3$$
where $\mc R \otimes \mc R \xrightarrow{\on{res}} \mc R$ is the diagonal restriction, which sends $f_1 \otimes f_2$ to $f_1f_2$.
The morphism $u''$ is $\wh G$-equivariant.
We have $\Theta(u') \circ \Theta(u) = \Theta(u'')$.
In fact, the following diagram is commutative
$${  \resizebox{15cm}{!}{ 
\!\!\!\!\!\!\!\!\!\!\!
\xymatrix{
\mc H_{\{0\}, V_1}   \ar[r]^{u}  \ar[dd]^{u''}
& \mc H_{\{0\}, \mc R \otimes V_2}  \ar[r]^{\simeq}  \ar[d]^{u'}
& \restr{ \mc H_{\{r, l, 0\}, \mc R \boxtimes V_2} }{\Delta}  \ar[r]^{F_{\{r\}}^{\deg v}}  \ar[d]^{u'}
& \restr{ \mc H_{\{r, l, 0\}, \mc R \boxtimes V_2} }{\Delta}  \ar[r]^{\simeq}  \ar[d]^{u'}
& \mc H_{\{0\}, \mc R \otimes V_2}  \ar[r]^{\on{ev}_{\mc R}}   \ar[d]^{u'}
& \mc H_{\{0\}, V_2}  \ar[d]^{u'} \\
& \mc H_{\{0\}, \mc R \otimes (\mc R \otimes V_3)}    \ar[r]^{\simeq}   \ar[d]^{\simeq}   \ar[dl]^{\on{res}}
& \restr{ \mc H_{\{r, l, 0\}, \mc R \boxtimes (\mc R \otimes V_3)} }{\Delta}  \ar[r]^{F_{\{r\}}^{\deg v}}    \ar[d]^{\simeq}
& \restr{ \mc H_{\{r, l, 0\}, \mc R \boxtimes (\mc R \otimes V_3)} }{\Delta}  \ar[r]^{\simeq}    \ar[d]^{\simeq}
& \mc H_{\{0\}, \mc R \otimes (\mc R \otimes V_3)}  \ar[r]^{\on{ev}_{\mc R}}    \ar[d]^{\simeq}
& \mc H_{\{0\}, \mc R \otimes V_3}    \ar[d]^{\simeq}   \\
\mc H_{\{0\}, \mc R \otimes V_3}  \ar[d]^{\simeq}
& \restr{ \mc H_{\{r', l', 0\}, (\mc R \otimes \mc R) \boxtimes V_3} }{\Delta}   \ar[r]^{\simeq}  \ar[ld]^{\on{res}}   \ar[rd]|-{F_{\{r'\}}^{\deg v}} 
& \restr{ \mc H_{\{r, l, r', l',  0\}, \mc R \boxtimes \mc R \boxtimes V_3} }{\Delta}  \ar[r]^{F_{\{r\}}^{\deg v}}   \ar[rd]|-{F_{\{r, r'\}}^{\deg v}} 
& \restr{ \mc H_{\{r, l, r', l', 0\}, \mc R \boxtimes \mc R \boxtimes V_3} }{\Delta}  \ar[r]^{\simeq}  \ar[d]^{F_{\{r'\}}^{\deg v}} 
& \restr{\mc H_{\{r', l', 0\}, (\mc R \otimes \mc R) \boxtimes V_3} }{\Delta} \ar[r]^{\on{ev}_{\mc R}}  \ar[d]
& \restr{ \mc H_{\{r', l', 0\}, \mc R \boxtimes V_3}  }{\Delta}  \ar[d]^{F_{\{r'\}}^{\deg v}} \\
\restr{ \mc H_{\{r', l', 0\}, \mc R \boxtimes V_3} }{\Delta}   \ar[rd]|-{F_{\{r'\}}^{\deg v}} 
& 
& \restr{ \mc H_{\{r', l', 0\}, (\mc R \otimes \mc R) \boxtimes V_3} }{\Delta}  \ar[r]^{\simeq}  \ar[ld]^{\on{res}}  \ar[rrd]^{\simeq}
& \restr{ \mc H_{\{r, l, r', l', 0\}, \mc R \boxtimes \mc R \boxtimes V_3} }{\Delta}  \ar[r]^{\simeq}
& \restr{\mc H_{\{r', l', 0\}, (\mc R \otimes \mc R) \boxtimes V_3} }{\Delta} \ar[r]^{\on{ev}_{\mc R}}   \ar[d]^{\simeq}
& \restr{ \mc H_{\{r', l', 0\}, \mc R \boxtimes V_3}  }{\Delta}  \ar[d]^{\simeq}  \\
& \restr{ \mc H_{\{r', l', 0\}, \mc R \boxtimes V_3} }{\Delta}  \ar[rrd]^{\simeq}
& 
& 
& \mc H_{\{0\}, \mc R \otimes \mc R \otimes V_3}  \ar[r]^{\on{ev}_{\mc R}}   \ar[dl]^{\on{res}} 
& \mc H_{\{0\}, \mc R \otimes V_3}  \ar[d]^{\on{ev}_{\mc R}}  \\
& 
& 
& \mc H_{\{0\}, \mc R \otimes V_3}  \ar[rr]^{\on{ev}_{\mc R}} 
& 
& \mc H_{\{0\}, V_3} 
} } }$$
where every cohomology with one leg (resp. three legs, five legs) is restricted on $v$ (resp. $\Delta(v)$ for $\Delta: X \rightarrow X^3$ diagonal morphism, $\Delta(v)$ for $\Delta: X \rightarrow X^5$ diagonal morphism). The composition of the upper horizontal line is $\Theta(u)$. The composition of the right vertical line is $\Theta(u')$. The composition of the left vertical and lower horizontal line is $\Theta(u'')$.

In the diagram, $\wh G^{\{r'\}} \times \wh G^{\{l'\}} \times \wh G^{\{0\}}$ acts on $(\mc R \otimes \mc R) \boxtimes V_3$ by $\wh G^{\{r'\}}$ (resp. $\wh G^{\{l'\}}$) acting diagonally on $\mc R \otimes \mc R$ by the right (resp. left) action and $\wh G^{\{0\}}$ acting on $V_3$.
The action of 
$\wh G^{\{r\}} \times \wh G^{\{l\}} \times \wh G^{\{r'\}} \times \wh G^{\{l'\}} \times \wh G^{\{0\}}$ on $\mc R \boxtimes \mc R \boxtimes V_3$ is by $\wh G^{\{r\}} \times \wh G^{\{l\}}$ acting on the first $\mc R$, $\wh G^{\{r'\}} \times \wh G^{\{l'\}}$ acting on the second $\mc R$ and $\wh G^{\{0\}}$ acting on $V_3$. 

%

\sssec{}
When $V_1 = V_2 = W \in \on{Rep}_{\mc O_E}(\wh G)$, Construction \ref{construction-Theta-I-singleton} gives a homomorphism
\begin{equation}
\Hom_{ \wh G }(W, \mc R \otimes W)  \xrightarrow{\Theta} \on{End}_{D^{(-)}_c( v , \mc O_E)}(\restr{\mc H_{\{0\}, W}}{ v})
\end{equation}

\sssec{}   \label{subsection-def-Tr-V}
For any $f \in \mc R$ invariant under the conjugation action of $\wh G$, it induces a map $f \otimes \Id_W: W \rightarrow \mc R \otimes W$ by $x \mapsto [g \mapsto f(g)x]$. 

In particular, let $V \in \on{Rep}_{\mc O_E}(\wh G)^{\on{free}}$ and $\on{Tr}_{V} \in \mc R$ be the function which sends $g \in \wh G$ to its trace $\on{Tr}_{V}(g)$. It is invariant under the conjugation action of $\wh G$. We have $\on{Tr}_{V} \otimes \Id_W \in \Hom_{ \wh G }(W, \mc R \otimes W) $.

\begin{lem}    \label{lem-Theta-Tr-V-equal-S-V-singleton}
We have 
\begin{equation}    \label{equation-Theta-Tr-V-equal-S-V-singelton}
\Theta(\on{Tr}_{V} \otimes \Id_W) = S_{V, v}.
\end{equation}
where $S_{V, v}$ is defined in \ref{subsection-def-S-V-v}, here we restrict it from $\mc H_{\{0\}, W}$ (which is over $X \sm N$) to $\restr{\mc H_{\{0\}, W}}{v}$.
\end{lem}
\dem
Let $V^*$ be the dual representation of $V$. We have a morphism of $\wh G^{\{r\}} \times \wh G^{\{l\}}$-representations:
$$m: V \boxtimes V^* \rightarrow \mc R, \quad (x , \xi) \mapsto [g \mapsto \xi(gx)]$$
where $\wh G^{\{r\}}$ acts on $V$ and $\wh G^{\{l\}}$ acts on $V^*$. 

We now view $m: V \otimes V^* \rightarrow \mc R$ as a morphism of $\wh G$-representations, where $\wh G$ acts via the diagonal morphism $\wh G \hookrightarrow \wh G^{\{r\}} \times \wh G^{\{l\}}$.

The following diagram is commutative
\begin{equation}    \label{equation-Theta-Tr-V-equal-S-V-singleton}
\xymatrixrowsep{2pc}
\xymatrixcolsep{5pc}
\xymatrix{
\restr{\mc H_{\{0\}, W}}{v }  \ar[d]^{\on{Tr}_V \otimes \Id_W \quad}
 \ar@{}[dr]|{(a)}  
& \restr{\mc H_{\{0\}, W}}{v}  \ar[d]^{\delta_V \otimes \Id_W \quad \quad }   \ar[l]^{=} \\
\restr{\mc H_{\{0\}, \mc R \otimes W}}{v }    \ar[d]_{\simeq }^{\chi_{\xi}^{-1} }
& \restr{\mc H_{\{0\}, ( V \otimes V^*) \otimes W}}{v }  \ar[d]_{\simeq  }^{\chi_{\xi}^{-1}  }   \ar[l]^{m \otimes \Id_W} \\
\restr{\mc H_{\{r, l, 0\} , \mc R \boxtimes W}}{\Delta^{\{r, l, 0\}}(v) }   \ar[d]^{F_{\{r\}}^{\on{deg}(v)}}   \ar@{}[dr]|{(b)}
& \restr{\mc H_{\{r, l, 0\} , ( V \boxtimes V^*) \boxtimes W}}{\Delta^{\{r, l, 0\}}(v) }  \ar[d]^{F_{\{r\}}^{\on{deg}(v)}}   \ar[l]^{m \boxtimes \Id_W}  \\
\restr{\mc H_{ \{r, l, 0\} , \mc R \boxtimes W}}{\Delta^{\{r, l, 0\}}(v) }   \ar[d]_{\simeq}^{\chi_{\xi}}
& \restr{\mc H_{ \{r, l, 0\} , ( V \boxtimes V^*) \boxtimes W}}{\Delta^{\{r, l, 0\}}(v) }  \ar[l]^{m \boxtimes \Id_W} \ar[d]_{ \simeq}^{\chi_{\xi}} \\
\restr{\mc H_{\{0\}, \mc R \otimes W}}{v }    
\ar[d]^{\on{ev}_{\mc R} \otimes \Id_W}   \ar@{}[dr]|{(c)}
& \restr{\mc H_{\{0\}, ( V \otimes V^*) \otimes W} }{v }    \ar[d]^{\on{ev}_V \otimes \Id_W}   \ar[l]^{m \otimes \Id_W} \\
\restr{\mc H_{\{0\}, W}}{v }  
& \restr{\mc H_{\{0\}, W}}{v }   \ar[l]^{=} 
}
\end{equation} 
In this diagram, $(a)$ commutes because the following diagram of $\wh G$-representations is commutative:
\begin{equation}
\xymatrixrowsep{2pc}
\xymatrixcolsep{5pc}
\xymatrix{
{\bf 1} \otimes W   \ar[r]^{\delta_{V} \otimes \Id_W \quad }   \ar[d]^{=}
& V \otimes V^* \otimes W  \ar[d]^{m \otimes \Id_W} \\
W  \ar[r]^{\on{Tr}_{V} \otimes \Id_W}
& \mc R \otimes W
}
\end{equation}
where $\wh G$ acts on $V \otimes V^* \otimes W$ diagonally, the action of $\wh G$ on $\mc R \otimes W$ is defined in  \ref{subsection-def-R-otimes-W-singleton}.

$(b)$ commutes because the partial Frobenius morphisms are functorial on $W$ (see \ref{subsection-partial-Frob-functorial}).

$(c)$ commutes because the following diagram of $\wh G$-representations is commutative:
\begin{equation}
\xymatrixrowsep{2pc}
\xymatrixcolsep{5pc}
\xymatrix{
V \otimes V^* \otimes W  \ar[d]^{m \otimes \Id_W}  \ar[r]^{\quad \on{ev}_{V} \otimes \Id_W  }  
& {\bf 1} \otimes W    \ar[d]^{=} \\
\mc R \otimes W  \ar[r]^{ \on{ev}_{\mc R} \otimes \Id_W}
& W  
}
\end{equation}

The composition of the left vertical line of (\ref{equation-Theta-Tr-V-equal-S-V-singleton}) is $\Theta(\on{Tr}_{V} \otimes \Id_W)$. The composition of the right vertical line of (\ref{equation-Theta-Tr-V-equal-S-V-singleton}) is $S_{V, v}$. We deduce from the commutativity of the diagram that $\Theta(\on{Tr}_{V} \otimes \Id_W) = S_{V, v}.$
\cqfd

\quad

\sssec{}   \label{subsection-Theta-u-equal-Frob-0-singleton}
Consider the morphism $$\on{act}: W \rightarrow \mc R \otimes W, \quad x \mapsto [g \mapsto g x]$$
It is $\wh G$-equivariant, for the diagonal action of $\wh G$ on $\mc R \otimes W$ defined in \ref{subsection-def-R-otimes-W-singleton}. Thus $\on{act} \in \Hom_{ \wh G }(W, \mc R \otimes W) $.

\begin{lem}   \label{lem-Theta-act-equal-Frob-0-singleton}
We have
\begin{equation}   \label{equation-Theta-act-equal-Frob-0-singleton}
\Theta( \on{act} ) = F^{\on{deg}(v)}
\end{equation}
where $F^{\on{deg}(v)}$ is the total Frobenius on $\mc H_{\{0\}, W}$.
\end{lem}
\dem
The image of act is fixed by the diagonal action of $\wh G \subset \wh G^{\{l, 0\}}$. Thus in the definition of $\Theta(\on{act})$ (c.f. ( \ref{equation-Theta-u-I-singleton})), we can replace $F_{\{r\}}^{\deg(v)}$ by $F_{\{r, l, 0\}}^{\deg(v)}$. 
The following diagram is commutative:
\begin{equation}   \label{equation-Theta-act-commute-Frob-singleton}
\xymatrixrowsep{2pc}
\xymatrixcolsep{4pc}
\xymatrix{
\restr{\mc H_{\{0\}, W}}{ v}   \ar[r]^{ \on{act} }     \ar[d]^{F^{\on{deg}(v)}}
& \restr{\mc H_{\{0\}, \mc R \otimes W}}{v}  \ar[r]_{\simeq \quad}^{\chi_{\zeta}^{-1} \quad }  \ar[d]^{F^{\on{deg}(v)}}
& \restr{\mc H_{\{r, l, 0\}, \mc R \boxtimes W}}{\Delta(v) }  \ar[d]^{F_{\{r, l, 0\}}^{\on{deg}(v)}} \\
\restr{\mc H_{\{0\}, W}}{ v}    \ar[r]^{\on{act}}  
& \restr{\mc H_{\{0\}, \mc R \otimes W}}{v}   \ar[r]_{\simeq \quad}^{\chi_{\zeta}^{-1} \quad }
& \restr{\mc H_{\{r, l, 0\} , \mc R \boxtimes W}}{\Delta(v)}   
}
\end{equation}
We deduce that 
$$
\begin{aligned}
\Theta(\on{act}) & \overset{\on{def}}{=} (\on{ev}_{\mc R} \otimes \Id_W) \circ \chi_{\zeta} \circ F_{\{r, l, 0\}}^{\on{deg}(v)}  \circ \chi_{\zeta}^{-1}  \circ \on{act} \\
& \overset{(\ref{equation-Theta-act-commute-Frob-singleton})}{=}(\on{ev}_{\mc R} \otimes \Id_W) \circ \chi_{\zeta} \circ \chi_{\zeta}^{-1}  \circ \on{act} \circ F^{\on{deg}(v)} 
\end{aligned}
$$
By definition, the composition
$$W \xrightarrow{\on{act}}  \mc R \otimes W \xrightarrow{\on{ev}_{\mc R} \otimes \Id_W}  W$$
is identity. By functoriality, $\mc H_{\{0\}, W} \xrightarrow{\on{act}} \mc H_{\{0\}, \mc R \otimes W} \xrightarrow{\on{ev}_{\mc R} \otimes \Id_W} \mc H_{\{0\}, W} $ is identity. We deduce that $\Theta(\on{act}) =  F^{\on{deg}(v)}$.
\cqfd

\begin{rem}
When $W$ is a free $\mc O_E$-module, we have a more explicit proof of Lemma \ref{lem-Theta-act-equal-Frob-0-singleton} (but this argument only works when $W$ is free). Let $W^*$ be its dual representation, the following diagram is commutative
\begin{equation}  
\xymatrixrowsep{2pc}
\xymatrixcolsep{5pc}
\xymatrix{
\restr{\mc H_{\{0\}, W}}{v }  \ar[d]^{\on{act} \quad}
 \ar@{}[dr]|{(a)}  
& \restr{\mc H_{\{0\}, W}}{v}  \ar[d]^{\Id_W \otimes \delta_W \quad \quad }   \ar[l]^{=} \\
\restr{\mc H_{\{0\}, \mc R \otimes W}}{v }    \ar[d]_{\simeq }^{\chi_{\xi}^{-1} }
& \restr{\mc H_{\{0\}, W \otimes W^* \otimes W}}{v }  \ar[d]_{\simeq  }^{\chi_{\xi}^{-1}  }   \ar[l]^{m \otimes \Id_W} \\
\restr{\mc H_{\{r, l, 0\} , \mc R \boxtimes W}}{\Delta^{\{r, l, 0\}}(v) }   \ar[d]^{F_{\{r\}}^{\on{deg}(v)}}   \ar@{}[dr]|{(b)}
& \restr{\mc H_{\{r, l, 0\} , W \boxtimes W^* \boxtimes W}}{\Delta^{\{r, l, 0\}}(v) }  \ar[d]^{F_{\{r\}}^{\on{deg}(v)}}   \ar[l]^{m \boxtimes \Id_W}  \\
\restr{\mc H_{ \{r, l, 0\} , \mc R \boxtimes W}}{\Delta^{\{r, l, 0\}}(v) }   \ar[d]_{\simeq}^{\chi_{\xi}}
& \restr{\mc H_{ \{r, l, 0\} , W \boxtimes W^* \boxtimes W}}{\Delta^{\{r, l, 0\}}(v) }  \ar[l]^{m \boxtimes \Id_W} \ar[d]_{ \simeq}^{\chi_{\xi}} \\
\restr{\mc H_{\{0\}, \mc R \otimes W}}{v }    
\ar[d]^{\on{ev}_{\mc R} \otimes \Id_W}   \ar@{}[dr]|{(c)}
& \restr{\mc H_{\{0\}, W \otimes W^* \otimes W} }{v }    \ar[d]^{\on{ev}_W \otimes \Id_W}   \ar[l]^{m \otimes \Id_W} \\
\restr{\mc H_{\{0\}, W}}{v }  
& \restr{\mc H_{\{0\}, W}}{v }   \ar[l]^{=} 
}
\end{equation} 
where $(a)$ commutes because the following diagram of $\wh G$-representations is commutative
\begin{equation}
\xymatrixrowsep{2pc}
\xymatrixcolsep{3pc}
\xymatrix{
W \otimes {\bf 1}   \ar[r]^{\Id_W \otimes \delta_W \quad \quad}   \ar[d]^{=}
& W \otimes W^* \otimes W  \ar[d]^{m \otimes \Id_W}
& x  \ar@{|->}[r]   \ar@{|->}[d]
& \sum_j x \otimes  e_j^* \otimes e_j  \ar@{|->}[d] \\
W  \ar[r]^{\on{act}}
& \mc R \otimes W
& x  \ar@{|->}[r]
& [g \mapsto \sum_j e_j^*(gx)  e_j = gx] 
}
\end{equation}
The diagram $(b)$ commutes because the partial Frobenius morphisms are functorial on $W$ (see \ref{subsection-partial-Frob-functorial}).
The diagram $(c)$ commutes because the following diagram of $\wh G$-representations is commutative
\begin{equation}
\xymatrixrowsep{2pc}
\xymatrixcolsep{5pc}
\xymatrix{
W \otimes W^* \otimes W  \ar[d]^{m \otimes \Id_W}  \ar[r]^{\quad \on{ev}_{W} \otimes \Id_W  }  
& {\bf 1} \otimes W    \ar[d]^{=} \\
\mc R \otimes W  \ar[r]^{ \on{ev}_{\mc R} \otimes \Id_W}
& W  
}
\end{equation}
The composition of left vertical line is $\Theta( \on{act} )$.
The composition of right vertical line is $(\on{ev}_W \otimes \Id_W) \circ \chi_{\zeta} 
\circ F_{\{r\}}^{\on{deg}(v)} \circ \chi_{\zeta}^{-1} \circ (\Id_W \otimes \delta_W)$. Here it is evident to see that this is equal to 
$(\on{ev}_W \otimes \Id_W) \circ \chi_{\zeta} 
\circ F_{\{r, l, 0\}}^{\on{deg}(v)} \circ \chi_{\zeta}^{-1} \circ (\Id_W \otimes \delta_W)$. Then we deduce that it is equal to $(\on{ev}_W \otimes \Id_W) \circ \chi_{\zeta} \circ\chi_{\zeta}^{-1} \circ (\Id_W \otimes \delta_W) \circ  F^{\on{deg}(v)}$.
By the "Zorro Lemma" , $(\on{ev}_W \otimes \Id_W) \circ (\Id_W \otimes \delta_W) = \Id$. Thus $\Theta(\on{act}) = F^{\on{deg}(v)}$.
\end{rem}

\quad

\noindent {\bf Proof of Proposition \ref{prop-Eichler-Shimura-general-M} (for $I$ singleton):} 

We have an isomorphism
\begin{equation}
\Psi: (\mc R \otimes \on{End}(W))^{\wh G} \isom \Hom_{\wh G}(W, \mc R \otimes W)
\end{equation}
which sends $[g \mapsto f_g \in \on{End}(W)]$ to $x \mapsto [g \mapsto f_g(x)]$. 
We have
$$(\mc R \otimes \on{End}(W))^{\wh G} \underset{\Psi}{\isom} \Hom_{\wh G}(W, \mc R \otimes W) \xrightarrow{\Theta} \on{End}_{D^{(-)}_c(  v , \mc O_E)}(\restr{\mc H_{\{0\}, W}}{v})$$
$$[g \mapsto g \in \on{End}(W)] \quad \quad \mapsto \quad \quad \on{act} \quad \quad \mapsto \quad \quad F^{\deg(v)}$$
$$[g \mapsto \on{Tr}_V(g)\Id_W \in \on{End}(W)] \quad \quad \mapsto \quad \quad \on{Tr}_V \otimes \Id_W \quad \quad \mapsto \quad \quad S_{V, v}$$

We have the Cayley-Hamilton theorem in $\on{End}(W)$. To see this, note that 
the morphism $\on{act}: W \rightarrow \mc R \otimes W$ is also $\wh G$-equivariant for the action of $\wh G$ on $\mc R \otimes W$ given by the right action of $\wh G$ on $\mc R$ and it is injective. 
Since $W$ is of finite type as $\mc O_E$-module, there exists $M \subset \mc R$ a sub-$\mc O_E$-module of finite type, invariant by the right action of $\wh G$, such that $\on{act}(W) \subset M \otimes W$. Moreover, since $\mc R$ is torsion free, $M$ is free.
For any $g \in \wh G$, the Cayley-Hamilton theorem (cf. \cite{bourbaki} Chapitre III. \S 8. 11) for the finite type free $\mc O_E$-module $M$ implies that 
\begin{equation}  \label{equation-Cayley-Hamilton-in-End-M-singleton}
\sum_{\alpha=0}^{\on{rk}M} (-1)^{\alpha} \on{Tr}_{\wedge^{\on{rk} M - \alpha} M}(g) g^{\alpha} =0  \quad \text{ in }   \on{End}(M)
\end{equation}
We deduce that (\ref{equation-Cayley-Hamilton-in-End-M-singleton}) holds in $\on{End}(M \otimes W)$, where $\wh G$ acts on $M \otimes W$ by its right action on $M$. Thus it holds in $\on{End}(\on{act}(W))$, so holds in $\on{End}(W)$.

Applying $\Theta \circ \Psi$ to (\ref{equation-Cayley-Hamilton-in-End-M-singleton}) and taking into account \ref{subsection-Theta-linear-compatible-with-composition-singleton}, we deduce that
\begin{equation}   \label{equation-sum-Theta-Tr-Theta-act-equal-0-singleton}
\sum_{\alpha=0}^{\on{rk}M} (-1)^{\alpha} \Theta(\on{Tr}_{\wedge^{\on{rk} M - \alpha} M} \otimes \Id_W)  \Theta(\on{act})^{\alpha} =0  
\end{equation}
in $\on{End}_{D^{(-)}_c(  v , \mc O_E)}(\restr{\mc H_{\{0\}, W}}{v}).$
By Lemma \ref{lem-Theta-Tr-V-equal-S-V-singleton} and Lemma \ref{lem-Theta-act-equal-Frob-0-singleton}, we deduce that (\ref{equation-sum-Theta-Tr-Theta-act-equal-0-singleton}) is
$$\sum_{\alpha=0}^{\on{rk}M} (-1)^{\alpha} S_{\wedge^{\on{rk} M - \alpha} M, v}  (F^{\on{deg}(v)})^{\alpha} =0 . $$
\cqfd

\quad

\begin{rem}
When $I = \{0\}$ is a singleton, with the notation in \cite{xiao-zhu} Section 6, we have $\on{Hom}_{\wh G}(V_1, \mc R \otimes V_2) = \on{Hom}_{\on{Coh}_{fr}^{\wh G}(\wh G)}(\wt V_1, \wt V_2)$. 
Construction \ref{construction-Theta-I-singleton} above is essentially $loc.cit.$ Lemma 6.2.5. The morphism $\Theta$ coincides with the morphism (6.2.1) in $loc.cit.$: $$S: \on{Hom}_{\on{Coh}_{fr}^{\wh G}(\wh G)}(\wt V_1, \wt V_2) \rightarrow \on{Mor}_{P^{\on{Corr}}(\on{Sht}^{\on{loc}})}(S(\wt V_1), S(\wt V_2)).$$  

Note that in $loc.cit.$, the authors considered correspondences supported on stacks of local shtukas $\on{Mor}_{P^{\on{Corr}}(\on{Sht}^{\on{loc}})}(S(\wt V_1), S(\wt V_2))$. In this paper, we consider morphisms of global cohomology $$\on{Hom}_{D^{(-)}_c( v , \mc O_E)}(\restr{\mc H_{\{0\}, V_1}}{v}, \restr{\mc H_{\{0\}, V_2}}{v}),$$ which by \cite{vincent} Section 6.4 comes from correspondences supported on stacks of local shtukas.
\end{rem}

\subsection{Proof of Proposition \ref{prop-Eichler-Shimura-general-M}: general case}   \label{subsection-proof-ES-general}

The proof is very similar to the case where $I$ is a singleton.

\sssec{}    \label{subsection-def-R-otimes-W}
Recall that in \ref{subsection-def-R-otimes-W-singleton}, $\mc R$ is equiped with an action of $\wh G^{\{r\}} \times \wh G^{\{l\}}$.

Let $I = \{0\} \cup \wt I$. Let $W \in \on{Rep}_{\mc O_E}(\wh G^I)$.
We equip $\mc R \boxtimes W$ with an action of $\wh G^{\{r\}} \times \wh G^{\{l\}} \times \wh G^I$, where $\wh G^{\{r\}} \times \wh G^{\{l\}}$ acts on $\mc R$ and $\wh G^I$ acts on $W$.

We equip $\mc R \otimes W$ with a diagonal action of $\wh G^I$, where 
the action of $\wh G^I = \wh G^{\{0\}} \times \wh G^{\wt I}$ on $\mc R$ is by conjugation of $\wh G^{\{0\}}$. In other words, if we view $\mc R \otimes W$ as the algebra of regular functions on $\wh G$ with coefficients in $W$, then the above action of $\wh G^I$ is 
$$
\begin{aligned}
\wh G^I \times (\mc R \otimes W) & \rightarrow  \mc R \otimes W, \\
\big( (h_i)_{i \in I}, f: \wh G \rightarrow W \big) & \mapsto  [g \mapsto (h_i)_{i \in I} f(h_0^{-1} g h_0)]
\end{aligned}
$$

\begin{construction}    \label{construction-Theta}
For any $V_1, V_2 \in \on{Rep}_{\mc O_E}(\wh G^I)$,
we construct a morphism 
\begin{equation}
\Hom_{ \wh G^I }(V_1, \mc R \otimes V_2)  \xrightarrow{\Theta} \on{Hom}_{D^{(-)}_c( v \times (X \sm N)^{\wt I}  , \mc O_E)}(\restr{\mc H_{I, V_1}}{v \times (X \sm N)^{\wt I} }, \restr{\mc H_{I, V_2}}{v \times (X \sm N)^{\wt I} })
\end{equation}
in the following way.
For any $u \in \Hom_{ \wh G^I }(V_1, \mc R \otimes V_2) $, let $\Theta(u)$ be the composition of morphisms
\begin{equation}   \label{equation-Theta-u-general-I}
\xymatrixrowsep{2pc}
\xymatrixcolsep{4pc}
\xymatrix{
\restr{\mc H_{I, V_1}}{v \times (X \sm N)^{\wt I} }   \ar[r]^{u}
& \restr{\mc H_{I, \mc R \otimes V_2}}{v \times (X \sm N)^{\wt I} }  \ar[r]_{\simeq \quad \quad}^{\chi_{\zeta}^{-1} \quad \quad }
& \restr{\mc H_{\{r, l\} \cup I, \mc R \boxtimes V_2}}{\Delta^{\{r, l, 0\}}(v) \times (X \sm N)^{\wt I}}  \ar[d]^{F_{\{r\}}^{\on{deg}(v)}} \\
\restr{\mc H_{I, V_2}}{v \times (X \sm N)^{\wt I} }   
& \restr{\mc H_{I, \mc R \otimes V_2}}{v \times (X \sm N)^{\wt I} }   \ar[l]_{\on{ev}_{\mc R} \otimes \Id_{V_2}}
& \restr{\mc H_{\{r, l\} \cup I, \mc R \boxtimes V_2}}{\Delta^{\{r, l, 0\}}(v) \times (X \sm N)^{\wt I}}   \ar[l]^{\simeq \quad \quad }_{\chi_{\zeta} \quad \quad}
}
\end{equation}
where 
\begin{itemize}
\item $\chi_{\zeta}$ is the fusion morphism (cf. \cite{vincent} Proposition 4.12) associated to 
\begin{equation*}
\zeta: \{r, l\} \cup I \rightarrow I, \quad
 r \mapsto 0, \; l \mapsto 0, \; i \mapsto i, \forall i \in I
\end{equation*}

\item $\Delta^{\{r, l, 0\}}: X \hookrightarrow X^{\{r, l, 0\}}$ is the diagonal morphism.


\end{itemize}

\end{construction}

\sssec{}   \label{subsection-Theta-linear-compatible-with-composition}
The morphism $\Theta$ is $\mc O_E$-linear. 
As in \ref{subsection-Theta-linear-compatible-with-composition-singleton}, the morphism $\Theta$ is compatible with composition. 

\sssec{}
When $V_1 = V_2 = W \in \on{Rep}_{\mc O_E}(\wh G^I)$, Construction \ref{construction-Theta} gives a homomorphism
\begin{equation}
\Hom_{ \wh G^I }(W, \mc R \otimes W)  \xrightarrow{\Theta} \on{End}_{D^{(-)}_c( (X \sm N)^{\wt I} \times v , \mc O_E)}(\restr{\mc H_{I, W}}{(X \sm N)^{\wt I} \times v})
\end{equation}

\sssec{}
Let $V \in \on{Rep}_{\mc O_E}(\wh G)^{\on{free}}$. Similar to \ref{subsection-def-Tr-V}, we have $\on{Tr}_{V} \otimes \Id_W \in \Hom_{ \wh G^I }(W, \mc R \otimes W) $. 

\begin{lem}    \label{lem-Theta-Tr-V-equal-S-V}
We have 
\begin{equation}    \label{equation-Theta-Tr-V-equal-S-V}
\Theta(\on{Tr}_{V} \otimes \Id_W) = S_{V, v}.
\end{equation}
where $S_{V, v}$ is defined in \ref{subsection-def-S-V-v}, here we restrict it from $\mc H_{I, W}$ (which is over $(X \sm N)^I$) to $\restr{\mc H_{I, W}}{v \times (X \sm N)^{\wt I}}$.
\end{lem}
\dem
The proof is similar to the proof of Lemma \ref{lem-Theta-Tr-V-equal-S-V-singleton}. We use the fact that the following diagrams of $\wh G^I$ representations is commutative:
\begin{equation}
\xymatrixrowsep{2pc}
\xymatrixcolsep{5pc}
\xymatrix{
{\bf 1} \otimes W   \ar[r]^{\delta_{V} \otimes \Id_W \quad }   \ar[d]^{=}
& V \otimes V^* \otimes W  \ar[d]^{m \otimes \Id_W} \\
W  \ar[r]^{\on{Tr}_{V} \otimes \Id_W}
& \mc R \otimes W
}
\end{equation}
\begin{equation}
\xymatrixrowsep{2pc}
\xymatrixcolsep{5pc}
\xymatrix{
V \otimes V^* \otimes W  \ar[d]^{m \otimes \Id_W}  \ar[r]^{\quad \on{ev}_{V} \otimes \Id_W  }  
& {\bf 1} \otimes W    \ar[d]^{=} \\
\mc R \otimes W  \ar[r]^{ \on{ev}_{\mc R} \otimes \Id_W}
& W  
}
\end{equation}
where $\wh G^I$ acts on $V \otimes V^* \otimes W$ by $\wh G^{\{0\}}$ acting diagonally on $V \otimes V^*$ and $\wh G^I$ acting on $W$.
\cqfd

\sssec{}   \label{subsection-Theta-u-equal-Frob-0}
Consider the morphism $$\on{act}: W \rightarrow \mc R \otimes W, \quad x \mapsto [g \mapsto g \cdot^0 x]$$
where $\cdot^0$ is the action of $\wh G^{\{0\}}$ on $W$ via the inclusion $$\wh G^{\{0\}} \hookrightarrow \wh G^{\{0\} \cup \wt I} = \wh G^I, \quad g \mapsto (g, 1, \cdots, 1)$$

The morphism $\on{act}$ is $\wh G^I$-equivariant, where the action of $\wh G^I$ on $\mc R \otimes W$ is defined in \ref{subsection-def-R-otimes-W}. Thus $\on{act} \in \Hom_{ \wh G^I }(W, \mc R \otimes W) $.

\begin{lem}   \label{lem-Theta-act-equal-Frob-0}
We have
\begin{equation}   \label{equation-Theta-act-equal-Frob-0}
\Theta( \on{act} ) = F_{\{0\}}^{\on{deg}(v)}
\end{equation}
where $F_{\{0\}}^{\on{deg}(v)}$ is the partial Frobenius morphism.
\end{lem}
\dem
The image of act is fixed by the diagonal action of $\wh G \subset \wh G^{\{l, 0\}}$. Thus in the definition of $\Theta(\on{act})$ (c.f. (\ref{equation-Theta-u-general-I})), we can replace $F_{\{r\}}^{\deg v}$ by $F_{\{r, l, 0\}}^{\deg v}$. 
The following diagram is commutative:
\begin{equation}   \label{equation-Frob-d-equal-Frob-0}
\xymatrix{
\restr{\mc H_{I, W}}{ v \times (X \sm N)^{\wt I} }   \ar[r]^{\on{act} \quad }  \ar[d]^{F_{\{0\}}^{\on{deg}(v)}}
& \restr{\mc H_{I, \mc R \otimes W}}{v \times (X \sm N)^{\wt I}}  \ar[r]_{\simeq \quad \quad }^{\chi_{\zeta}^{-1} \quad \quad }   \ar[d]^{ F_{\{ 0\}}^{\deg v}  }
& \restr{\mc H_{\{r, l\} \cup I, \mc R \boxtimes W}}{\Delta^{\{r, l, 0\}}(v) \times (X \sm N)^{\wt I}}  \ar[d]^{F_{\{r, l, 0\}}^{\on{deg}(v)}} \\
\restr{\mc H_{I, W}}{v \times (X \sm N)^{\wt I}}   \ar[r]^{\on{act} \quad }
& \restr{\mc H_{I, \mc R \otimes W}}{v \times (X \sm N)^{\wt I}}   \ar[r]_{\simeq \quad \quad }^{\chi_{\zeta}^{-1} \quad \quad }
& \restr{\mc H_{\{r, l\} \cup I, \mc R \boxtimes W}}{\Delta^{\{r, l, 0\}}(v) \times (X \sm N)^{\wt I}}   
}
\end{equation}
We deduce that $$
\begin{aligned}
\Theta(\on{act}) & \overset{\on{def}}{=} (\on{ev}_{\mc R} \otimes \Id_W) \circ \chi_{\zeta} \circ F_{\{r, l, 0\}}^{\on{deg}(v)}  \circ \chi_{\zeta}^{-1}  \circ \on{act} \\
& \overset{(\ref{equation-Frob-d-equal-Frob-0})}{=}(\on{ev}_{\mc R} \otimes \Id_W) \circ \chi_{\zeta} \circ \chi_{\zeta}^{-1}  \circ \on{act} \circ F_{\{0\}}^{\on{deg}(v)} 
\end{aligned}
$$
By definition, the composition
$$W \xrightarrow{\on{act}}  \mc R \otimes W \xrightarrow{\on{ev}_{\mc R} \otimes \Id_W}  W$$
is identity. We deduce (\ref{equation-Theta-act-equal-Frob-0}).
\cqfd

\quad

%
%

\noindent {\bf Proof of Proposition \ref{prop-Eichler-Shimura-general-M}:} 
The proof is very similar to the case where $I$ is a singleton. 
We have an isomorphism
\begin{equation}  
\Psi: (\mc R \otimes \on{End}(W))^{\wh G^I} \isom \Hom_{\wh G^I}(W, \mc R \otimes W)
\end{equation}
which sends $[g \mapsto f_g \in \on{End}(W)]$ to $x \mapsto [g \mapsto f_g(x)]$. In particular, $\Psi$ sends $[g \mapsto \on{Tr}_V(g) \in \on{End}(W)]$ to $\on{Tr}_V(g) \otimes \Id_W$ and sends $[g \mapsto g \cdot^0 \in \on{End}(W)]$ to $\on{act}$, 
where $\cdot^0$ is the action of $\wh G^{\{0\}}$ on $W$ via the inclusion $\wh G^{\{0\}} \hookrightarrow \wh G^{\{0\} \cup \wt I} = \wh G^I, \quad g \mapsto (g, 1, \cdots, 1)$.


The morphism $\on{act}: W \rightarrow \mc R \otimes W$ is $\wh G$-equivariant and injective, where $\wh G$ acts on $W$ by $\cdot^0$ and acts on $\mc R \otimes W$ by the right action on $\mc R$. 
As the case where $I$ is a singleton, there exists a free $\mc O_E$-module $M$ of finite type, 
such that $\on{act}(W) \subset M \otimes W$ (where $\wh G$ acts on $M \otimes W$ by its right action on $M$.). 
For any $g \in \wh G$, the Cayley-Hamilton theorem (cf. \cite{bourbaki} Chapitre III. \S 8. 11) for the finite type free $\mc O_E$-module $M$ implies that 
\begin{equation}  \label{equation-Cayley-Hamilton-in-End-M}
\sum_{\alpha=0}^{\on{rk}M} (-1)^{\alpha} \on{Tr}_{\wedge^{\on{rk} M - \alpha} M}(g) g^{\alpha} =0 
\end{equation}
in $\on{End}(M)$, thus in $\on{End}(M \otimes W)$, where $\wh G$ acts on $M \otimes W$ by its right action on $M$. We deduce that (\ref{equation-Cayley-Hamilton-in-End-M}) holds in $\on{End}(\on{act}(W))$ via $\on{act}(W) \subset M \otimes W$. Then it holds in $\on{End}(W)$ because of the injectivity of act.


Applying $\Theta \circ \Psi$ to (\ref{equation-Cayley-Hamilton-in-End-M}) and taking into account \ref{subsection-Theta-linear-compatible-with-composition}, 
Lemma \ref{lem-Theta-Tr-V-equal-S-V} and Lemma \ref{lem-Theta-act-equal-Frob-0}, we deduce 
$$\sum_{\alpha=0}^{\on{rk}M} (-1)^{\alpha} S_{\wedge^{\on{rk} M - \alpha} M, v}  (F_{\{0\}}^{\on{deg}(v)})^{\alpha} =0 $$
in $\on{End}_{D^{(-)}_c( (X \sm N)^{\wt I} \times v , \mc O_E)}(\restr{\mc H_{I, W}}{(X \sm N)^{\wt I} \times v}).$
\cqfd

%

\quad

\section{Excursion operators}

In Section \ref{section-partial-Frob-ES-integral-coef}, we proved the Eichler-Shimura relations. This will be used in \ref{subsection-specialization-morphism} to prove that the specialization morphism
$$
\mf{sp}^*: \restr{ \mc H_{G, N, I, W}^{\Oe} }{\Delta(\ov{\eta})} \rightarrow  \restr{ \mc H_{G, N, I, W}^{\Oe}  }{\ov{\eta^I}} 
$$
is a bijection. Then we use a variant of Drinfeld's lemma to equip $\restr{ \mc H_{G, N, I, W}^{\Oe}  }{\ov{\eta^I}}$ with an action of $\on{Weil}(\eta, \ov{\eta})^I$ in \ref{subsection-Drinfeld-lemma} and construct the excursion operators in \ref{subsection-excursion-operators}.

\subsection{Specialization morphism}    \label{subsection-specialization-morphism}

\sssec{}   \label{subsection-fix-specialization-from-eta-I-bar-to-delta-eta-bar}
As in \cite{vincent} Section 8 or \cite{coho-filt} 3.1.1, we denote by $\Delta: X \rightarrow X^I$ the diagonal morphism. We denote by $F^I$ the function field of $X^I$ and $\eta^I = \on{Spec}(F^I)$ the generic point of $X^I$. We fix an algebraic closure $\ov{F^I}$ of $F^I$ and denote by $\ov{\eta^I}=\on{Spec}(\ov{F^I})$ the geometric point over $\eta^I$. Moreover, we fix a specialization map 
\begin{equation}  \label{equation-specialisation-ov-eta-I-to-Delta-ov-eta}
\mf{sp}:\ov{\eta^I} \rightarrow \Delta(\ov{\eta}).
\end{equation}
It induces the homomorphism of specialization:
\begin{equation}    \label{equation-specialisation-H-G}
\mf{sp}^*: \restr{ \mc H_{G, N, I, W}^{\Oe} }{\Delta(\ov{\eta})} \rightarrow  \restr{ \mc H_{G, N, I, W}^{\Oe}  }{\ov{\eta^I}} 
\end{equation}
\begin{prop}   \label{prop-specialisation-is-bij}
The morphism (\ref{equation-specialisation-H-G}) is a bijection.
\end{prop}
\dem
The proof of the surjectivity is the same as \cite{coho-filt} Proposition 3.1.2.

The proof of the injectivity is the same as \cite{vincent} Proposition 8.32, except that we replace \cite{vincent} Proposition 7.1 (the Eichler-Shimura relation) by our Proposition \ref{prop-Eichler-Shimura-general-M}.
\cqfd

\subsection{Drinfeld's lemma}      \label{subsection-Drinfeld-lemma}

\sssec{}
As in \cite{vincent} Remarque 8.18 (also recalled in \cite{coho-filt} 3.2.1), we define the group $\on{FWeil}(\eta^I, \ov{\eta^I})$. 
We have a surjective morphism (depending on the choice of $\mf{sp}$)
\begin{equation}   \label{equation-FWeil-to-Weil-I}
\Psi: \on{FWeil}(\eta^I, \ov{\eta^I}) \twoheadrightarrow \on{Weil}(\eta, \ov{\eta})^I .
\end{equation}
Let $Q$ be its kernel (which does not depend on the choice of $\mf{sp}$).

\begin{lem} \label{lemma-Drinfeld-O-finite-type} (\cite{drinfeld-compact} Proposition 6.1, \cite{vincent} Lemme 8.2, recalled in \cite{coho-filt} Lemma 3.2.8)
A continuous action of $\on{FWeil}(\eta^I, \ov{\eta^I})$ on a $\mc O_E$-module of finite type factors through $\on{Weil}(\eta, \ov{\eta})^I$.
\end{lem}

\sssec{}
An action of $\on{FWeil(\eta^I, \ov{\eta^I})}$ on an $\mc O_E$-module $M$ is said to be continuous if $M$ is a union of $\Oe$-submodules of finite type which are stable under $\pi_1^{\on{geom}}(\eta^I, \ov{\eta^I})$ and on which the action of $\pi_1^{\on{geom}}(\eta^I, \ov{\eta^I})$ is continuous.

\begin{lem} (see \cite{coho-filt} Lemma 3.2.13 for $E$-coefficients version)   \label{lem-Drinfeld-partial-Frob-to-pi-I-module-of-finite-type}
Let $A$ be a finitely generated $\Oe$-algebra. 
Let $M$ be a $A$-module of finite type. Then a continuous $A$-linear action of $\on{FWeil(\eta^I, \ov{\eta^I})}$ on $M$ factors through $\on{Weil}(\eta, \ov{\eta})^I$.
\end{lem}

\dem (The proof is very similar to the proof of \cite{coho-filt} Lemma 3.2.13.)
For any maximal ideal $\mf m$ of $A$, 
for any $n \in \N$, the quotient $A / \mf m^n$ is of finite type (and of torsion) as a $\Oe$-module.
Since $M$ is an $A$-module of finite type, $M / \mf m^n M$ is an $A / \mf m^n$-module of finite type. Thus $M / \mf m^n M$ is an $\Oe$-module of finite type.

Applying Lemma \ref{lemma-Drinfeld-O-finite-type} to $M / \mf m^n M$, we deduce that the action of $Q$ on $M / \mf m^n M$ is trivial. 
Since $A$ is Noetherian, for any $q \in Q$ and $x \in M$, we have 
$$q \cdot x - x \in \bigcap_{\mf m \text{ max ideal }} \, \bigcap_{n=1}^{\infty} \mf m^n M \, \overset{(a)}\subset \bigcap_{\mf m \text{ max ideal }} \on{Ker} (M \rightarrow M_{\mf m}) \overset{(b)}= \{ 0 \} .$$ 
where $M_{\mf m}$ is the localization of $M$ on $A-\mf m$. (a) follows from \cite{matsumura} Theorem 8.9 and (b) follows from $loc. cit.$ Theorem 4.6. 
We deduce that $q \cdot x=x$. Thus the action of $Q$ on $M$ is trivial. 
\cqfd

\sssec{}
We have a continuous action of $\on{FWeil(\eta^I, \ov{\eta^I})}$ on $\restr{\mc H_{G, N, I, W}^{j, \, \Oe}}{\ov{\eta^I}} $ (depending on the choice of $\ov{\eta^I}$ and $\mf{sp}$) which combines the action of $\pi_1(\eta^I, \ov{\eta^I})$ and the action of the partial Frobenius morphisms (defined in \ref{subsection-partial-Frob-coef-integral}). 

\begin{prop}   \label{prop-FWeil-on-H-G-factors-through-Weil}  
The action of $\on{FWeil(\eta^I, \ov{\eta^I})}$ on $\restr{\mc H_{G, N, I, W}^{j, \, \Oe}}{\ov{\eta^I}} $ factors through $\on{Weil}(\eta, \ov{\eta})^I$.
\end{prop}
\dem
The action of $\on{FWeil(\eta^I, \ov{\eta^I})}$ on $\restr{\mc H_{G, N, I, W}^{j, \, \Oe}}{\ov{\eta^I}} $ commutes with the action of the Hecke algebra $\ms H_{G, u}^{\Oe}$. The Hecke algebra $\ms H_{G, u}^{\Oe}$
is finitely generated as a $\Oe$-algebra and is commutative.   
By Theorem \ref{coho-cht-Oe-is-Hecke-mod-type-fini-en-u}, $\restr{\mc H_{G, N, I, W}^{j, \, \Oe}}{\ov{\eta^I}} $ is of finite type as a $\ms H_{G, u}^{\Oe}$-module.
Applying Lemma \ref{lem-Drinfeld-partial-Frob-to-pi-I-module-of-finite-type} to $A=\ms H_{G, u}^{\Oe}$ and $M = \restr{\mc H_{G, N, I, W}^{j, \, \Oe}}{\ov{\eta^I}} $, we deduce the proposition.
\cqfd

\subsection{Excursion operators}    \label{subsection-excursion-operators}

\begin{construction}     \label{const-Weil-I-on-H-Delta-eta}  (for $\Ql$-coefficients, see \cite{coho-filt} Construction 3.4.1)
Let $(\gamma_i)_{i \in I} \in \on{Weil}(\eta, \ov{\eta})^I$.
We construct an action of $(\gamma_i)_{i \in I}$ on $\restr{ \mc H_{G, N, I, W}^{j, \, \Oe}   }{\Delta(\ov{\eta})}$ for any $j \in \Z$ as the composition of morphisms
\begin{equation}
\xymatrixrowsep{2pc}
\xymatrixcolsep{4pc}
\xymatrix{ 
    \restr{ \mc H_{G, N, I, W}^{j, \, \Oe}   }{\Delta(\ov{\eta})}  \ar[r]^{\mf{sp}^{*}}_{\sim}
&   \restr{ \mc H_{G, N, I, W}^{j, \, \Oe}  }{\ov{\eta^I}} \ar[d]^{(\gamma_i)_{i\in I}} \\  
  \restr{ \mc H_{G, N, I, W}^{j, \, \Oe}   }{\Delta(\ov{\eta})} 
&       \restr{ \mc H_{G, N, I, W}^{j, \, \Oe}  }{\ov{\eta^I}}  \ar[l]_{(\mf{sp}^{*})^{-1}}^{\sim}
    }
\end{equation}
where the isomorphism $\mf{sp}^*$ is defined in \ref{subsection-fix-specialization-from-eta-I-bar-to-delta-eta-bar} and Proposition \ref{prop-specialisation-is-bij}, the action of $\on{Weil}(\eta, \ov{\eta})^I$ on $\restr{\mc H_{G, N, I, W}^{j, \, \Oe} }{\ov{\eta^I}}$ is defined in Proposition \ref{prop-FWeil-on-H-G-factors-through-Weil}.
\end{construction}  

\begin{lem}  \label{lem-excursion-operator-independ-of-eta-I-bar-and-sp}   (for $\Ql$-coefficients, see \cite{coho-filt} Lemma 3.4.2)
For any $j \in \Z$, the action of $\on{Weil}(\eta, \ov{\eta})^I$ on $\restr{ \mc H_{G, N, I, W}^{j, \, \Oe}   }{\Delta(\ov{\eta})}$ defined in Construction \ref{const-Weil-I-on-H-Delta-eta} is independent of the choice of $\ov{\eta^I}$ and $\mf {sp}$. 
\end{lem}
\dem
The proof is the same as in $loc.cit.$ Lemma 3.4.2, except that now we let $\ms I$ be an ideal of $\ms H_{G, u}^{\Oe}$ such that $\ms H_{G, u}^{\Oe} / \ms I \cdot \ms H_{G, u}^{\Oe} $ is an $\Oe$-module of finite type. Then $\restr{ \mc H_{G, N, I, W}^{j, \, \Oe}   }{ \eta^I } / \ms I \cdot \restr{ \mc H_{G, N, I, W}^{j, \, \Oe}   }{ \eta^I }$ is an $\Oe$-constructible (and smooth) sheaf over $\eta^I$.
\cqfd

\begin{construction}     \label{construction-excursion-operator-on-cohomology}   (for $\Ql$-coefficients, see \cite{coho-filt} Construction 3.7.6)
Let $J$ be a finite set and $V \in \on{Rep}_{\Oe}(\wh G^J)$. 
Let $j \in \Z$. 

Let $I$ be a finite set and $W \in \on{Rep}_{\Oe}(\wh G^I)$.
Let $x \in W$ and $\xi \in W^*$ be invariant under the diagonal action of $\wh G$. Let $(\gamma_i)_{i \in I} \in \on{Weil}(\eta, \ov{\eta})^I $. We construct an excursion operator $S_{I, W, x, \xi, (\gamma_i)_{i \in I}} ^{\Oe}$ acting on $\restr{\mc  H_{G, N, J, V}^{j, \, \Oe}   }{\Delta^J(\ov{\eta})}$ as the composition of morphisms:
$$
\xymatrixrowsep{2pc}
\xymatrixcolsep{4pc}
\xymatrix{ 
 \restr{ \mc H_{G, N, J, V}^{j, \, \Oe}  }{\Delta^J(\ov{\eta}) }  \ar[r]^{  \mc C_{x}^{\sharp} \quad \quad } 
&    \restr{ \mc H_{G, N, J \sqcup I, V \boxtimes W }^{j, \, \Oe}   }{  \Delta^{J \sqcup I}(\ov{\eta})  }  \ar[d]^{ (\gamma_i)_{i\in I} } \\  
  \restr{ \mc H_{G, N, J, V}^{j, \, \Oe} }{\Delta^J(\ov{\eta}) }  
&  \restr{ \mc H_{G, N, J \sqcup I, V \boxtimes W }^{j, \, \Oe}   }{ \Delta^{J \sqcup I}(\ov{\eta}) }  \ar[l]_{  \mc C_{ \xi}^{\flat }   \quad \quad }
    }
$$
where $\mc C_{x}^{\sharp} $ and $ \mc C_{ \xi}^{\flat }$ are the creation and annihilation operators defined in \cite{vincent} Définition 5.1 and Définition 5.2, the action of $\on{Weil}(\eta, \ov{\eta})^I$ on $ \restr{ \mc H_{G, N, J \cup I, V \boxtimes W}^{j, \, \Oe}  }{ \Delta^{J \cup I}(\ov{\eta}) } $ is via $\on{Weil}(\eta, \ov{\eta})^I \hookrightarrow \on{Weil}(\eta, \ov{\eta})^{J \sqcup I}, (\gamma_i)_{i \in I}  \mapsto ((1)_{i \in J}, (\gamma_i)_{i \in I} )$.
\end{construction}

\sssec{}
In particular, when $J=\emptyset$ and $V = \bf 1$ is the trivial representation, Construction \ref{construction-excursion-operator-on-cohomology} gives the excursion operators acting on $C_c(G(F) \backslash G(\mb A) / K_N\Xi, \Oe)$.

\begin{lem}
The following diagram commutes:
$$
\xymatrixrowsep{2pc}
\xymatrixcolsep{6pc}
\xymatrix{
\restr{ \mc H_{G, N, J, V \otimes_{\Oe} E }^{j, \, \Ql}  }{\Delta^J(\ov{\eta}) }   \ar[r]^{ S_{I, W, x, \xi, (\gamma_i)_{i \in I}}  }
& \restr{ \mc H_{G, N, J, V\otimes_{\Oe} E }^{j, \, \Ql}  }{\Delta^J(\ov{\eta}) }    \\
\restr{ \mc H_{G, N, J, V}^{j, \, \Oe}  }{\Delta^J(\ov{\eta}) }    \ar[r]^{  S_{I, W, x, \xi, (\gamma_i)_{i \in I}}^{\Oe}   }   \ar[u]
& \restr{ \mc H_{G, N, J, V}^{j, \, \Oe}  }{\Delta^J(\ov{\eta}) }   \ar[u]
}
$$
where the upper line is defined in \cite{coho-filt} Construction 3.7.6.
\end{lem}

\subsection{Compatibility of the excursion operators and the constant term morphisms}

\sssec{}
Let $P$ be a parabolic subgroup of $G$ and $M$ its Levi quotient. Just as in Construction \ref{construction-excursion-operator-on-cohomology}, for any $\nu \in \wh \Lambda_{Z_M / Z_G}^{\Q}$, we define excursion operators acting on $\restr{\mc  H_{M, N, J, V}^{' \, j, \, \nu, \, \Oe}  }{\Delta^J(\ov{\eta})}$. 

\begin{prop} 
The following diagram is commutative:
$$
\xymatrixrowsep{2pc}
\xymatrixcolsep{6pc}
\xymatrix{
\restr{ \mc H_{G, N, J, V}^{j, \, \Oe}  }{\Delta^J(\ov{\eta}) }  \ar[r]^{\quad S_{I, W, x, \xi, (\gamma_i)_{i \in I}}^{G, \, \Oe} \quad}  \ar[d]^{C_{G}^{P, \, \nu, \, j}}
& \restr{ \mc H_{G, N, J, V}^{j, \, \Oe}  }{\Delta^J(\ov{\eta}) }   \ar[d]^{C_{G}^{P, \, \nu, \, j}} \\
\restr{\mc  H_{M, N, J, V}^{' \, j, \, \nu, \, \Oe}  }{\Delta^J(\ov{\eta})} \ar[r]^{\quad S_{I, W, x, \xi, (\gamma_i)_{i \in I}}^{M, \, \Oe} \quad}
& \restr{\mc  H_{M, N, J, V}^{' \, j, \, \nu, \, \Oe}  }{\Delta^J(\ov{\eta})}
}
$$
where the constant term morphism $C_{G}^{P, \, \nu, \, j}$ is defined in Remark \ref{rem-CT-cohomology-nu-Oe}.
\end{prop}
\dem
The proof is the same as the proof of \cite{coho-filt} Proposition 4.1.3.
\cqfd

\quad

\end{document}